\documentclass[3p]{elsarticle}

\usepackage[utf8]{inputenc}
\usepackage{amsmath}
\usepackage{graphicx}
\usepackage{xcolor}
\usepackage{soul}
\usepackage{bigints}
\usepackage{physics}
\usepackage{bm}
\usepackage[section]{placeins}
\usepackage{comment}
\usepackage{acronym}
\usepackage{float}
\usepackage{subcaption}
\usepackage{array}
\usepackage{amssymb}
\usepackage{braket}
\usepackage{hyperref}
\hypersetup{
    colorlinks=true,
    linkcolor=blue,
    filecolor=magenta,      
    urlcolor=cyan,
    pdftitle={Overleaf Example},
    pdfpagemode=FullScreen,
    }

\usepackage{xargs}  
\usepackage[disable,colorinlistoftodos,prependcaption,textsize=tiny]{todonotes}
\newcommand{\longtime}{long-time }
\renewcommand{\norm}[1]{\left\lVert#1\right\rVert}
\newcommand{\brac}[1]{\left(#1\right)}
\newcommand{\h}{\mathbf{H}}

\newcommand{\D}{\Delta}

\newcommand{\A}{\mathbf{A}}
\newcommand{\B}{\mathbf{B}}
\newcommand{\G}{\mathbf{G}}
\newcommand{\Gk}{\G_k}
\newcommand{\w}{\mathbf{w}}
\newcommand{\uVec}{\mathbf{u}}
\newcommand{\vVec}{\mathbf{v}}

\newcommand{\gammab}{\boldsymbol{\gamma}}
\newcommand{\gammabk}{\gammab_k}
\newcommand{\xib}{\boldsymbol{\xi}}
\newcommand{\nonnormal}{non-normal }
\newcommand{\timestep}{time-step }

\renewcommand{\tr}{\textcolor{red}}

\newcommand{\Newmark}{Newmark-$\beta$ }
\newcommand{\HHTA}{HHT-$\alpha$ }

\newcommand{\gNa}{\boldsymbol{\nabla}_{\mathbf{X}} N_A}
\newcommand{\gNb}{\boldsymbol{\nabla}_{\mathbf{X}} N_B}
\newcommand{\Hu}{\mathbf{H(u)}}
\newcommand{\Ibar}{\bar{I}_1}
\newcommand{\trC}{\text{tr}(\mathbf{C})}
\newcommand{\Jtt}{J^{-\frac{2}{3}}}

\newcommand{\dudx}[2]{\frac{\partial{u_{#1}}}{\partial{X_{#2}}}}
\newcommand{\dudxm}[4]{\frac{1}{3}\,\frac{\partial{u_{#1}}}{\partial{X_{#2}}}\,\frac{\partial{u_{#3}}}{\partial{X_{#4}}}}
\newcommand{\dudxp}[4]{\frac{1}{2}\,\left(\frac{\partial{u_{#1}}}{\partial{X_{#2}}}\,+ \,\frac{\partial{u_{#3}}}{\partial{X_{#4}}}\right)}
\renewcommand{\Finv}{\mathbf{F}^{-1}}
\newcommand{\OmegaB}{\boldsymbol{\Omega}}
\newcommand{\OmegaNStep}{\tilde{\OmegaB}(t_n,t_{n-1})}

\newcommand{\FinvT}{\mathbf{F}^{-T}}

\newcolumntype{P}[1]{>{\centering\arraybackslash}p{#1}}
\newcolumntype{M}[1]{>{\centering\arraybackslash}m{#1}}

 \newcommandx{\error}[2][1=]{\todo[linecolor=red,backgroundcolor=red!25,bordercolor=red,#1]{#2}}
 \newcommandx{\info}[2][1=]{\todo[linecolor=blue,backgroundcolor=blue!25,bordercolor=blue,#1]{#2}}
 \newcommandx{\change}[2][1=]{\todo[linecolor=orange,backgroundcolor=orange!25,bordercolor=orange,#1]{#2}}
 \newcommandx{\improv}[2][1=]{\todo[linecolor=violet,backgroundcolor=violet!25,bordercolor=violet,#1]{#2}}

\journal{Journal of the Mechanics and Physics of Solids}
\begin{document}

\begin{frontmatter}



\title{Exponential time propagators for elastodynamics}


\author[inst1]{Paavai Pari}

\affiliation[inst1]{organization={Department of Mechanical Engineering},
            addressline={University of Michigan}, 
            city={Ann Arbor},
            postcode={48105}, 
            state={Michigan},
            country={USA}}

\author[inst1]{Bikash Kanungo}
\author[inst1,inst2]{Vikram Gavini}

\affiliation[inst2]{organization={Department of Materials Science and Engineering},
            addressline={University of Michigan}, 
            city={Ann Arbor},
            postcode={48105}, 
            state={Michigan},
            country={USA}}

\begin{abstract}
We propose a computationally efficient and systematically convergent approach for elastodynamics simulations. We recast the second-order dynamical equation of elastodynamics into an equivalent first-order system of coupled equations, so as to express the solution in the form of a Magnus expansion. With any spatial discretization, it entails computing the exponential of a matrix acting upon a vector. We employ an adaptive Krylov subspace approach to inexpensively and and accurately evaluate the action of the exponential matrix on a vector. In particular, we use an \textit{apriori} error estimate to predict the optimal Kyrlov subspace size required for each time-step size.  We show that the Magnus expansion truncated after its first term provides quadratic and superquadratic convergence in the time-step for nonlinear and linear elastodynamics, respectively. We demonstrate the accuracy and efficiency of the proposed method for one linear (linear cantilever beam) and three nonlinear (nonlinear cantilever beam, soft tissue elastomer, and hyperelastic rubber) benchmark systems. For a desired accuracy in energy, displacement, and velocity, our method allows for $10-100\times$ larger time-steps than conventional time-marching schemes such as Newmark-$\beta$ method. Computationally, it translates to a $\sim$$1000\times$ and $\sim$$10-100\times$ speed-up over conventional time-marching schemes for linear and nonlinear elastodynamics, respectively.
\end{abstract}


\begin{highlights}
\item We deomonstrate that exponential time propagators are a promising technique for fast and accurate \longtime elastdynamic simulations.  
\item Exponential time propagators can be constructed for nonlinear dynamical systems using truncated Magnus expansion.
\item Within a time iteration, the exponential propagator can be accurately evaluated using a Krylov subspace technique. The subspace size required is much smaller to the system size for a sufficiently accurate propagation.
\end{highlights}


\end{frontmatter}

\section{Introduction} \label{sec:introduction}
Elastodynamics plays a crucial role in a variety of fields, such as in structural engineering to perform fatigue and fracture analysis, seismology for predicting the dynamic response of structures under seismic loads, nuclear sciences for capturing material behavior under dynamic thermal loads, material sciences for characterization and development of new materials, and aerospace engineering for material optimization. Accordingly, time integration of linear and nonlinear dynamical systems in solid mechanics is one of the widely studied fields. Much research effort has been devoted over the years to develop stable and efficient time integration schemes; see, for example~\cite{newmark1959method,bathe2006finite,simo1992exact,chung1993time,bathe2005composite,hilber1977improved,wilson1972nonlinear}. Among the various time-integration schemes, direct integration techniques remain the most popular. Here, the time derivatives in the governing equations are approximated using finite difference and the equations of motion are satisfied at discrete time intervals. Within each interval, a certain type of variation is assumed based on the Taylor series approximation for the displacement, velocity, and acceleration fields. Depending on the type of variation assumed for each of the fields, several numerical integration schemes are available and are broadly classified into explicit and implicit schemes~\cite{dukkipati2009matlab}. In an explicit scheme, such as the central difference method or the explicit Runge-Kutta schemes, the current state is approximated using previously determined values. They are easier to implement and have low computational cost given that they avoid any linear system of equations at each time-step. However, they are severely restricted in their time-step according to the stability criterion~\cite{wriggers2008nonlinear}. In contrast, implicit schemes are stable at larger time-steps, but incur higher computational cost per~\timestep. This is because, in implicit schemes, the current state depends on the solutions of the past and future states of the system. This entails performing a fixed-point iteration within each time-step, and for each of these iterations, a linear system of equations needs to be solved. Variants of the implicit scheme, extensively used in practice, are the \Newmark method~\cite{newmark1959method}, Wilson-$\theta$ method~\cite{wilson1972nonlinear} and \HHTA method~\cite{hilber1977improved}. However, all of the aforementioned methods do not perform well in the presence of high nonlinearity and/or large deformations, and become unstable in \longtime calculations~\cite{simo1992exact, simo1992discrete, bathe2005composite}. To mitigate the instability, several classes of methods have been proposed. Some methods introduce controllable numerical dissipation to damp out the spurious high frequency modes that come from spatial discretization to improve accuracy~\cite{chung1993time,hilber1977improved}. Alternatively, methods that are both energy- and momentum-conserving have also been developed to facilitate a stabler time propagation scheme. However, these methods, even with their improved characteristics, are not significantly more accurate than the simple trapezoidal rule and still require small \timestep sizes for reasonable accuracy~\cite{simo1992exact}. Hence, we are in need of a good numerical integration scheme, which exhibits superior accuracy and stability even at large \timestep sizes, operates at minimal computational cost, preserves the qualitative properties of the underlying system, and also offers good parallel scalability.
Recently, with improved numerical techniques to compute the matrix exponential~\cite{moler2003nineteen} along with advances in parallel computing, there is increasing interest in the use of exponential time propagators~\cite{hochbruck2010exponential}. Exponential propagators offer a wide array of benefits compared to their direct integration counterparts. First, exponential propagators have been shown to be a better choice for solving stiff ODEs~\cite{michels2014exponential}, which is the case in elastodynamics. Second, these methods allow for the use of larger \timestep sizes compared to state-of-the-art direct integration methods. This advantage arises from their ability to capture the dominant modes of the system more accurately, even with larger time-step sizes~\cite{michels2014exponential, loffeld2013comparative, hochbruck1998exponential}. Third, they are often designed to preserve the geometric properties of the underlying system, such as symplecticity and conservation laws, which translate into \longtime numerical stability and better capture qualitative characteristics~\cite{hairer2006structure}. For these reasons, exponential propagators have found extensive applications in many areas, including quantum chemistry and molecular dynamics~\cite{hochbruck1999exponential,nauts1983new, nettesheim1999numerical,kosloff1994propagation, chen2017exponential,michels2014exponential}.
In second-order initial value problems, such as in elastodynamics, one can still use exponential time propagators by recasting the second-order equation into a coupled system of first-order equations and constructing an exponential propagator for the resulting first-order system.
However, developing specialized exponential propagators for elastodynamics is still in its early stages, requiring detailed investigations, error analysis and assessing its computational efficiency. Recent efforts ~\cite{michels2014exponential,chen2017exponential} employing exponential propagators for elastodynamics rely on linearizing the nonlinear ODEs, either around the initial or current state. The linearization, in turn, limits the ability to take large time-steps, especially in presence of strong nonlinearity.  
In contrast, we propose a Magnus exponential propagator, wherein, after recasting the second-order ODE into a first-order system, we do not linearize the problem, but rather express the nonlinear term as an action of a nonlinear exponential operator upon the composite state-space vector. The nonlinear exponential operator is expressed in terms of a power series, namely the Magnus series, of a nonlinear Hamiltonian-like operator. To the best of our knowledge, such an approach has not been studied thus far for elastodynamics. A clear potential benefit of the Magnus exponential propagator is that it provides a systematic way to construct higher-order time propagators by inclusion of more terms in the Magnus series. 
%

In this work, we develop a second-order Magnus propagator for elastodynamics by only including the leading-order term in the Magnus series. We provide a rigorous error estimate to prove the second-order convergence of our approach as well as coroborate it numerically. 
The proposed approach mainly requires computing the action of an exponential operator on a state vector rather than evaluating the matrix exponential. To this end, we use the Krylov subspace method, as it takes advantage of this distinction and provides a significant computational benefit. We further use an \textit{a priori} error estimate to adaptively determine the optimal subspace size required for a desired accuracy at a given \timestep. 
We prove the sympleciticy for the propagator and show energy conservation in the linear model. We further demonstrate key numerical benefits of the proposed approach considering benchmark systems in linear and nonlinear hyperelastodynamics and provide performance comparisons against state-of-the-art methods such as \Newmark and HHT-$\alpha$. In linear elastodynamics, we demonstrate that the exponential propagator is more accurate in the displacement field even at timesteps that are $\mathcal{O}(10^4)$ larger than that of \Newmark method. In the velocity field, we see that the accuracies afforded by the exponential propagator are nearly intractable by the other methods. The advantage of the exponential propagator method extends to the nonlinear problem as well, providing staggering $100-400\times$ speed-up in computational time for certain benchmark problems. We also demonstrate good parallel scalability by conducting strong scaling studies. For a benchmark system with $\sim 4 \times 10^5$ degrees of freedom, the linear problem offers $\sim 74\%$ efficiency in $64$ processors, and the nonlinear problem offers $\sim 80\%$ efficiency in $128$ processors.

The rest of the paper is organized as follows. In Section \hyperref[sec:formulation]{2}, we present the formulation for elastodynamics and the associated  spatial discretization. In Section \hyperref[sec:exponential propagator]{3}, we explain in detail the construction of the Magnus-based exponential propagator. In Section \hyperref[sec:numerical implementations]{4}, we discuss the various aspects of the numerical implementation used to efficiently evaluate the propagator. In Section \hyperref[sec:analysis]{5}, we provide the analysis on accuracy and stability. In Section \hyperref[sec:numerical results]{6}, we demonstrate the efficacy of the proposed exponential propagator by studying benchmark systems in elastodynamics. In Section \hyperref[sec:conclusion and future work]{7}, we summarize the observations and discuss the future scope of this work.

\section{Formulation} \label{sec:formulation}

We provide a brief overview of the nonlinear elastodynamics equations, followed by a discussion on the various spatial and temporal techniques employed in solving these equations. 


\subsection{Governing Equations} \label{subsec:governing equations}
In this section, we discuss the governing  equations of nonlinear elastodynamics. We choose hyperelastic material models for our study modeled in the Lagrangian framework. The second Piola-Kirchoff stress tensor, $\mathbf{S}$, is related to the Green-Lagrangian strain, $\mathbf{E}$, as
\begin{align}\label{eq:Constitutive equation}
\centering
    \mathbf{S} = \frac{\partial \psi}{\partial \mathbf{E}}\text{, with } \mathbf{E} := \frac{1}{2} \left( \mathbf{F}^T\mathbf{F} - \mathbf{I}\right)\text{, } \mathbf{F} = \mathbf{I} + \grad_\mathbf{X}\bm{u}
\end{align}
where $\psi$ is the strain  energy function, $\mathbf{F}$ is the deformation gradient, $\bm{u}$ is the displacement vector, $\mathbf{I}$ is the identity tensor, and $\grad_\mathbf{X}(\,.\,)$ is the gradient operator in the reference configuration.
Using the above relations, the equations governing elastodynamics are given by
{
\newcommand{\xt}{\left(\mathbf{X},t\right)}
\newcommand{\xo}{\left(\mathbf{X},0\right)}
\newcommand{\x}{\left(\mathbf{X}\right)}
\newcommand{\gradX}{\grad_\mathbf{X}}
\newcommand{\domain}{\Omega_0}
\newcommand{\totalB}{\partial \Omega_0}
\newcommand{\DirichletB}{\partial \Omega_0^D}
\newcommand{\NeumannB}{\partial \Omega_0^N}
\newcommand{\Xindomain}{\mathbf{X} \in \domain}
\newcommand{\XinDirichlet}{\mathbf{X} \in \DirichletB}
\newcommand{\XinNeumann}{\mathbf{X} \in \NeumannB}
\newcommand{\tinspan}{t \in [0, T]}
\renewcommand{\u}{\bm{u}}
\renewcommand{\v}{\bm{v}}
\begin{alignat}{2}
        \label{eq: strong form}
        \gradX \cdot \left( \mathbf{F}\xt \mathbf{S} \xt\right) + \rho_0\bm{b} & = \rho_0 \bm{\ddot{u}}\xt &&\quad \forall\text{ } \Xindomain, \tinspan \\
       \u\xt & = \u_D\xt &&\quad \forall\text{ } \XinDirichlet, \tinspan  \\
       (\mathbf{F}\xt\mathbf{S}\xt) \cdot  \bm{n}_0 & = \bm{t}_0  &&\quad \forall\text{ } \XinNeumann, \tinspan \\
       \u\xo &= \u_0\x && \quad \forall \text{ } \Xindomain \\
       \v\xo &= \v_0\x && \quad \forall \text{ } \Xindomain
\end{alignat}
\renewcommand{\x}{\mathbf{X}}
where $\domain$ represents the domain in the reference configuration; $\x$ represents an arbitrary point in the domain $\domain$; $\totalB$ is the boundary of domain $\domain$; $\DirichletB$ and $\NeumannB$ are the Dirichlet and Neumann boundaries, respectively, such that $\totalB = \DirichletB \cup \NeumannB, \DirichletB \cap \NeumannB = \emptyset$; $\bm{n}_0$ is the unit outward normal on the boundary $\totalB$; $t$ is the time variable and $T$ is the final time; $\rho_0$ is the density in the reference configuration; $\rho_0 \bm{b}$ is the body force in the reference configuration; $\u_D$ is the prescribed displacement on the Dirichlet boundary $\DirichletB$; $\bm{t}_0$ is the prescribed traction in the Neumann boundary $\NeumannB$; $\u, \v, \ddot{\u}$ are the displacement, velocity and acceleration fields, respectively; $\u_0, \v_0$ are the initial displacement and velocity fields, respectively.

}

In this work, we employ the St.~Venant-Kirchhoff model (purely geometric nonlinearity) and the Yeoh model (material nonlinearity) to demonstrate the efficacy of the proposed approach for nonlinear models. We note that the choice of model does not affect the general applicability of the proposed method or the analysis that will be discussed subsequently.

\subsection{Finite Element Formulation} \label{subsec:finite element formulation}
{ 
\newcommand{\xt}{\left(\mathbf{X},t\right)}
\newcommand{\xo}{\left(\mathbf{X},0\right)}
\newcommand{\x}{\mathbf{X}}
\newcommand{\gradX}{\grad_\mathbf{X}}
\newcommand{\domain}{\Omega_0}
\newcommand{\elemDomain}{\Omega_e}
\newcommand{\totalB}{\partial \Omega_0}
\newcommand{\DirichletB}{\partial \Omega_0^D}
\newcommand{\NeumannB}{\partial \Omega_0^N}
\renewcommand{\u}{\bm{u}}
\renewcommand{\v}{\bm{v}}
\renewcommand{\w}{\bm{w}}
\newcommand{\mbu}{\mathbf{u}}
\newcommand{\mbv}{\mathbf{v}}
\newcommand{\mbw}{\mathbf{w}}
\newcommand{\numElem}{n_e}
\newcommand{\elemUnion}{\bigcup_{e=1}^{\numElem}}
\newcommand{\shapeFuncSumI}{\sum\limits_{I=1}^{n}}
\newcommand{\shapeFuncSumJ}{\sum\limits_{J=1}^{n}}

In this work, we approximate the spatial domain using a finite element discretization. However, we note that the choice of the spatial discretization scheme does not affect the applicability of the exponential propagator. Multiplying equation \eqref{eq: strong form} with a test function $\w(\x)$ and integrating over the spatial domain $\domain$ gives the weak formulation of elastodynamics:
\begin{equation}\label{eq:weak form}
    \begin{split}
        \int_{\domain} \w \cdot \rho_0\ddot{\u} dV & = \int_{\domain}\w \cdot \left[ \gradX \cdot \left( \mathbf{F}\mathbf{S}\right)\right] dV + \int_{\domain}\w \cdot \left( \rho_0 \bm{b}\right) dV \\
        & = -\int_{\domain} \left( \gradX\w \right) \colon \left( \mathbf{F}\mathbf{S} \right) dV + \int_{\NeumannB} \w \cdot \bm{t}_0 dA +  \int_{\domain}\w \cdot \left( \rho_0 \bm{b}\right) dV \\
        & = -\int_{\domain} \left( \mathbf{F}^T\gradX\w \right) \colon  \mathbf{S} dV + \int_{\NeumannB} \w \cdot \bm{t}_0 dA +  \int_{\domain}\w \cdot \left( \rho_0 \bm{b}\right) dV \\
        & = -\int_{\domain} \delta\mathbf{E} \colon  \mathbf{S} dV + \int_{\NeumannB} \w \cdot \bm{t}_0 dA +  \int_{\domain}\w \cdot \left( \rho_0 \bm{b}\right) dV
    \end{split}
\end{equation}
 where $dV$, $dA$ represent the differential domain volume and area, respectively; $\delta \mathbf{E}$ represents the directional derivative tensor of the Green-Lagrangian strain~($\mathbf{E}$) along the direction of $\w$; `$\colon$' represents the tensor contraction operation defined as $A \colon B = \text{trace}(A^TB)$.
 
 In the finite element framework, the spatial domain $\domain$ is discretized into $n_e$ non-overlapping subdomains ($\elemDomain$)---finite elements.
%
%
Within each finite element~($\elemDomain$), the kinematic variables as well as the geometry are interpolated isoparametrically using a suitable finite-element basis. The interpolation functions forming the basis are referred to as shape functions. 
We refer to~\cite{bathe1982finite,hughes2012finite} for the various choices of finit-elements and their approximation properties. In our study, we use hexahedral finite elements with quadratic Lagrange interpolation polynomials for our shape functions. The interpolation within each finite element~($\elemDomain$) is given by
\begin{equation}\label{eq:expansion in finite-element basis}
    \begin{split}
        \u_e(\xi) \approx \u^h_e(\xi) & = \sum_{I=1}^{n}N_I(\xi) \mbu_{I,e}\,, \\
        \w_e(\xi) \approx \w^h_e(\xi) & = \sum_{I=1}^{n}N_I(\xi) \mbw_{I,e}\,, \\
    \end{split}
\end{equation}
where $\xi$ is an arbitrary point within the element $\elemDomain$; $N_I(\xi)$ represents the value of the $I^{th}$ shape function at $\xi$; $\u_e^h, \w_e^h$ are the elemental-level approximations of the fields $\u, \w$ in the finite-element space; $\mathbf{u}_{I,e}, \mathbf{w}_{I,e}$ represent the $n$ nodal  values of $\u, \w$ respectively. The various terms in equation~\eqref{eq:weak form} can be expressed in the finite-element basis as
\begin{equation}\label{eq:weak form term 1}
    \int_{\domain} \w \cdot \rho_0\ddot{\u} dV = \elemUnion \shapeFuncSumI \shapeFuncSumJ \mbw_I^T \left( \int_{\elemDomain} \rho_0 (N_I \cdot N_J)\mathbf{I}_{3\times3} d\Omega  \right) \ddot{\mbu}_J := \mbw^T\mathbf{M}\ddot{\mbu} \,;
\end{equation}

\begin{equation} \label{eq:weak form term 2}
    \int_{\domain} \delta\mathbf{E} \colon  \mathbf{S} dV  = \elemUnion \shapeFuncSumI \mbw_I^T \left( \int_{\elemDomain} \mathbf{B}_{I,e}^T \mathbf{S}_e d\Omega \right) := \mbw^T\mathbf{R}(\mbu) \,,
\end{equation}

\newcommand{\N}[1]{\frac{\partial N_I}{\partial X_{#1}}}
\begin{equation*}
    \text{where } \mathbf{B}_{I,e}  = \begin{bmatrix}
        F_{11}\N{1} & F_{21}\N{1} & F_{31}\N{1} \\ F_{12}\N{2} & F_{22}\N{2} & F_{32}\N{2} \\ F_{13}\N{3} & F_{23}\N{3} & F_{33}\N{3} \\ F_{11}\N{2} + F_{12}\N{1} & F_{21}\N{2} + F_{22}\N{1} & F_{31}\N{2} + F_{32}\N{1} \\ F_{12}\N{3} + F_{13}\N{2} & F_{22}\N{3} + F_{23}\N{2} & F_{32}\N{3} + F_{33}\N{2} \\ F_{11}\N{3} + F_{13}\N{1} & F_{21}\N{3} + F_{23}\N{1} & F_{31}\N{3} + F_{33}\N{1}  
    \end{bmatrix} 
    \end{equation*}
    \begin{equation*}
    \text{ and } \mathbf{S}_e  = \{S_{11}, S_{22}, S_{33}, S_{12}, S_{23}, S_{13}\}^T \,;
\end{equation*}
\begin{equation} \label{eq:weak form term 3}
    \int_{\NeumannB} \w \cdot \bm{t}_0 dA +  \int_{\domain}\w \cdot \left( \rho_0 \bm{b}\right) dV = \elemUnion \shapeFuncSumI \mbw_I^T \left( \int_{\elemDomain} \rho_0 \bm{b} N_I d\Omega\right) + \bigcup_{r=1}^{n_r} \shapeFuncSumI \mbw_I^T \left( \int_{\Gamma_r} N_I \bm{t}_{0}d\Gamma \right) := \mbw^T\mathbf{P}
\end{equation}
where $n_r$ represents the number of element boundaries with Neumann boundary conditions, $\Gamma_r$ is the surface of the element with the Neumann boundary, $m$ represents the number of nodes lying on the boundary; $\mathbf{M},\mathbf{P}, \mathbf{R}(\mathbf{u}), \ddot{\mathbf{u}}$ denote the mass matrix, external force vector, internal force vector, acceleration vector, respectively, after assembly to the global system. 

Substituting equations \eqref{eq:weak form term 1},\eqref{eq:weak form term 2},\eqref{eq:weak form term 3} in the weak formulation \eqref{eq:weak form}, we arrive at the following semi-discrete equation for elastodynamics:
 \begin{equation}\label{eq:weak form discrete}
       - \mathbf{R}(\mathbf{u}) + \mathbf{P} = \mathbf{M}\ddot{\mathbf{u}}\,.
\end{equation}
For evaluating the integrals in equations \eqref{eq:weak form term 2} and \eqref{eq:weak form term 3}, we use the Gauss quadrature rule, while the integral in equation \eqref{eq:weak form term 1} is evaluated using Gauss-Legendre-Lobatto quadrature rule. The quadratic Lagrange polynomials in combination with the Gauss-Legendre-Lobatto quadrature rule renders the mass matrix diagonal. This allows for transformations involving $\mathbf{M}$, discussed in the next section, to be computed trivially. 
}

\subsection{Proposed reformulation} \label{subsec:proposed reformulation for the exponential propagator} 
{
\renewcommand{\H}{\mathbf{H}(\mathbf{u})}
\renewcommand{\u}{\mathbf{u}}
\newcommand{\Ru}{\mathbf{R}(\mathbf{u})}

In the semi-discrete equation\eqref{eq:weak form discrete}, we decompose the vector $\Ru$ into a matrix-vector product of the form $\H\u$. Going forward, we will refer to the matrix $\H$ as the H-matrix. Since the decomposition is not straightforward, some decompositions for standard hyperelastic models like Yeoh, St.~Venant-Kirchhoff models have been derived in \ref{sec:appendix:H formulation} for reference. Thus, we have the following equivalent equation to solve:
\begin{equation}\label{eq:discrete dynamic equilibrium using H}
       - \mathbf{H}(\mathbf{u})\mathbf{u} + \mathbf{P} = \mathbf{M}\ddot{\mathbf{u}} \,.
    \end{equation}
Equation \eqref{eq:discrete dynamic equilibrium using H} can be reformulated as follows
\begin{equation} \label{eq:discrete dynamic equilibrium using Hbar}
    	  \ddot{\bar{\mathbf{u}}} + \bar{\mathbf{H}}(\mathbf{u}) \bar{\mathbf{u}} = \bar{\mathbf{P}}\,, 
    	\end{equation}
where
    	\begin{equation}
    	    \bar{\mathbf{u}} = \mathbf{M}^{1/2}\mathbf{u} \,,
    	\end{equation}
    	\begin{equation}
    	   \bar{\mathbf{H}}(\mathbf{u}) = \mathbf{M}^{-1/2}\mathbf{H}(\mathbf{u})\mathbf{M}^{-1/2} \,,
    	\end{equation}
    \begin{equation}
    	   \bar{\mathbf{P}} = \mathbf{M}^{-1/2}\mathbf{P}\,.
    	\end{equation}
This reformulation preseves any symmetry property that the H-matrix may possess. We note that the evaluation of $\mathbf{M}^{1/2}$ and $\mathbf{M}^{-1/2}$ are trivial owing to the diagonal structure of the mass matrix. 

Next, we recast the second order differential equation \eqref{eq:discrete dynamic equilibrium using Hbar} into an equivalent system of first order equations. Defining $\bar{\mathbf{v}} = \dot{\bar{\mathbf{u}}}$, we get

    	\begin{equation}\label{eq: 1st order reformulation}
    	 \frac{d}{dt}\begin{bmatrix}	 {\bar{\mathbf{u}}}(t) \\	 {\bar{\mathbf{v}}}(t) \end{bmatrix}  = \begin{bmatrix}
    	       0 & \mathbf{I} \\
    - \bar{\mathbf{H}}(\mathbf{u}(t)) & 0
    	     \end{bmatrix} \begin{bmatrix}
    	 \bar{\mathbf{u}}(t) \\
    	 \bar{\mathbf{v}}(t)
    	\end{bmatrix} \qquad (\text{when }\bar{\mathbf{P}} = 0)
    	\end{equation}
    	
    		\begin{equation}
    	\frac{d}{dt}\begin{bmatrix}	 \bar{\mathbf{u}}(t) \\	 \bar{\mathbf{v}}(t) \\ 1 \end{bmatrix} = \begin{bmatrix}
    	       0 & \mathbf{I} & 0 \\
    - \bar{\mathbf{H}}(\mathbf{u}(t)) & 0 & \bar{\mathbf{P}} \\
    0 & 0 & 0
    	     \end{bmatrix} \begin{bmatrix}
    	 \bar{\mathbf{u}}(t) \\
    	 \bar{\mathbf{v}}(t) \\
    	 1
    	\end{bmatrix} \qquad (\text{when }\bar{\mathbf{P}} \neq 0)
    	\end{equation}
Substituting
    \begin{subequations}\label{eq:A expression}
    \begin{align}
  \bar{\mathbf{w}}(t) = \begin{bmatrix}
    	 \bar{\mathbf{u}}(t) \\
    	 \bar{\mathbf{v}}(t)
    	\end{bmatrix} \quad \text{ and } \qquad \mathbf{A}[\mathbf{u}](t) = \begin{bmatrix}
    	 0 & \mathbf{I} \\
    - \bar{\mathbf{H}}(\mathbf{u}(t)) & 0
    	\end{bmatrix} \qquad (\text{when }\bar{\mathbf{P}} = 0) \label{eq:A expression at P=0}\\
  \bar{\mathbf{w}}(t) = \begin{bmatrix}
    	 \bar{\mathbf{u}}(t) \\
    	 \bar{\mathbf{v}}(t) \\
    	 1
    	\end{bmatrix} \quad \text{ and } \qquad \mathbf{A}[\mathbf{u}](t) = \begin{bmatrix}
    	 0 & \mathbf{I} & 0 \\
    - \bar{\mathbf{H}}(\mathbf{u}(t)) & 0 & \bar{\mathbf{P}} \\
    0 & 0 & 0
    	\end{bmatrix} \qquad (\text{when }\bar{\mathbf{P}} \neq 0)
    \label{eq:A expression at P != 0} 
     \end{align}
    \end{subequations}
we arrive at
    \begin{equation}\label{eq:1st order dynamic equation}
        \dot{\bar{\mathbf{w}}}(t) = \mathbf{A}[\mathbf{u}](t)\bar{\mathbf{w}}(t) \,.
    \end{equation}
Thus, we have reformulated equation \eqref{eq:discrete dynamic equilibrium using Hbar} into a first order ODE. 

 }

\section{Exponential Propagator} \label{sec:exponential propagator}


        \subsection{Magnus expansion} \label{subsec: Magnus expansion}

            We now construct the full-discrete solution for the elastodynamics problem by temporally discretizing the semi-discrete equation recast into a first-order initial value problem (cf.equation~\eqref{eq:1st order dynamic equation}). To this end, we construct an exponential time propagator by invoking the Magnus ansatz~\cite{MagnusReviewPaper} on the first order initial value problem. Thus, the solution of equation \eqref{eq:1st order dynamic equation} can be written as

     \begin{equation}\label{MagnusEquation}
        \bar{\mathbf{w}}(t) = \operatorname{exp}(\OmegaB(t,t_0))\bar{\mathbf{w}}(t_0)\,,
    \end{equation}
where $\bar{\mathbf{w}}(t_0)$ and $\bar{\mathbf{w}}(t)$  denote the initial and evolved states of the system, respectively; $\operatorname{exp}(\OmegaB(t,t _0))$ is the exponential propagator where $\OmegaB(t,t_0)$ is given by
    \begin{equation}\label{OmegaBtt0Expr}
\OmegaB(t,t_0) = \int_{t_{0}}^{t} \mathbf{A}(\tau)d\tau - \frac{1}{2} \int_{t_{0}}^{t} \left[ \int_{t_{0}}^{\tau} \mathbf{A}(\sigma) d\sigma , \mathbf{A}(\tau) \right] d\tau + ...
    \end{equation}
with $\left[\mathbf{X}, \mathbf{Y} \right] = \mathbf{X}\mathbf{Y} - \mathbf{Y}\mathbf{X}$ representing the commutator.
    
 
\subsection{Time discretization} \label{subsec:time discretization}
        
The propagator $\operatorname{exp}(\OmegaB(t,t_0))$ directly takes the initial state at $t_0$ to the desired state $t$ in one time propagating step. However, for practical applications requiring a reasonable accuracy, the evaluation of the propagator with a single large time-step would involve an intractable computational cost. This issue can be avoided by time discretizing the propagator and thus breaking down the propagation into smaller time-steps. To this end, we use the following  composition property of the propagator:
  \begin{equation}\label{compositionProperty}
\operatorname{exp}(\OmegaB(t_2,t_0)) = \operatorname{exp}(\OmegaB(t_2,t_1)) \operatorname{exp}(\OmegaB(t_1,t_0))\,.
    \end{equation}
This property can be exploited to achieve temporal discretization described as
 \begin{equation}\label{temporalDisc}
\operatorname{exp}(\OmegaB(t_n,t_0)) = \prod_{i=0}^{n-1} \operatorname{exp}(\OmegaB(t_{i+1},t_{i}))\,,
    \end{equation}
where $t_n = t$, $\Delta t = \frac{t_n-t_0}{n}$ and $t_{i+1} = t_{i} + \Delta t$.
    
The smaller the time-step size $\Delta t$, more accurate is the time propagation. The composition property enables us to contain the error within a time-step and the accumulated error over the entire propagation within the desired tolerance. However, the number of time-steps scale inversely with the time-step size, increasing the computational cost. Thus, from a computational efficiency standpoint, it is desirable to have a method that strikes a balance beteween the time-step size and computational cost per time-step to achieve the desired accuracy. 


\subsection{Approximations} \label{subsec:approximations}

The exponent $\OmegaB(t,t_0)$ (cf.~equation \eqref{OmegaBtt0Expr}) is evaluated as a series, and the computational complexity increases exponentially with the number of terms in the series. Thus, for all practical purposes, we truncate the series to a few terms to obtain suitable approximations. Truncating the series at the first term gives a second-order accurate time propagator. We show in Section~\ref{subsec:accuracy} that for 
\begin{equation}\label{truncOmegaB}
 \bar{\OmegaB}(t,t_0) = \int_{t_{0}}^{t} \mathbf{A}(\tau)d\tau \,,
\end{equation}
%
\begin{equation}\label{timeDiscError}
     \norm{\left( \exp{(\OmegaB(t_n,t_0))} - \exp{(\bar{\OmegaB}(t_n,t_0))}  \right) \bar{\mathbf{w}}(t_{0})} \leq C (\Delta t)^2 \,,
\end{equation}
where $C$ is problem-specific bounded positive constant. We note that truncating the series after two terms gives a fourth-order accurate scheme. However, our numerical studies showed that the fourth-order method was overall computationally more expensive than the second-order method. Thus, we restricted our study to second-order Magnus propagators.


\section{Numerical Implementations} \label{sec:numerical implementations}    
We present the various numerical aspects related to the implementation of the proposed exponential propagator method. 

       

        \subsection{Evaluating the integral in the exponential} \label{subsec:evaluating the integral expression}
            
As discussed in section \ref{subsec:approximations}, truncating the Magnus expansion after the first term leads to a second-order accurate method. Thus, in order for the second order convergence to be retained, the integral is approximated using a mid-point quadrature rule, which is also second order accurate.  
Thus, 
\begin{equation}\label{eq:quadratureApprox}
 \bar{\OmegaB}(t_n,t_{n-1}) = \int_{t_{n-1}}^{t_n} \mathbf{A}(\tau)d\tau \approx \mathbf{A}\left(t_{n-1}+\frac{\Delta t}{2}\right)(\Delta t) = \tilde{\OmegaB}(t_n,t_{n-1})\,.
\end{equation}

\subsection{Predictor-corrector to evaluate \texorpdfstring{$\mathbf{A}\left(t_{n-1}+\frac{\Delta t}{2}\right)$}{A(tn-1 + delta t/2}} \label{subsec:pred-corr}
           
Evaluation of $\mathbf{A}\left(t_{n-1}+\frac{\Delta t}{2}\right)$ requires knowledge of the H-matrix at a future time instance that is not yet known. A first-order approximation of the H-matrix is sufficient to construct a second-order exponential propagator, and this is achieved by employing a predictor-corrector scheme~\cite{PredCorr}.
The values of $\mathbf{\bar{H}}$ evaluated at the midpoints of previous time-steps are used to construct a linear extrapolation of $\mathbf{\bar{H}}\left(t_{n-1}+\frac{\Delta t}{4}\right)$, which is in turn used to propagate $\mathbf{\bar{w}}_{n-1}$ to $\mathbf{\bar{w}}_{n-\frac{1}{2}}$. Subsequently, $\mathbf{\bar{w}}_{n-\frac{1}{2}}$ is used to construct $\mathbf{\bar{H}}\left(t_{n-1}+\frac{\Delta t}{2}\right)$, which is then used to propagate $\mathbf{\bar{w}}_{n-1}$ to $\mathbf{\bar{w}}_{n}$. We note that, in this approach, the H-matrix is computed only once per time-step of propagation.

\subsection{Evaluation of \texorpdfstring{$\exp(\OmegaNStep)\w_{n-1}$}{eA}} \label{subsec:evaluation of e^A}

Exponential propagator entails the evaluation of $\exp(\A)\w$. Many methods have been proposed over the years to efficiently evaluate the exponential of a matrix. We refer to~\cite{moler2003nineteen} for a comprehensive discussion. A common classification of the methods developed, as mentioned in~\cite{moler2003nineteen}, includes series-based methods, decomposition techniques, polynomial methods, splitting methods, and Krylov subspace-based approximations. We note that evaluation of $\exp(\A)$ by matrix decomposition techniques is prohibitively expensive for large system sizes, and can be numerically intractable for non-normal matrices. The polynomial methods are computationally expensive and suffer severe round-off errors. They are also not suitable for large system sizes, and, as they are based on the evaluation of eigenvalues, they are ill suited for non-normal matrices. The series methods, although computationally less expensive than decomposition techniques and polynomial methods, are highly prone to round-off errors for large norms of $\A$ with non-normality leading to unstable solutions. Among the series-based methods, scaling and squaring is the least susceptible to round-off errors. However, it is again ill suited for large matrix sizes. The efficiency of split-order techniques relies on the ability to split the exponent into sub-parts that can be inexpensively exponentiated. It is commonly employed to solve the Schr\"{o}ndinger equation where splitting the Laplace operator and potential energy results in faster exponentiation schemes. However, such splitting operations are not trivial for exponents that arise in elastodynamics. 

We note that the evaluation of the exponential of a matrix acting on a vector is a far different problem compared to the evaluation of the exponential of a matrix.  Although the $\A$ matrix, computed in the finite-element basis is sparse, $\exp(\A)$ is a dense matrix in general. Thus, for large-scale problems, the storage and handling of $\exp(\A)$ would be infeasible and hard to parallelize. The Krylov subspace method is an approach which takes advantage of this distinction, and compared to the aforementioned methods, most suitable in evaluating $\exp(\A)\w$~\cite{hochbruck1997krylov}. The approach is computationally less expensive and is well suited for large-scale problems and is not affected by the non-normality, thus making it a desirable choice. Krylov subspace techniques also exhibit good parallel scalability, since they rely solely on matrix-vector products for their construction.

The Krylov subspace method is based on the construction of a subspace $\mathbf{V}_m$ spanned by the vectors $=\{\w,\A\w,{\A}^2\w,...,{\A}^m\w\}$, where $m$ is the size of the subspace considered. The columns of $\mathbf{V}_m$ are orthogonal, denoted $\{\mathbf{v}_1,\mathbf{v}_2,...,\mathbf{v}_m\}$ and $\mathbf{v}_1 = \w/\norm{\w}$. Such a subspace can be effectively constructed by employing algorithms like Lanczos for a symmetric $\A$ matrix and Arnoldi techniques for a non-symmetric $\A$ matrix. The algorithm generates two matrices $\mathbf{V}_m$ and $\mathbf{H}_m$, where $\mathbf{H}_m$ is a square matrix of dimension $m$, and the following relation holds: 
\begin{equation}
    \A\mathbf{V}_m = \mathbf{V}_m\mathbf{H}_m + h_{m+1,m} \mathbf{v}_{m+1} \mathbf{e}_m^T \,,
\end{equation}
where $h_{m+1,1} = \left[\mathbf{H}_{m+1}\right]_{m+1,1}$, $\mathbf{v}_{m+1}$ is the last column of the matrix $\mathbf{V}_{m+1}$ and $\mathbf{e}_m$ is the $m$-th unit vector.
From this equation, we can observe that the $\mathbf{V}_m^T\A\mathbf{V}_m = \mathbf{H}_m$ and this property is exploited in approximating $\exp(\A)\w$ as
\begin{align}
    \exp(\A)\w & \approx \norm{\w}\mathbf{V}_m\mathbf{V}_m^T\exp(\A)\mathbf{V}_m\mathbf{e}_1 \\
    & \approx \norm{\w}\mathbf{V}_m\exp(\mathbf{V}_m^T\A\mathbf{V}_m)\mathbf{e}_1 \\
    & = \norm{\w}\mathbf{V}_m\exp(\mathbf{H}_m)\mathbf{e}_1\,.
\end{align}

Using this approximation~\cite{gallopoulos1992efficient}, the evaluation of $\exp(\A)\w$ reduces to the computation of $\exp(\mathbf{H}_m)$. Typically, the subspace convergence is super-linear, and the $m$ required for the desired accuracy is much smaller compared to the system size. 
When $\norm{\mathbf{H}_m} \leq 1$, then $\exp(\mathbf{H}_m)$ can be efficiently calculated using Taylor series methods, while scaling and squaring can be used when $\norm{\mathbf{H}_m} > 1$. We do not use matrix decomposition techniques to evaluate $\exp(\mathbf{H}_m)$ as the matrix is non-normal and can result in instabilities while computing the eigendecomposition. An estimate for the error incurred by this approximation for a subspace size of $m$ is given by~\cite{hochbruck1998exponential},
       
       \begin{equation}\label{eq:subspace error estimate}
            \norm{\exp(\A)\w - \norm{\w}\mathbf{V}_m\exp(\mathbf{H}_m)\mathbf{e}_1} \approx \beta_{m+1,m}\norm{\w}|\left[ \exp(\mathbf{H}_m) \right]_{m,1}| = \epsilon_m,
       \end{equation}
where $\beta_{m+1,m}$ is the $(m+1,m)$ entry of $\mathbf{H}_{m+1} = \mathbf{V}_{m+1}^T\A\mathbf{V}_{m+1}$. This estimate provides an efficient and adaptive way to determine the size of the subspace for a desired error tolerance.
       

\section{Analysis} \label{sec:analysis}

In this section, we present the numerical analysis of the exponential propagator. In particular, we prove that the exponential propagator is second-order accurate for nonlinear elastodynamics. We also demonstrate energy conservation for linear elastodynamics, and discuss symplecticity.


\subsection{Accuracy}\label{subsec:accuracy}
         
We note that in the case of linear elastodynamics, there is no truncation error in the Magnus expansion, as the expression is exact. In the nonlinear problem, however, the truncation of the Magnus propagator to the leading order results in a second-order accurate method. In this section, we derive this error bound. To this end, we introduce the following:
\begin{align}
        S_0^{k} = e^{\OmegaB_k}e^{\OmegaB_{k-1}}...e^{\OmegaB_1} = \prod_{l=0}^{k-1} e^{\OmegaB_{k-l}} \quad  \mbox{ for } 0 < k \leq n, \qquad S_0^0 = I \label{eq:S0k}\\
          R_k^n = e^{\tilde{\OmegaB}_n}e^{\tilde{\OmegaB}_{n-1}}...e^{\tilde{\OmegaB}_{k+1}} = \prod_{l=0}^{n-k-1}e^{\tilde{\OmegaB}_{n-l}} \quad \mbox{ for } 0 \leq k < n, \qquad R_n^n = I \label{eq:Rkn}
          \end{align}
          \begin{gather}
          \text{where } \OmegaB_k = \int_{{t_{k-1}}}^{{t_{k}}} \A(\tau)d\tau - \frac{1}{2} \int_{{t_{k-1}}}^{{t_{k}}} \left[ \int_{{t_{k-1}}}^{\tau} \A(\sigma) d\sigma , \A(\tau) \right] d\tau + ... \,, \label{eq:OmegaB k}\\
    \tilde{\OmegaB}_k = \A\left(\frac{t_{k-1} + t_{k}}{2} \right)  \Delta t \,, \label{eq:OmegaB tilde k}\\
    \A(t) = \begin{bmatrix} 0 & I \\ - \bar{\mathbf{H}}(\mathbf{u}(t)) & 0 \end{bmatrix}\,,\quad    \bar{\mathbf{H}}(\mathbf{u}) = \mathbf{M}^{-1/2}\mathbf{H}(\mathbf{u})\mathbf{M}^{-1/2} \,, \label{eq:Adef} \\
     \mbox{ and }  \bar{\w}(t) = \begin{bmatrix}
    \mathbf{M}^{1/2}\uVec(t) \\
    \mathbf{M}^{1/2}\vVec(t)
    \end{bmatrix} \,.
\end{gather}
%





 
\vspace{0.25cm}
\noindent The time discretization error can be expressed as follows
\begin{equation}
    \bar{\w}(t_n) - \bar{\w}_n = \left( R_n^n S_0^n - R_0^n S_0^0 \right)\bar{\w}(0) = \sum_{k=1}^{n}\left( R_k^n S_0^k - R_{k-1}^n S_0^{k-1} \right) \bar{\w}(0) \,.
\end{equation}
This expression can be further manipulated as
 \begin{equation}
 \begin{split}
    \bar{\w}(t_n) - \bar{\w}_n  & = \sum_{k=1}^{n}\left( R_k^n S_0^k - R_{k-1}^n S_0^{k-1} \right) \bar{\w}(0) \\
    & = \sum_{k=1}^{n}\left( R_k^n e^{\OmegaB_k} S_0^{k-1} - R_k^n e^{\tilde{\OmegaB}_k} S_0^{k-1} \right) \bar{\w}(0) \\
    & = \sum_{k=1}^n R_k^n \left( e^{\OmegaB_k} - e^{\tilde{\OmegaB}_k} \right) S_0^{k-1} \bar{\w}(0) \\
    & =  \sum_{k=1}^n R_k^n \left( e^{\OmegaB_k} - e^{\tilde{\OmegaB}_k} \right) \bar{\w}(t_{k-1}) \,.
    \end{split}
\end{equation}
Henceforth, $\norm{.}$ represents the Euclidean norm for a vector and $2$-norm for a matrix. We arrive at the following bound
\begin{equation}\label{eq:boundNormBarW}
\begin{split}
    \norm{\bar{\w}(t_n) - \bar{\w}_n} & = \norm{\sum_{k=1}^n R_k^n \left( e^{\OmegaB_k} - e^{\tilde{\OmegaB}_k} \right) \bar{\w}(t_{k-1})}\\
    & \leq \sum_{k=1}^n \norm{R_k^n \left( e^{\OmegaB_k} - e^{\tilde{\OmegaB}_k} \right) \bar{\w}(t_{k-1})} \\
    & \leq \sum_{k=1}^n \norm{R_k^n} \norm{\left( e^{\OmegaB_k} - e^{\tilde{\OmegaB}_k} \right) \bar{\w}(t_{k-1})}\,.
    \end{split}
\end{equation}


\vspace{0.25cm}
 
\noindent Using $\bar{\OmegaB}_k = \int_{t_{k-1}}^{t_k} \A(\tau)d\tau$, we split $\left( e^{\OmegaB_k} - e^{\tilde{\OmegaB}_k} \right) \bar{\w}(t_{k-1})$ as follows, arriving at
 \begin{equation}
     \left( e^{\OmegaB_k} - e^{\tilde{\OmegaB}_k} \right) \bar{\w}(t_{k-1}) = \left( e^{\OmegaB_k} - e^{\bar{\OmegaB}_k} \right) \bar{\w}(t_{k-1}) + \left( e^{\bar{\OmegaB}_k} - e^{\tilde{\OmegaB}_k} \right) \bar{\w}(t_{k-1}) \,,
 \end{equation}
\begin{equation}\label{nonlinearConvError_split}
     \norm{\left( e^{\OmegaB_k} - e^{\tilde{\OmegaB}_k} \right) \bar{\w}(t_{k-1})} \leq \norm{\left( e^{\OmegaB_k} - e^{\bar{\OmegaB}_k} \right) \bar{\w}(t_{k-1})} + \norm{\left( e^{\bar{\OmegaB}_k} - e^{\tilde{\OmegaB}_k} \right) \bar{\w}(t_{k-1})}\,.
 \end{equation}
We first seek to obtain a bound for $\norm{\left( e^{\OmegaB_k} - e^{\bar{\OmegaB}_k} \right) \bar{\w}(t_{k-1})}$. To this end, we define
\begin{equation}\label{xiDef}
    \xib_k(t) = e^{\B_k(t)} \bar{\w}(t_{k-1}),\text{ where }\B_k(t) = \int_{t_{k-1}}^{t} \A(\tau)d\tau \text{ and }\B_k(t_k) = \bar{\OmegaB}_k \,.
\end{equation}
Before proceeding, we note that for an exponential operator of the form $e^{\mathbf{L}(t)}$, the derivative with respect to $t$ is given by~\cite{MagnusReviewPaper}
%
\begin{equation}
    \frac{d}{d t} e^{\mathbf{L}(t)} = \text{dexp}_{\mathbf{L}(t)}(\dot{\mathbf{L}}(t))e^{\mathbf{L}(t)}\,,
\end{equation}
where $\text{dexp}_{\mathbf{X}}(\mathbf{Y}) = \sum_{k=0}^{\infty} \frac{1}{(k+1)!} \text{ad}^k_{\mathbf{X}}(\mathbf{Y})$. The operator $\text{ad}^k_{\mathbf{X}}(\mathbf{Y})$ is defined as 
\begin{equation}
  \text{ad}^k_{\mathbf{X}}(\mathbf{Y}) = \text{ad}_{\mathbf{X}}(\text{ad}^{k-1}_{\mathbf{X}}(\mathbf{Y}))\,,
\end{equation}
with $\text{ad}^1_{\mathbf{X}}(\mathbf{Y}) = \mathbf{XY - YX}$ and $\text{ad}^0_{\mathbf{X}}(\mathbf{Y}) = \mathbf{Y}$.\\ 

\noindent Thus, upon differentiating $\xib_k(t)$, we get
\begin{equation}
    \dot{\xib}_k(t) =  \text{dexp}_{\B_k(t)}(\dot{\B}_k(t))e^{\B_k(t)} \bar{\w}(t_{k-1}) = \mathbf{G}_k(t)\xib_k(t)\,,
\end{equation}
where
\begin{equation} \label{Gdef}
\begin{split}
 \mathbf{G}_k(t) & =  \text{dexp}_{\mathbf{B}_k(t)}(\dot{\mathbf{B}}_k(t)) = \sum_{k=0}^{\infty} \frac{1}{(k+1)!} \text{ad}^{k}_{\B_k(t)}\dot{\B}_k(t)\\
   & =  \begin{bmatrix}
    0 & \mathbf{I} \\
    -\bar{\h}(t) & 0
    \end{bmatrix} + \frac{(t - t_{k-1})^2}{2}\begin{bmatrix}
    -\frac{\bar{\h}'(t)}{2} & 0 \\
    0 & \frac{\bar{\h}'(t)}{2}
    \end{bmatrix} + O((t - t_{k-1})^3) \,.
   \end{split}
\end{equation}

\noindent Next, we define $\gammab_k(t)$ as follows,
\begin{equation}\label{gammaDef}
    \gammab_k(t) = \bar{\w}(t) - \xib_k(t)\
\end{equation}
Then
\begin{equation}\label{gammaFullDef}
    \gammab_k(t_k) = \bar{\w}(t_k) - \xib_k(t_k) = \left( e^{\OmegaB_k} - e^{\bar{\OmegaB}_k} \right) \bar{\w}(t_{k-1})\,
\end{equation}

\noindent We know from Eq.~\eqref{gammaDef} and Eq.~\eqref{xiDef} that
    \begin{equation}\label{gammaIC}
        \gammabk(t_{k-1}) = 0 \,.
    \end{equation}

\noindent Upon differentiating $\gammabk(t)$, we get

\begin{equation}\label{gammaDotdef}
    \dot{\gammab}_k(t) = \Gk(t)\gammabk(t) + \left[ \A(t) - \Gk(t) \right] \bar{\w}(t) \,.
\end{equation}

\noindent Applying one iteration of the Picard scheme gives

\begin{equation}\label{1PicardIter}
    \gammabk^1(t) = \gammabk(t_{k-1}) + \int_{t_{k-1}}^t \left[ \Gk(s)\gammabk(t_{k-1}) + \left[ \A(s) -\Gk(s) \right]\bar{\w}(s) \right]ds
\end{equation}

\noindent where $\gammabk^1(t)$ is the first approximation of the Picard iteration. Substituting Eq.~\eqref{gammaIC} and Eq.~\eqref{Gdef} in Eq.~\eqref{1PicardIter}, we get

\begin{equation}\label{1PicardIterWithIC}
\begin{split}
    \gammabk^1(t) & = \int_{t_{k-1}}^t \left[ \A(s) - \Gk(s)\right]\bar{\w}(s) ds \\
    & = \int_{t_{k-1}}^t \frac{(s-t_{k-1})^2}{4} \begin{bmatrix}
    \bar{\h}'(t) & 0 \\
    0 & -\bar{\h}'(t)
    \end{bmatrix}\bar{\w}(s)ds + \int_{t_{k-1}}^t O((s - t_{k-1})^3)ds
\end{split}
\end{equation}

\noindent Thus,

\begin{equation}\label{1PicardIterNorm}
\begin{split}
    \norm{\gammabk^1(t)} & \leq \norm{\int_{t_{k-1}}^t \frac{(s-t_{k-1})^2}{4} \begin{bmatrix}
    \bar{\h}'(t) & 0 \\
    0 & -\bar{\h}'(t)
    \end{bmatrix}\bar{\w}(s)ds} + O((t - t_{k-1})^4) \\
    & \leq C\max_{\substack{t_{k-1}\leq p \leq t}}\norm{\bar{\h}'(p)}\max_{\substack{t_{k-1}\leq p \leq t}}\norm{\bar{\w}(p)}\frac{(t - t_{k-1})^3}{3!}
\end{split}
\end{equation}
where $C>0$ is a constant. Next, consider the norm of the difference between the $n^{th}$ and $(n-1)^{th}$  approximations of the Picard iteration:
\begin{equation*}\label{PicardDiffn}
\begin{split}
    \norm{\gammabk^n(t) - \gammabk^{n-1}(t)} & = \norm{\int_{t_{k-1}}^t \Gk(s) \left[ \gammabk^{n-1}(s) - \gammabk^{n-2}(s)\right] ds} \\
    & \leq \max_{\substack{t_{k-1}\leq p \leq t}}\norm{\Gk(p)}\int_{t_{k-1}}^t \norm{\gammabk^{n-1}(s) - \gammabk^{n-2}(s)} ds \\
    & \text{substituting the relation recursively, we get} \\
    & \leq \left[ \max_{\substack{t_{k-1}\leq p \leq t}}\norm{\Gk(p)} \right]^{n-1} \underbrace{\int \int ... \int}_{n-1} \norm{\gammabk^{1}(s) - \gammabk^{0}(s)} (ds)^{n-1} \\
    & \text{substituting Eq.~\eqref{1PicardIterNorm}, we get}\\
    & \leq C \max_{\substack{t_{k-1}\leq p \leq t}}\norm{\bar{\h}'(p)}\max_{\substack{t_{k-1}\leq p \leq t}}\norm{\bar{\w}(p)} (t - t_{k-1})^3\frac{\left[ \max\limits_{t_{k-1}\leq p \leq t}\norm{\Gk(p)} (t - t_{k-1}) \right]^{n-1} }{(n+2)!}
    \end{split}
\end{equation*}

\noindent With the Picard iteration, the norm of $\gammabk(t)$ converges to 
%
%
\begin{equation}\label{PicardnthNorm}
        \lim_{n \to \infty}\norm{\gammabk^n(t)} \leq C \max_{\substack{t_{k-1}\leq p \leq t}}\norm{\bar{\h}'(p)}\max_{\substack{t_{k-1}\leq p \leq t}}\norm{\bar{\w}(p)} (t - t_{k-1})^3  \exp{\max_{\substack{t_{k-1}\leq p \leq t}}\norm{\Gk(p)}(t - t_{k-1})} \,.
\end{equation}

Before proceeding, let us consider the conditions for convergence of the Magnus series expansion. For a given operator $\A(t)$, the Magnus series given by Eq.~\eqref{OmegaBtt0Expr} converges when
\begin{equation}\label{eq:MagnusConvergenceCriteria}
    \int_{t_0}^t\norm{\A(s)}ds < r_c\,,
\end{equation}
where $r_c$ is the radius of convergence. Several bounds to the actual radius of convergence have been obtained in the literature~\cite{MagnusReviewPaper, pechukas1966exponential,karasev1976infinite,blanes1998magnus,moan1998efficient}. Most commonly $r_c$ is taken to be $\ln{2}$.
 Thus, we have 
\begin{equation}\label{eq:ad(B_k(t)) bound}
    \norm{\text{ad}_{\B_k(t)}} \leq 2\norm{\B_k(t)} \leq 2\norm{\int_{t_{k-1}}^{t} \A(\tau)d\tau } \leq 2 \int_{t_{k-1}}^t \norm{\A(\tau)}d\tau \leq 2r_c
\end{equation}
Using this criterion, we bound $\norm{\Gk(t)}$ as follows
\begin{equation}\label{eq:bound G_k(t)}
\begin{split}
    \norm{\Gk(t)} & = \norm{\text{dexp}_{\mathbf{B}_k(t)}(\dot{\mathbf{B}}_k(t))} = \norm{\sum_{i=0}\frac{\text{ad}_{\B_k(t)}^i}{(i+1)!}\dot{\mathbf{B}}_k(t)} \\ 
    & \leq \sum_{i=0}\norm{\frac{(\text{ad}_{\B_k(t)})^i}{(i+1)!}}\norm{\dot{\mathbf{B}}_k(t)} \\
    & \leq \frac{\exp{2r_c}-1}{2r_c}\norm{\A(t)} \\
    \end{split}
\end{equation}
Thus, combining Eq.~\eqref{gammaFullDef}, Eq.~\eqref{PicardnthNorm} and Eq.~\eqref{eq:bound G_k(t)}, we get
\begin{equation}\label{nonlinearConvError_1}
\begin{split}
     \norm{\left( e^{\OmegaB_k} - e^{\bar{\OmegaB}_k} \right) \bar{\w}(t_{k-1})} & = \norm{\gammab_k(t_k)} \\ & \leq C \max\limits_{t_{k-1}\leq p \leq t_k}\norm{\bar{\h}'(p)}\max\limits_{t_{k-1}\leq p \leq t_k}\norm{\bar{\w}(p)} (t_k - t_{k-1})^3 \\
     & \quad\quad \times \exp{\kappa \max\limits_{t_{k-1}\leq p \leq t}\norm{\A(p)}(t_k - t_{k-1})}\,, \end{split}
\end{equation}
where $\kappa = \frac{\exp{2r_c}-1}{2r_c}$ . 
%
%
%
%
%
\vspace{0.25cm}
 
We next seek to bound $\norm{\left( e^{\bar{\OmegaB}_k} - e^{\tilde{\OmegaB}_k} \right) \bar{\w}(t_{k-1})}$. This bound describes the error from using  the mid-point quadrature to approximate the integral (cf.~Eq.\eqref{eq:quadratureApprox}). We note
%
\begin{equation}\label{expTimeQuadError}
    \begin{split}
        \norm{\left( e^{\bar{\OmegaB}_k} - e^{\tilde{\OmegaB}_k} \right) \bar{\w}(t_{k-1})} & \leq \norm{\left( e^{\bar{\OmegaB}_k} - e^{\tilde{\OmegaB}_k} \right)} \norm{\bar{\w}(t_{k-1})} \\
        & = \norm{\int_{0}^{1} \frac{d}{dx}\left( e^{(1-x)\tilde{\OmegaB}_k} e^{x\bar{\OmegaB}_k} \right)dx} \norm{\bar{\w}(t_{k-1})} \\
        & = \norm{\int_{0}^{1} \left( e^{(1-x)\tilde{\OmegaB}_k}\left(\bar{\OmegaB}_k - \tilde{\OmegaB}_k  \right) e^{x\bar{\OmegaB}_k} \right)dx} \norm{\bar{\w}(t_{k-1})} \\
        & \leq \left( \int_{0}^{1} \norm{\left( e^{(1-x)\tilde{\OmegaB}_k}\left(\bar{\OmegaB}_k - \tilde{\OmegaB}_k  \right) e^{x\bar{\OmegaB}_k} \right)}dx \right) \norm{\bar{\w}(t_{k-1})} \\
        & \leq \left(\int_{0}^{1} \norm{ e^{(1-x)\tilde{\OmegaB}_k}}\norm{\left(\bar{\OmegaB}_k - \tilde{\OmegaB}_k  \right)} \norm{e^{x\bar{\OmegaB}_k}}dx \right) \norm{\bar{\w}(t_{k-1})} \\
        & \leq \norm{\left(\bar{\OmegaB}_k - \tilde{\OmegaB}_k  \right)} \left(\int_{0}^{1} \norm{e^{(1-x)\tilde{\OmegaB}_k}} \norm{e^{x\bar{\OmegaB}_k}} dx \right) \norm{\bar{\w}(t_{k-1})} \\
    \end{split}
\end{equation}
From the Taylor series expansion, the quadrature error corresponding to using a midpoint quadrature is given by $\abs{\int_a^b f(x) dx - f(\frac{a+b}{2})(b-a)} \leq \frac{1}{24}(b-a)^3\max\limits_{a \leq c \leq b}\abs{f''(c)}$. Thus, 
\begin{equation}\label{timeQuadError}
\begin{split}
    \norm{\left(\bar{\OmegaB}_k - \tilde{\OmegaB}_k  \right)} & = \norm{\int_{t_{k-1}}^{t_k}\begin{bmatrix} \mathbf{0} & \mathbf{I} \\ -\mathbf{\bar{H}}(\tau) & \mathbf{0} \end{bmatrix}d\tau - \begin{bmatrix} \mathbf{0} & \mathbf{I} \\ -\mathbf{\bar{H}}\left( \frac{t_{k-1} + t_{k}}{2} \right)\Delta t  & \mathbf{0}\end{bmatrix}} \\
    & = \norm{\int_{t_{k-1}}^{t_k}  \mathbf{\bar{H}}\left( \frac{t_{k-1} + t_{k}}{2} \right)\Delta t - \mathbf{\bar{H}}(\tau) d\tau } \\
    & \leq \frac{1}{24}(t_k - t_{k-1})^3 \max_{\substack{t_{k-1}\leq p \leq t}}\norm{\bar{\h}''(p)} \\
  \end{split}  
\end{equation}
%
%

%
\noindent We now seek to bound $\norm{e^{(1-x)\tilde{\OmegaB}_k}}$ and $\norm{e^{x\bar{\OmegaB}_k}}$ in Eq.~\eqref{expTimeQuadError}. To this end, we note
\begin{equation}\label{exp1mxOmegaB}
    \begin{split}
        \norm{e^{x\bar{\OmegaB}_k}} & \leq e^{\norm{x\bar{\OmegaB}_k}} \\
        & \leq e^{|x|\norm{\bar{\OmegaB}_k}} \\
        & \Big(\text{ we know that } \\ & \norm{\bar{\OmegaB}_k} = \norm{\int_{t_{k-1}}^{t_k}\A(\tau)d\tau} \leq \int_{t_{k-1}}^{t_k}\norm{\A(\tau)}d\tau \leq  \max\limits_{t_{k-1}\leq p \leq t}\norm{\A(p)}(t_k - t_{k-1}) \qquad \Big) \\
        &  \leq \exp{|x| \max\limits_{t_{k-1}\leq p \leq t}\norm{\A(p)}(t_k - t_{k-1})}   \quad  \forall x \in [0,1] \,,\\
        & \leq \exp{ \max\limits_{t_{k-1}\leq p \leq t}\norm{\A(p)}(t_k - t_{k-1})}\,,\\
    \end{split}
\end{equation}
and
\begin{equation}\label{expxOmegaB}
    \begin{split}
        \norm{e^{(1-x)\tilde{\OmegaB}_k}} & \leq e^{\norm{(1-x)\tilde{\OmegaB}_k}} \\
        & \leq e^{|(1-x)|\norm{\tilde{\OmegaB}_k}} \\
        &\Big( \text{ we know that } \norm{\tilde{\OmegaB}_k} = \norm{\A\left(\frac{t_{k-1} + t_{k}}{2} \right)}(t_k-t_{k-1}) \leq  \max\limits_{t_{k-1}\leq p \leq t}\norm{\A(p)}(t_k - t_{k-1}) \quad\Big) \\
        &  \leq \exp{|(1-x)| \max\limits_{t_{k-1}\leq p \leq t}\norm{\A(p)}(t_k - t_{k-1})
)}, \quad  \forall x \in [0,1] \,,\\
        & \leq \exp{ \max\limits_{t_{k-1}\leq p \leq t}\norm{\A(p)}(t_k - t_{k-1})}\,.
    \end{split}
\end{equation}
Substituting Eq.~\eqref{timeQuadError}, Eq.~\eqref{exp1mxOmegaB}, Eq.~\eqref{expxOmegaB} in Eq.~\eqref{expTimeQuadError}, we obtain
\begin{equation}\label{nonlinearConvError_2}
    \begin{split}
     \norm{\left( e^{\bar{\OmegaB}_k} - e^{\tilde{\OmegaB}_k} \right) \bar{\w}(t_{k-1})} & \leq \frac{1}{24}(t_k - t_{k-1})^3 \max_{\substack{t_{k-1}\leq p \leq t}}\norm{\bar{\h}''(p)} \norm{\bar{\w}(t_{k-1})} \\
     & \times \exp{2\max\limits_{t_{k-1}\leq p \leq t}\norm{\A(p)}(t_k - t_{k-1})} \,.
    \end{split}
\end{equation}
Using Eq.~\eqref{nonlinearConvError_1} and Eq.~\eqref{nonlinearConvError_2} in Eq.~\eqref{nonlinearConvError_split} provide the bound for $\norm{\left( e^{\OmegaB_k} - e^{\tilde{\OmegaB}_k} \right) \bar{\w}(t_{k-1})}$ in Eq.~\eqref{eq:boundNormBarW}. 
\bigskip

\noindent Using the definitions of $R_k^n$ from equations \eqref{eq:Rkn}-\eqref{eq:Adef}, we have,
 \begin{equation} 
   R_k^n = e^{\tilde{\OmegaB}_n}e^{\tilde{\OmegaB}_{n-1}}...e^{\tilde{\OmegaB}_{k+1}}\,,
\end{equation}
where
\begin{equation}
    \tilde{\OmegaB}_l = \mathbf{A}\left(\frac{t_{l-1} + t_{l}}{2} \right)  (t_l - t_{l-1}) \,.
\end{equation}
%
%
%
%
%

We note that 
\begin{equation}
    \norm{e^{\mathbf{A}\Delta t}} \leq  e^{w(\mathbf{A})\Delta t}\,,
\end{equation}
where $w(\mathbf{A})$ is the numerical abscissa of $\mathbf{A}$~\cite{trefethen2020spectra}, given as
\begin{equation}
    w(\mathbf{A}) = \text{sup}\left\{\text{Re}\left((\mathbf{q}^\dagger)^T \mathbf{A}\mathbf{q}\right)~ | ~ (\mathbf{q}^\dagger)^T\mathbf{q} = 1\right\}\,.
\end{equation}
Loosely speaking, the numerical abscissa measures the largest positive real part of the pseudo-spectrum of the matrix. Using the above bound, we can bound the norm of $R_k^n$ as
{
\renewcommand{\O}[1]{\tilde{\OmegaB}_{#1}}
 \begin{equation}\label{eq:RknInequality}
 \begin{split}
 \norm{R_k^n} & = \norm{e^{\tilde{\OmegaB}_n}e^{\tilde{\OmegaB}_{n-1}}...e^{\tilde{\OmegaB}_{k+1}}} \leq \exp{\max\limits_{t_{k}\leq p \leq t_{n}}w({\A(p)})(t_n - t_{k})} \,.
\end{split}
\end{equation}
}
We note that the above is not a tight bound and one can arrive at a tighter bound using some knowledge of the $\A(t)$. For instance, the pseudo-spectrum of $\A(t)$ is given by the eigen-spectrum of $\sqrt{-\mathbf{\bar{H}}(t)}$. For linear elastodynamics, wherein $\mathbf{\bar{H}}(t)$ is a positive semi-definite matrix, all eigenvalues of $\sqrt{-\mathbf{\bar{H}}(t)}$ are imaginary, leading to $w(\A(t))=0$. For nonlinear elastidynamics, $w(\A(t))$ will be goverened by the strength of nonlinearity, and hence, one can exploit the nature of the nonlinearity to find a tighter bound for $R_k^n$.

Substituting equations Eq.~\eqref{nonlinearConvError_1}, Eq.~\eqref{nonlinearConvError_2}, Eq.~\eqref{eq:RknInequality} and Eq.~\eqref{nonlinearConvError_split} in Eq.~\eqref{eq:boundNormBarW} we bound the total error at the $k^{th}$ time integration step as follows,
\begin{equation}
\begin{split}
    \norm{R_k^n}\norm{\left( e^{\OmegaB_k} - e^{\tilde{\OmegaB}_k} \right) \bar{\w}(t_{k-1})} & \leq \left\{C \max\limits_{t_{k-1}\leq p \leq t_k}\norm{\bar{\h}'(p)}\max\limits_{t_{k-1}\leq p \leq t_k}\norm{\bar{\w}(p)} (t_k - t_{k-1})^3 \right\} \\
     & \quad\quad \times \exp{\kappa \max\limits_{t_{k-1}\leq p \leq t_k}\norm{\A(p)}(t_k - t_{k-1})} \exp{\max\limits_{t_{k}\leq p \leq t_{n}}w(\A(p))(t_n - t_{k})}  \\
     & + \left\{\frac{1}{24}(t_k - t_{k-1})^3 \max_{\substack{t_{k-1}\leq p \leq t_k}}\norm{\bar{\h}''(p)} \norm{\bar{\w}(t_{k-1})} \right\} \\
     & \times \exp{2\max\limits_{t_{k-1}\leq p \leq t_k}\norm{\A(p)}(t_k - t_{k-1})} \exp{\max\limits_{t_{k}\leq p \leq t_{n}}w(\A(p))(t_n - t_{k})}  \\
     & \leq \left\{C \left[\max\limits_{t_{k-1}\leq p \leq t_k}\norm{\bar{\h}'(p)} + \max_{\substack{t_{k-1}\leq p \leq t_k}}\norm{\bar{\h}''(p)} \right]\max\limits_{t_{k-1}\leq p \leq t_k}\norm{\bar{\w}(p)} (t_k - t_{k-1})^3 \right\} \\
     & ~~~~\times  \exp{\max{(\kappa,2)} \max\limits_{t_{k-1}\leq p \leq t_n}\norm{\A(p)}(t_k - t_{k-1})} \exp{\max\limits_{t_{k}\leq p \leq t_{n}}w(\A(p))(t_n - t_{k})}\,.
\end{split}
\end{equation}

\noindent Assuming equal time-step sizes of $\Delta t$ for all time-steps, the bound for global error is given by

\begin{equation}\label{eq:totalError}
    \begin{split}
        \norm{\bar{\w}(t_n) - \bar{\w}_n} & \leq \Bigg\{C \left[\max\limits_{t_{0}\leq p \leq t_n}\norm{\bar{\h}'(p)} + \max_{\substack{t_{0}\leq p \leq t_n}}\norm{\bar{\h}''(p)} \right]\max\limits_{t_{0}\leq p \leq t_n}\norm{\bar{\w}(p)} \\
     & \quad\quad \times  \exp{\max{(\kappa,2)} \max\limits_{t_{0}\leq p \leq t_n}\norm{\A(p)}(t_n - t_{0})} \exp{\max\limits_{t_{0}\leq p \leq t_{n}}w(\A(p))(t_n - t_{0})}\Bigg\} (\Delta t)^2 \,.
    \end{split}
\end{equation}

Assuming that $\h(t)$ is twice-differentiable and $\norm{\h(t)},\norm{\h'(t)},\norm{\h''(t)}$, and $w(\A(t))$ are bounded for all ${t \in [t_0,t_n]}$, the global error demonstrates a second-order time propagator. 





  \subsection{Energy conservation in linear elastodynamics
  }\label{subsec:energy conservation}
We recall the reformulated governing equation of elastodynamics as given in \eqref{eq:discrete dynamic equilibrium using Hbar},
    	\begin{equation*} 
    	  \ddot{\bar{\mathbf{u}}} + \bar{\mathbf{H}}(\mathbf{u}) \bar{\mathbf{u}} = \bar{\mathbf{P}} \,,
    	\end{equation*}
where
   	\begin{equation}
    	    \bar{\mathbf{u}} = \mathbf{M}^{1/2}\mathbf{u}\,, \quad \bar{\mathbf{H}}(\mathbf{u}) = \mathbf{M}^{-1/2}\mathbf{H}(\mathbf{u})\mathbf{M}^{-1/2} \,, \quad \bar{\mathbf{P}} = \mathbf{M}^{-1/2}\mathbf{P} \,.\notag
    	\end{equation}
%
%
%
As we are interested in demonstrating energy conservation in linear elastodynamics, we consider the case where $\bar{\mathbf{P}}=0$ (no external loading) and $\bar{\mathbf{H}}(\mathbf{u}) = \bar{\mathbf{H}}$ (linear elastodynamics). Thus, for this case, the governing equation reduces to
\begin{equation}
      \ddot{\bar{\mathbf{u}}} + \bar{\mathbf{H}} \bar{\mathbf{u}} = 0 \,.
\end{equation}
We note that $\bar{\mathbf{H}}$ is a symmetric, positive definite matrix.  Refromulating this second order equation into a system of first order equations gives the equivalent of equation \eqref{eq: 1st order reformulation} for linear elastodynamics with no external load,
\begin{equation}
    	 \frac{d}{dt}\begin{bmatrix}	 {\bar{\mathbf{u}}}(t) \\	 {\bar{\mathbf{v}}}(t) \end{bmatrix}  = \begin{bmatrix}
    	       0 & \mathbf{I} \\
    - \bar{\mathbf{H}} & 0
    	     \end{bmatrix} \begin{bmatrix}
    	 \bar{\mathbf{u}}(t) \\
    	 \bar{\mathbf{v}}(t) 
    	\end{bmatrix} \,.
\end{equation}

We now seek to demonstrate the energy conservation afforded by the exponential propagator approach for linear elastodynamics. To this end, let $E_{n}$ denotes the energy at time $t_n$, which is given by 
\begin{equation}\label{energyExpr}
    E_{n} = \frac{1}{2}\bar{\mathbf{v}}^T\bar{\mathbf{v}} + \frac{1}{2}\bar{\mathbf{u}}^T\bar{\mathbf{H}}\bar{\mathbf{u}} \,.
\end{equation} 
We note that Eq.~\eqref{energyExpr} can be expressed as
\begin{equation}\label{energyw_nExpr}
    E_n = \frac{1}{2}\bar{\mathbf{w}}_{n}^T\bar{\mathbf{G}}\bar{\mathbf{w}}_n \,,
\end{equation}
where
\begin{equation}
    \bar{\mathbf{G}} = \begin{bmatrix} \bar{\mathbf{H}} & 0 \\ 0 & \mathbf{I}
    \end{bmatrix},
    \bar{\mathbf{w}}_n = \begin{bmatrix} \bar{\mathbf{u}}(t_n) \\ \bar{\mathbf{v}}(t_n)\end{bmatrix}\,.
\end{equation}
Now, $\bar{\mathbf{w}}_n$ can be written as
\begin{equation}\label{w_nExprEnergy}
    \bar{\mathbf{w}}_n = e^{\mathbf{\OmegaB}(t_n - t_0)}\bar{\mathbf{w}}_0\,, \quad\text{ where } \mathbf{\OmegaB} = \begin{bmatrix} 0 & \mathbf{I} \\ -\bar{\mathbf{H}} & 0 \end{bmatrix}\,.
\end{equation}
Without loss of generality, we assume $t_0 = 0$. Substituting equation \eqref{w_nExprEnergy} in equation \eqref{energyw_nExpr} we get
\begin{equation}
    E_n = \frac{1}{2}\left[\bar{\mathbf{w}}_0^Te^{\mathbf{\OmegaB}^Tt_n}\bar{\mathbf{G}}e^{\mathbf{\OmegaB}t_n}\bar{\mathbf{w}}_0 \right]
\end{equation}
We note that
{
\newcommand{\sqH}{\sqrt{\bar{\h}}}
\newcommand{\tn}{t_n}
\begin{equation}\label{eq:e^OmegaBt}
   e^{\mathbf{\OmegaB}\tn} = \begin{bmatrix} \cos{(\sqH \tn)} & \frac{1}{\sqH} \sin{(\sqH \tn)}  \\   -\sqH \sin{(\sqH \tn)} & \cos{(\sqH \tn)}  \end{bmatrix}
\end{equation}
}
and

%
{\newcommand*{\wo}{\bar{\mathbf{w}}_0}
\newcommand*{\woT}{\bar{\mathbf{w}}_0^T}
\renewcommand*{\G}{\bar{\mathbf{G}}}
\newcommand{\OT}{\mathbf{\OmegaB}^T}
\renewcommand{\O}{\mathbf{\OmegaB}}
\newcommand{\EOT}[1]{e^{\mathbf{\OmegaB}^T#1}}
\newcommand{\EO}[1]{e^{\mathbf{\OmegaB}#1}}
\newcommand{\sqH}{\sqrt{\bar{\h}}}
\newcommand{\dt}{t_n}

\begin{equation*}
\begin{split}
    \EOT{\dt}\G\EO{\dt} & = \begin{bmatrix} \cos{(\sqH \dt)} &  -\sqH \sin{(\sqH \dt)}  \\  \frac{1}{\sqH} \sin{(\sqH \dt)} & \cos{(\sqH \dt)}  \end{bmatrix} \begin{bmatrix} \bar{\mathbf{H}} & 0 \\ 0 & \mathbf{I} 
    \end{bmatrix} \begin{bmatrix} \cos{(\sqH \dt)} & \frac{1}{\sqH} \sin{(\sqH \dt)}  \\   -\sqH \sin{(\sqH \dt)} & \cos{(\sqH \dt)}  \end{bmatrix} \\
    & = \begin{bmatrix} \cos{(\sqH \dt)} &  -\sqH \sin{(\sqH \dt)}  \\  \frac{1}{\sqH} \sin{(\sqH \dt)} & \cos{(\sqH \dt)}  \end{bmatrix} \begin{bmatrix} \bar{\mathbf{H}}\cos{(\sqH \dt)} & \sqH \sin{(\sqH \dt)}  \\   -\sqH \sin{(\sqH \dt)} & \cos{(\sqH \dt)}  \end{bmatrix} \\
    & = \begin{bmatrix}
        \bar{\mathbf{H}} & 0 \\ 0 & \mathbf{I}
    \end{bmatrix} = \G\,.
    \end{split}
\end{equation*}
}
Thus, $E_n = E_0$ for any $n$, and we have shown energy conservtion.



\subsection{Symplecticity}
A linear mapping $A$ is symplectic if and only if it satisfies the following condition~\cite{hairer2006structure}
\begin{equation}\label{eq:symplectic condition}
    A^T\begin{bmatrix}
        0 & I \\ -I & 0 
    \end{bmatrix} A \,
    = \,
    \begin{bmatrix}
        0 & I \\ -I & 0 
    \end{bmatrix}\,.
\end{equation}
The exponential propagator in each propagating step acts as a linear mapping given by,
{
\newcommand{\sqH}{\sqrt{\bar{\mathbf{H}}_k}}
\newcommand{\dt}{\Delta t}
\begin{equation}\label{eq:Jacobian}
   e^{\mathbf{\tilde{\OmegaB}_k}\dt} = \begin{bmatrix} \cos{(\sqH \dt)} & \frac{1}{\sqH} \sin{(\sqH \dt)}  \\   -\sqH \sin{(\sqH \dt)} & \cos{(\sqH \dt)}  \end{bmatrix} \,.
\end{equation} 
where $\bar{\mathbf{H}}_k = \bar{\mathbf{H}}\left( \frac{t_{k-1} + t_{k}}{2} \right)$.
}
The linear map given by \eqref{eq:Jacobian} satisfies the condition \eqref{eq:symplectic condition}, thus making the exponential propagator a symplectic integrator. Symplecticity provides a stable and accurate \longtime behavior for the exponential propagator.

\section{Numerical results} \label{sec:numerical results}

In this section, we demonstrate the various numerical aspects of the exponential propagator approach, such as characteristics of the subspace approximation, rate of convergence, computational efficiency, and parallel scalability. We also demonstrate the suitability of the exponential propagator approach for \longtime elastodynamic simulations. To this end, we study the proposed exponential propagator approach compared with state-of-the-art methods such as the \Newmark method~\cite{newmark1959method}~($\beta=0.25, \gamma=0.5$), the \HHTA method~\cite{hilber1977improved}~($\alpha=-0.33, \beta=0.25, \gamma=0.5$), and energy-momentum conserving methods~\cite{gonzalez2000exact,simo1992exact} in protypical hyperelastic systems. In our benchmark studies, we focus on geometric models that are of practical significance in real-world applications. These models include the cantilever beam and the plate. We use hexahedral finite elements with quadratic Lagrange interpolation polynomials to spatially discretize the geometry. To generate the unstructured meshes, we use the software CUBIT™ Geometry~\cite{osti1457612}'s Mesh Generation Toolkit. The description and setup of our benchmark systems are as follows:

\begin{figure}[h]
     \captionsetup[subfigure]{justification=centering}
     \begin{subfigure}[b]{0.49\textwidth}
         \centering
         \includegraphics[scale=0.12]{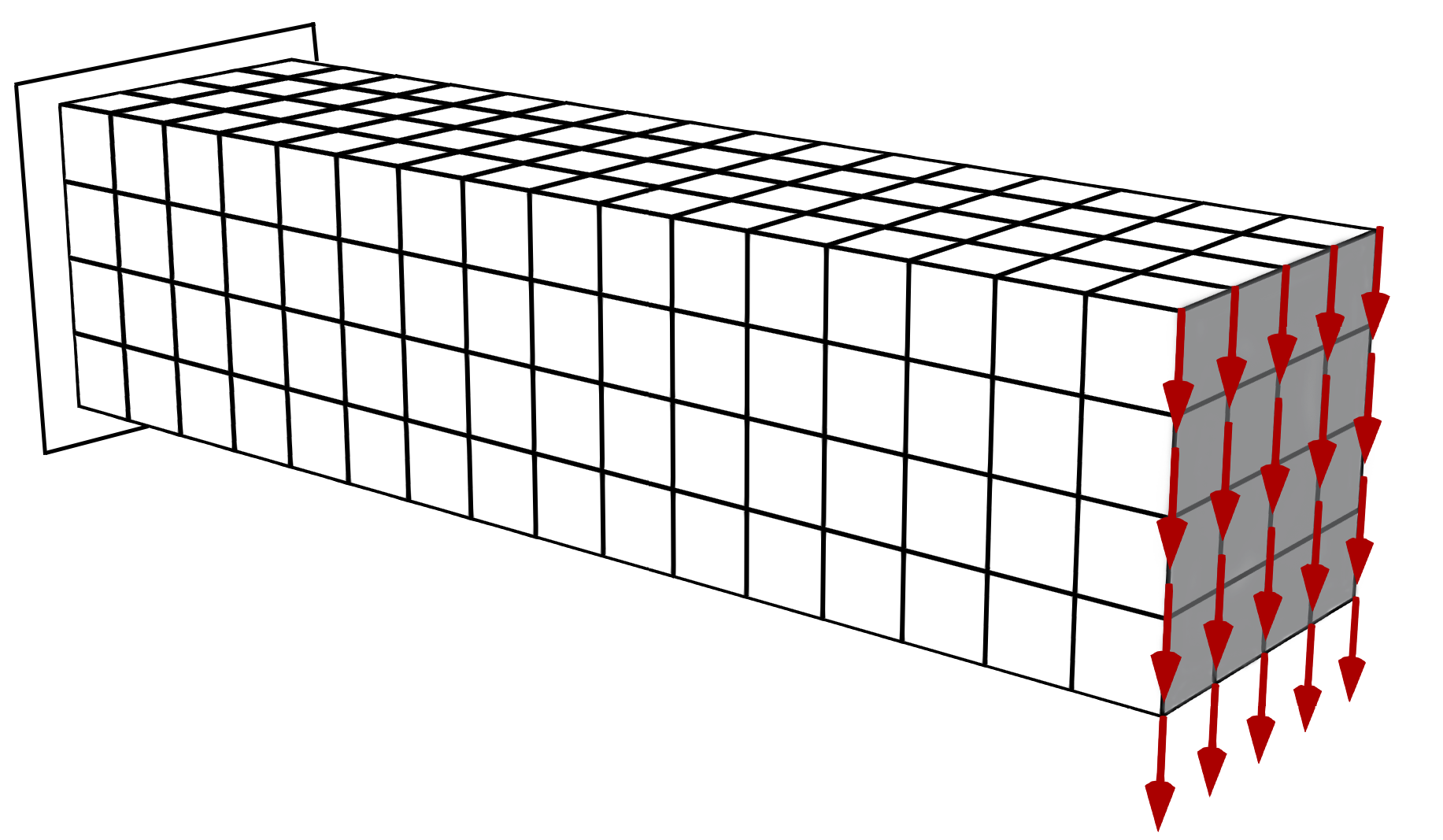}
         \caption{Undeforomed cantilever beam subjected to an initial loading.}
         \label{fig:Cantilever undeformed}
     \end{subfigure}
     \begin{subfigure}[b]{0.49\textwidth}
         \centering
     \includegraphics[scale=0.12]{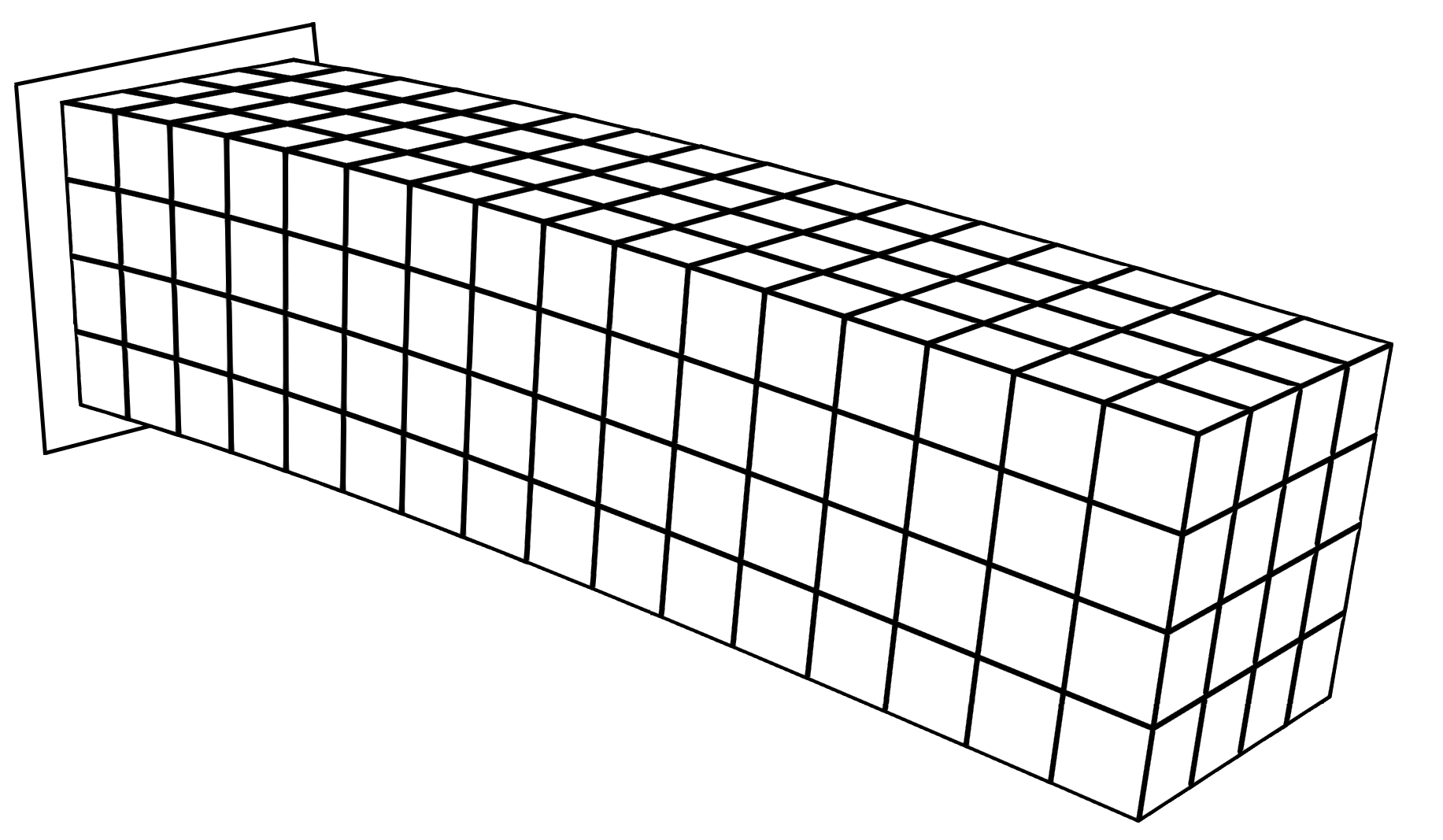}
         \caption{Deformed cantilever beam under the initial loading.}
         \label{fig:Cantilever deformed}
     \end{subfigure} 
        \caption{Geometry and finite-element mesh of the cantilever beam used for the analysis.}
        \label{fig:cantilever411_0p25_schematic}
\end{figure}


\begin{figure}[h]
     \captionsetup[subfigure]{justification=centering}
     \begin{subfigure}[b]{0.49\textwidth}
         \centering
         \includegraphics[scale=0.13]{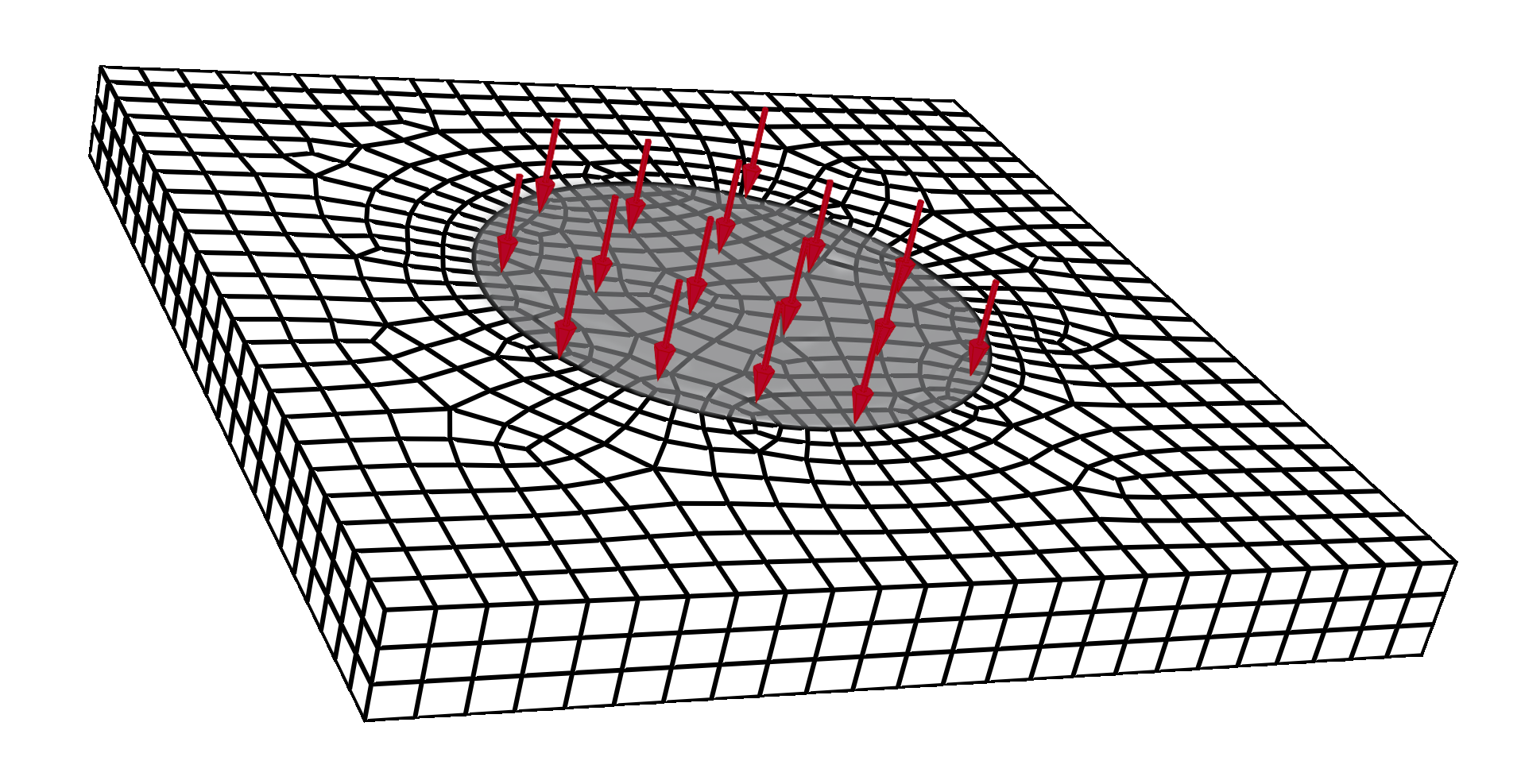}
         \caption{Undeformed plate subjected to an initial loading.}
         \label{fig:Plate undeformed}
     \end{subfigure}
     \begin{subfigure}[b]{0.49\textwidth}
         \centering
     \includegraphics[scale=0.13]{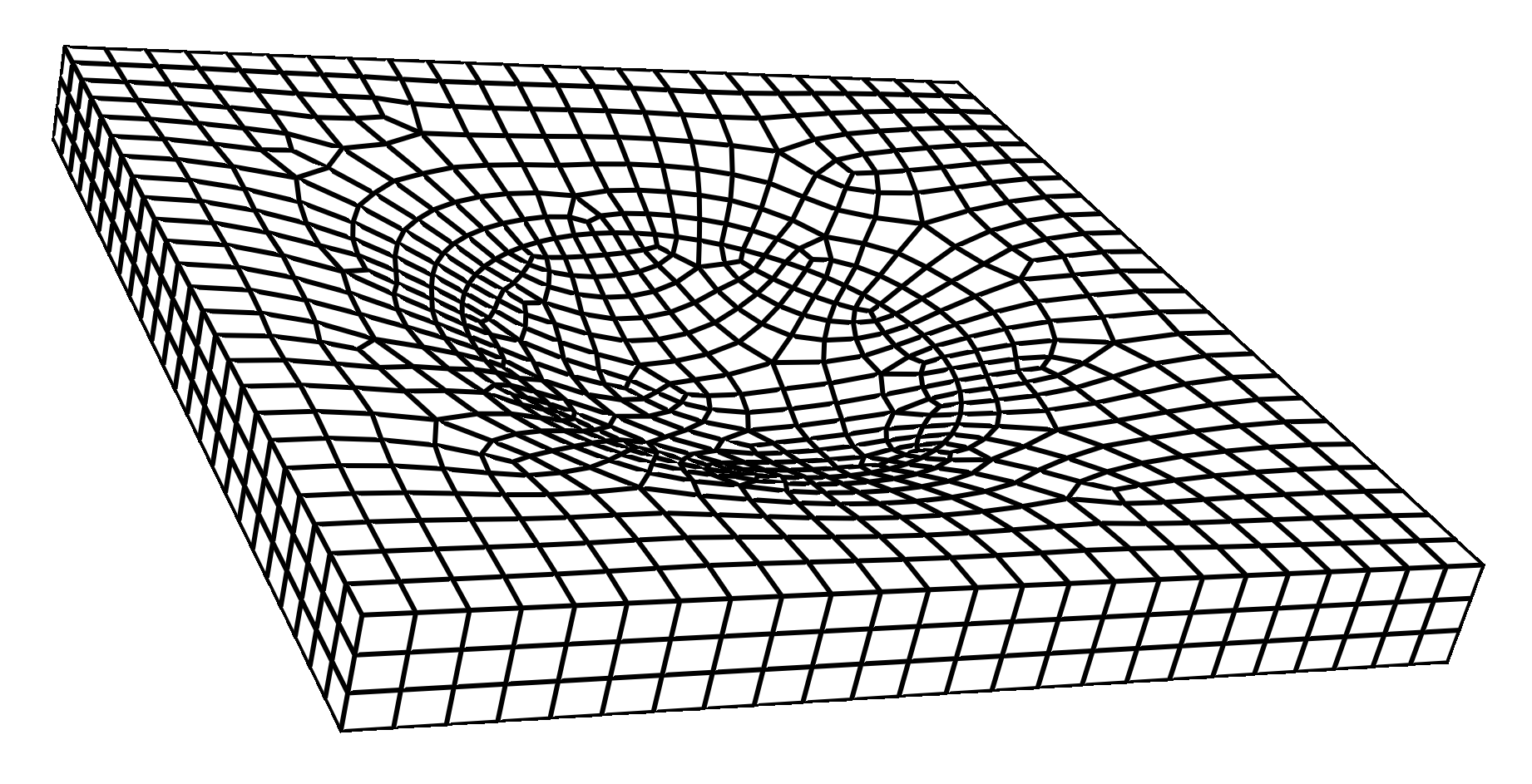}
         \caption{Deformed plate under the initial loading.}
         \label{fig:Plate deformed}
     \end{subfigure} 
        \caption{Geometry and finite-element mesh of the plate used for the analysis.}
        \label{fig:plate_schematic}
\end{figure}

\begin{itemize}
    \item A cantilever beam is initially displaced by a load applied to the surface at the free end and the evolution over time is studied. Fig.(\ref{fig:Cantilever undeformed}) shows the schematic of the undeformed geometry. The shaded region depicts the surface upon which an initial distributed load is applied, and Fig.(\ref{fig:Cantilever deformed}) shows the schematic of the deformed geometry under the load. \label{cantilever model}
    \item A square plate is initially deformed by a distributed circular load applied in the center and the dynamics is studied. Fig.(\ref{fig:Plate undeformed}) shows the schematic of the undeformed geometry. The shaded region depicts the surface upon which an initial distributed load is applied, and Fig.(\ref{fig:Plate deformed}) shows the schematic of the deformed geometry under the load. \label{plate model}
\end{itemize}

All numerical studies are conducted on a parallel computing cluster with the following configuration: 2$\times$~3.0 GHz Intel Xeon Gold 6154 (Skylake) CPU nodes with 36 processors (cores) per node, 187 GB memory per node, and Infiniband HDR100 networking between all nodes for MPI communications.

\subsection{Subspace characteristics}\label{subsec:subspace characteristics}
Taking advantage of the exponential propagator approach relies on efficient evaluation of the matrix exponential. For this purpose, we study the charactertistics of the Krylov subspace-based approximation and its efficacy in evaluating the action of a matrix exponential on a vector. First, we study the behavior of the subspace approximation error with increasing subspace sizes across various values of time-step sizes. To isolate the error from the subspace approximation, after one time propagation step, we contrast the propagated solution acquired through the Krylov subspace-based approximation using a subspace of size $m$~($\w^h_m(\D t)$) with the one obtained through direct exponentiation~($\w^h(\D t)$). Here, by direct exponentiation, we refer to the explicit Taylor series evaluation of the matrix exponential to machine precision. The relative error in the subspace approximation~($\varepsilon_m$) is measured as follows:
 \begin{equation}\label{eq:subspace error}
           \varepsilon_m = \frac{\norm{\w^h(\D t) - \w^h_m(\D t)}_2}{\norm{\w^h(\D t)}_2}
 \end{equation}
Since the evaluation of matrix exponential acting on a vector is performed within a time-step of propagation, wherein the problem is treated as linear elasticity, it is reasonable to limit our analysis of the subspace approximation error to the case of linear elastodynamics. In the context of linear elastodynamics, the exponential propagator approach does not impose any stability limit on the \timestep sizes. This enables us to study subspace characteristics for a wide range of \timestep sizes. We further present the results for a system with relatively fewer degrees of freedom, where direct exponentiation is tractable. 
\begin{figure}[!htbp]
     \centering
     \begin{subfigure}[b]{1.0\textwidth}
         \centering
         \includegraphics[width=0.7\textwidth]{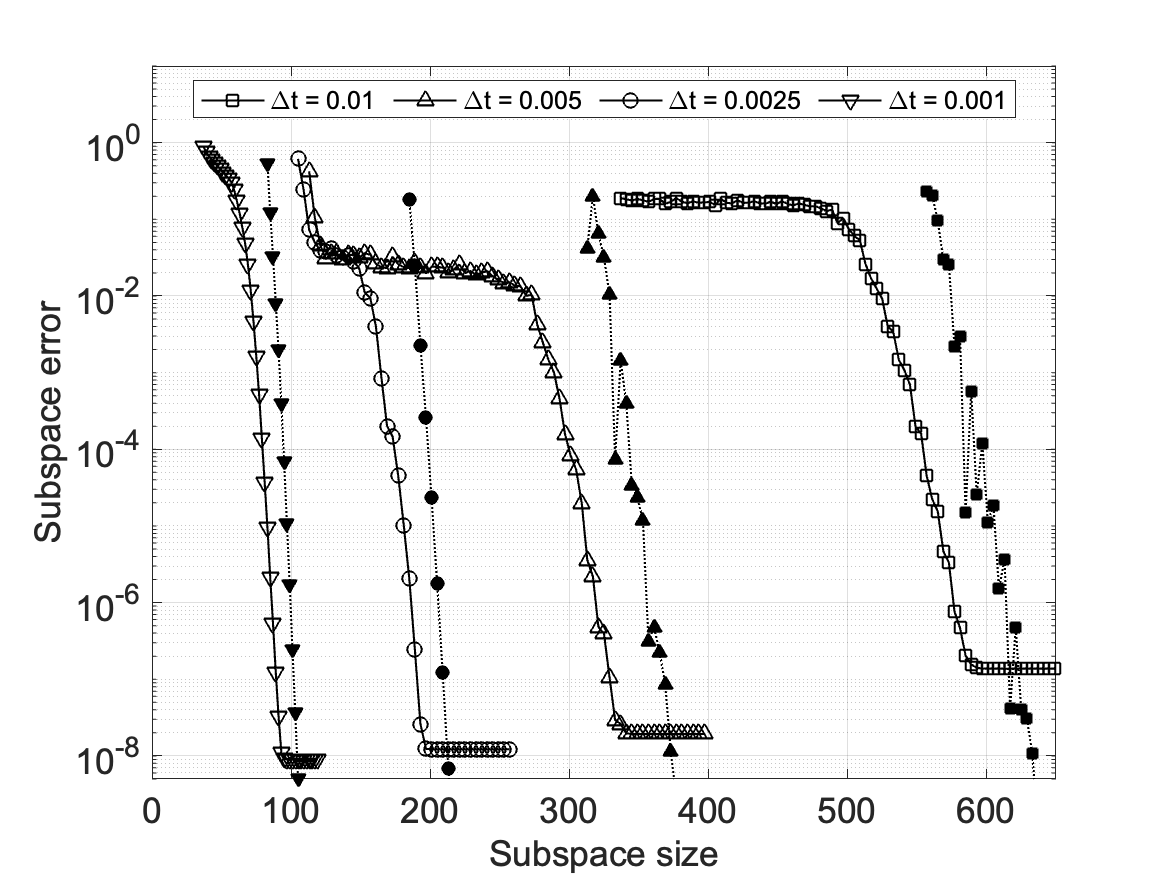}
     \end{subfigure}
     \begin{subfigure}[b]{1.0\textwidth}
         \centering
         \includegraphics[width=0.7\textwidth]{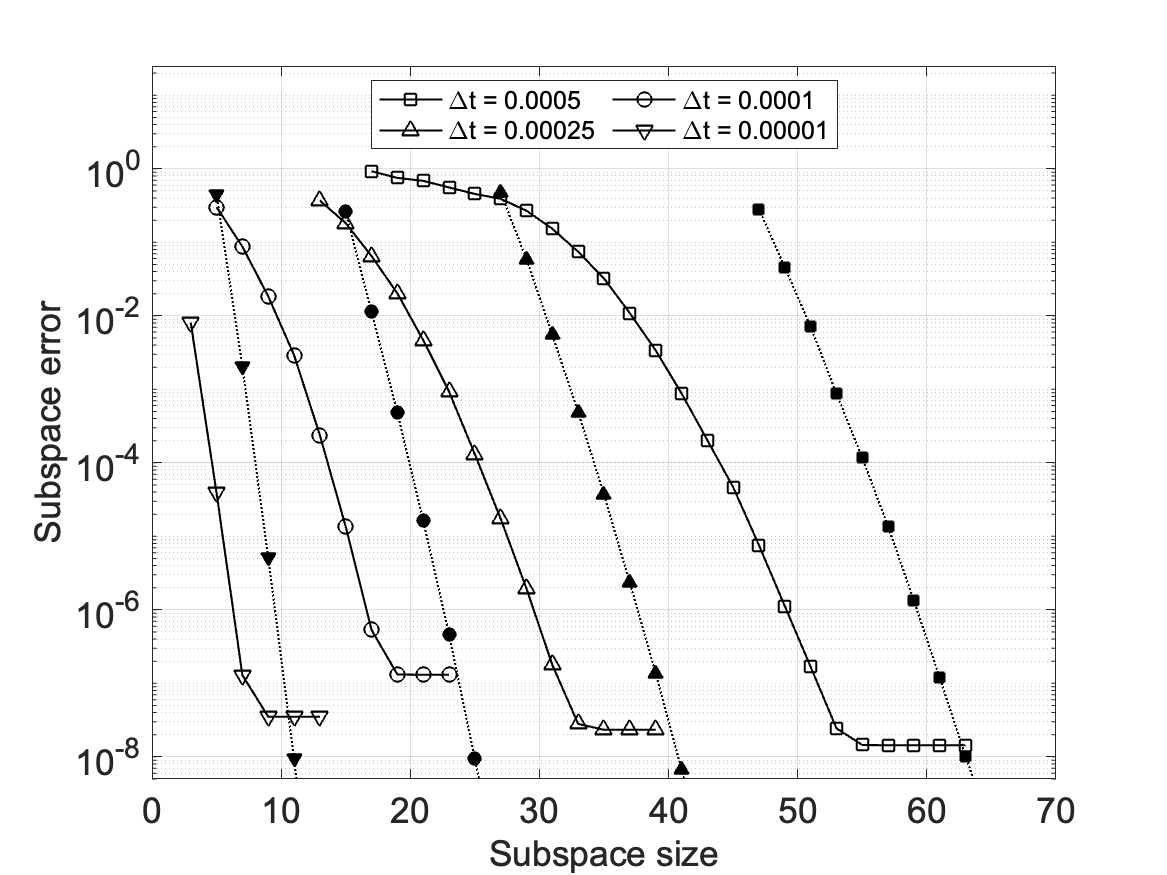}
     \end{subfigure} 
        \caption{Subspace error for various subspace sizes~($m$) and time-step sizes~($\Delta t$). The hollow markers denote the actual error~($\varepsilon_m$) and the corresponding solid markers denote the respective error estimates~($\epsilon_m$).}
        \label{fig:subspace error vs subspace size}
\end{figure}
\begin{figure}[h]
    \centering
    \includegraphics[width=0.7\textwidth]{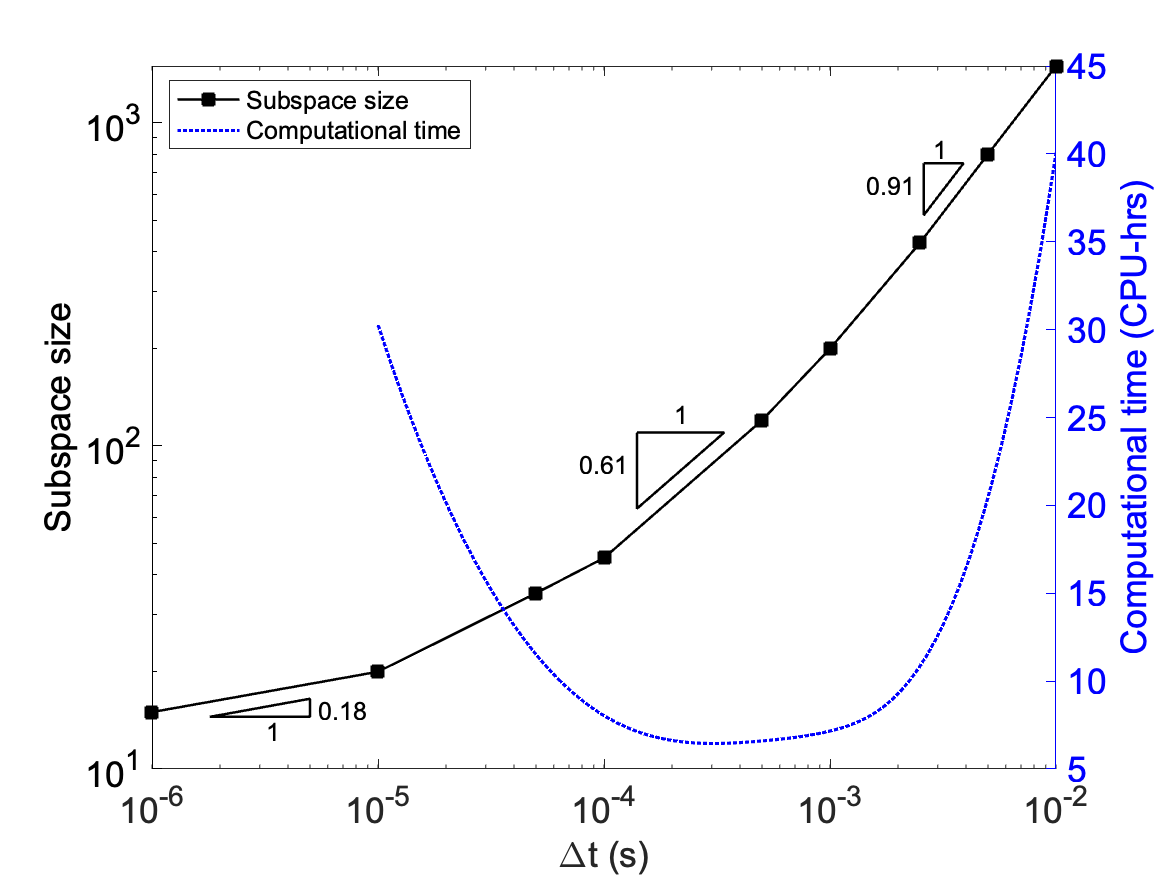}
    \caption{Subspace size, computational time~(CPU-hrs) vs time-step size~($\Delta t$).}
    \label{fig:subspace size, computational time vs time-step size}
\end{figure}
Figure~\ref{fig:subspace error vs subspace size} is a plot of the subspace approximation error~($\varepsilon_m$) given in eq.(\ref{eq:subspace error}) for various combinations of \timestep size and subspace size, studied on a structural steel cantilever beam modeled using linear elasticity with $1200$ DoFs. 
Here, we consider a problem with fewer DoFs since exact exponentiation is nearly intractable for larger \nonnormal matrices. 
From the figure, we can observe that the drop in subspace error is exponential beyond a threshold subspace size, and this threshold value increases with increasing \timestep size, i.e., a higher \timestep size requires a larger subspace for the same desired accuracy. We also note that the error estimate~($\epsilon_m$) given in eq.(\ref{eq:subspace error estimate}) accurately traces an upper bound of the actual error and can thus serve as a good metric to truncate the Arnoldi iteration. 

We also compare the relative computational costs incurred at various \timestep sizes for the same desired accuracy. This provides insight into the optimal choice for \timestep size that allows maximum computational benefit. For each \timestep size, the dimension of the subspace that gives the desired accuracy is determined using the previous study. Then the time propagation is performed at this subspace size until a final time, and the corresponding computational time is measured in node-hours. Figure~\ref{fig:subspace size, computational time vs time-step size} presents a plot showing the size of the subspace required at different \timestep sizes and the associated computational time, measured as described. This was carried out on a structural steel cantilever beam modeled using linear elasticity with $29,376$ DoFs using 16 processors, run to a final time of $t_n=1s$. It is worth noting that at lower time-step sizes, the subspace size needed for achieving a certain level of accuracy increases at a slower pace than the \timestep size. However, at higher \timestep sizes, the subspace requirements grow at a rate similar to the \timestep size. This observation is corroborated by the trends observed in computational time. At higher \timestep sizes, when the rate of growth of the necessary subspace size aligns with the rate of growth of the \timestep size, the benefit of using a larger \timestep size is offset by the increased expense incurred in generating and exponentiating a larger Hessenberg~($\mathbf{H}_m$) matrix. Therefore, it is most advantageous to select the largest time-step size within the regime where the rate of growth of required subspace size does not outpace the rate of growth of time-step size in order to realize maximum computational benefit. As we go further, we will see that at these optimal time-step sizes, we are able to achieve substantial speed-up and superior accuracy compared to the conventional techniques used in elastodynamics. 
%
        \subsection{Rates of convergence and computational cost}\label{subsec:rate of convergence and computational efficiency}

We now present the accuracy and computational efficiency afforded by the exponential propagator over some of the widely used time propagation schemes, namely the \Newmark and \HHTA methods. To this end, we first begin our comparison of the temporal convergence, measured as the relative $l^2$ errors of the propagated vectors, computed using different time-step sizes, at time $t_n$. Since we lack the knowledge of the actual solution at time $t_n$, the solution computed using a very high temporal refinement of \Newmark and exponential propagator calculation is used as a reference for measuring the errors. These reference solutions in displacement and velocity are denoted as $u^{h}_{\text{ref}}(t_n)$ and $v^{h}_{\text{ref}}(t_n)$, resepctively. If $u^{h}_{\Delta t}(t_n)$ and $v^{h}_{\Delta t}(t_n)$ denote the solutions in displacement and velocity respectively, computed using a time-step size of $\Delta t$, then the relative error in displacement is computed as

\begin{equation}
    \frac{\norm{\mathbf{u}^{h}_{\Delta t}(t_n) - \mathbf{u}^{h}_{\text{ref}}(t_n)}_2}{\norm{\mathbf{u}^{h}_{\text{ref}}(t_n)}_2} \,,
\end{equation}
and the relative error in velocity is computed as
\begin{equation}
    \frac{\norm{\mathbf{v}^{h}_{\Delta t}(t_n) - \mathbf{v}^{h}_{\text{ref}}(t_n)}_2}{\norm{\mathbf{v}^{h}_{\text{ref}}(t_n)}_2}\,.
\end{equation}
Since the underlying spatial discretization remains fixed, these error metrics solely reflect the error from temporal discretization. For a given time-step size, the exponential propagator is initially run for a few time-steps to estimate the size of the subspace required. Further, for a specific desired accuracy, our numerical experiments suggest that the required subspace size does not change much during the entire time propagation. Thus, we fix the subspace size for the duration of the time propagation to avoid the overhead costs associated with adaptively tuning the subspace size. 
We now present the results from the various benchmark studies conducted to assess the accuracy and computational efficiency of the exponential propagator method. 
            
\subsubsection{Linear cantilever}\label{subsubsec:computational efficiency:linear cantilever}
            
\begin{figure}[h]
    \centering
    \includegraphics[width=0.7\textwidth]{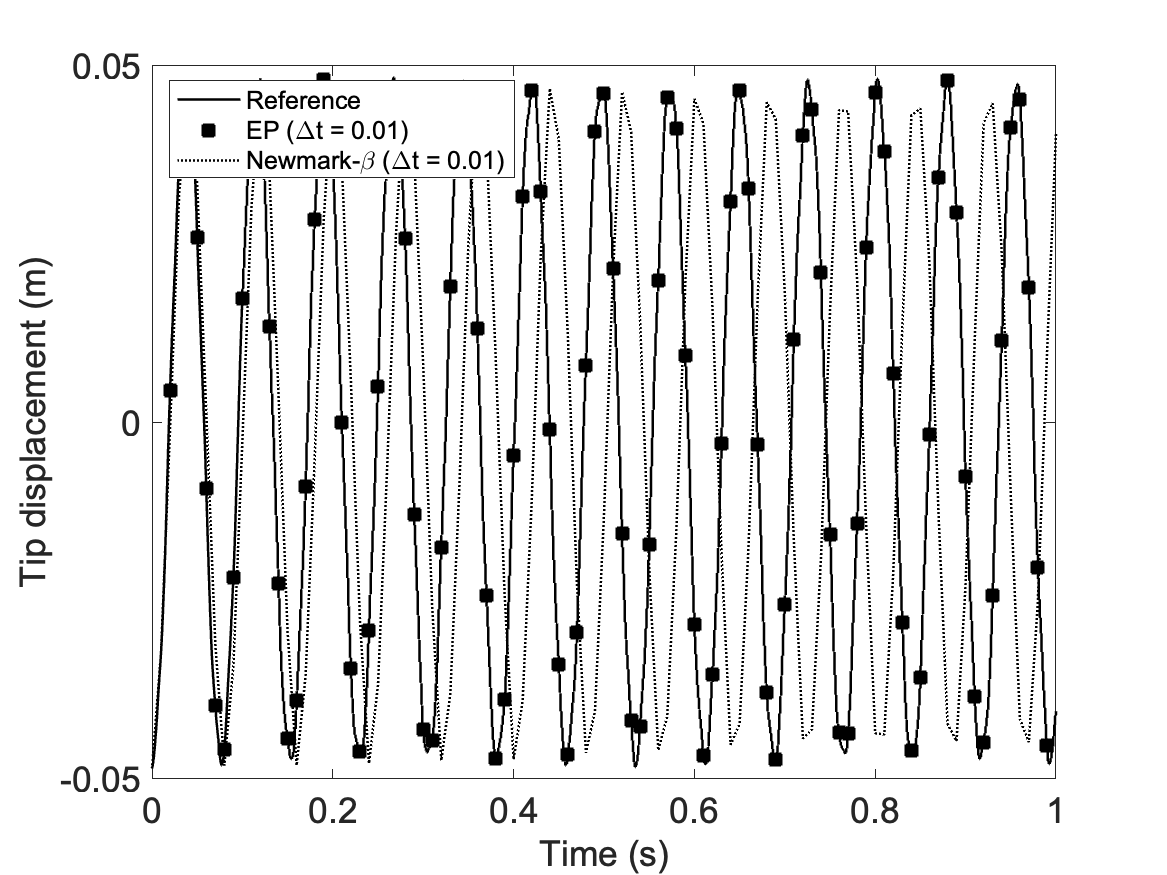}
    \caption{Evolution of tip displacement with time computed using the exponential propagator~(EP) and \Newmark methods.}
    \label{fig:linear tip displacement}
\end{figure}

The first benchmark system we consider is the steel cantilever beam, a representative problem to assess the performance of the exponential propagator for linear elastodynanics. 
We study the time evolution of the beam for a time of 1.0~s, i.e, $t_n$ = 1.0~s. In order to validate the accuracy of the exponential propagator, we provide in Fig.\ref{fig:linear tip displacement} the tip displacement of the cantilever with time as computed by the exponential propagator with a time-step size $\Delta t$=0.01~s, and compare against the reference calculation. The figure not only demonstrates a close agreement between exponential propagator and the reference calculation, but it also clearly establishes the superior accuracy affored by the propagator in comparison to a \Newmark calculation using the same time-step size. 

\begin{figure}[h]
    \centering
     \includegraphics[width=0.7\textwidth]{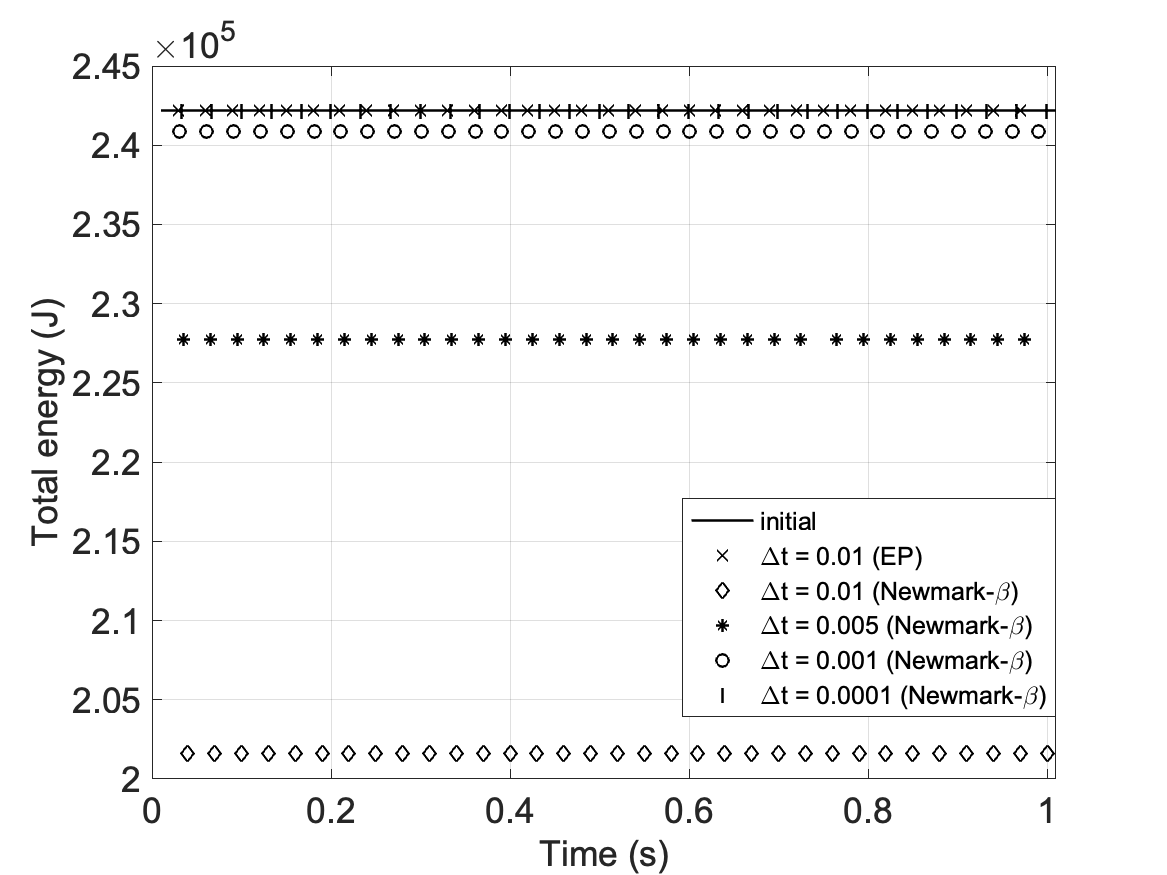}
    \caption{Time histories of the total energy computed using the exponential propagator~(EP) and \Newmark method at various time-step sizes~($\Delta t$).}
    \label{fig:energy conservation comparison}
\end{figure}

\begin{figure}[h]
    \centering
     \includegraphics[width=0.7\textwidth]{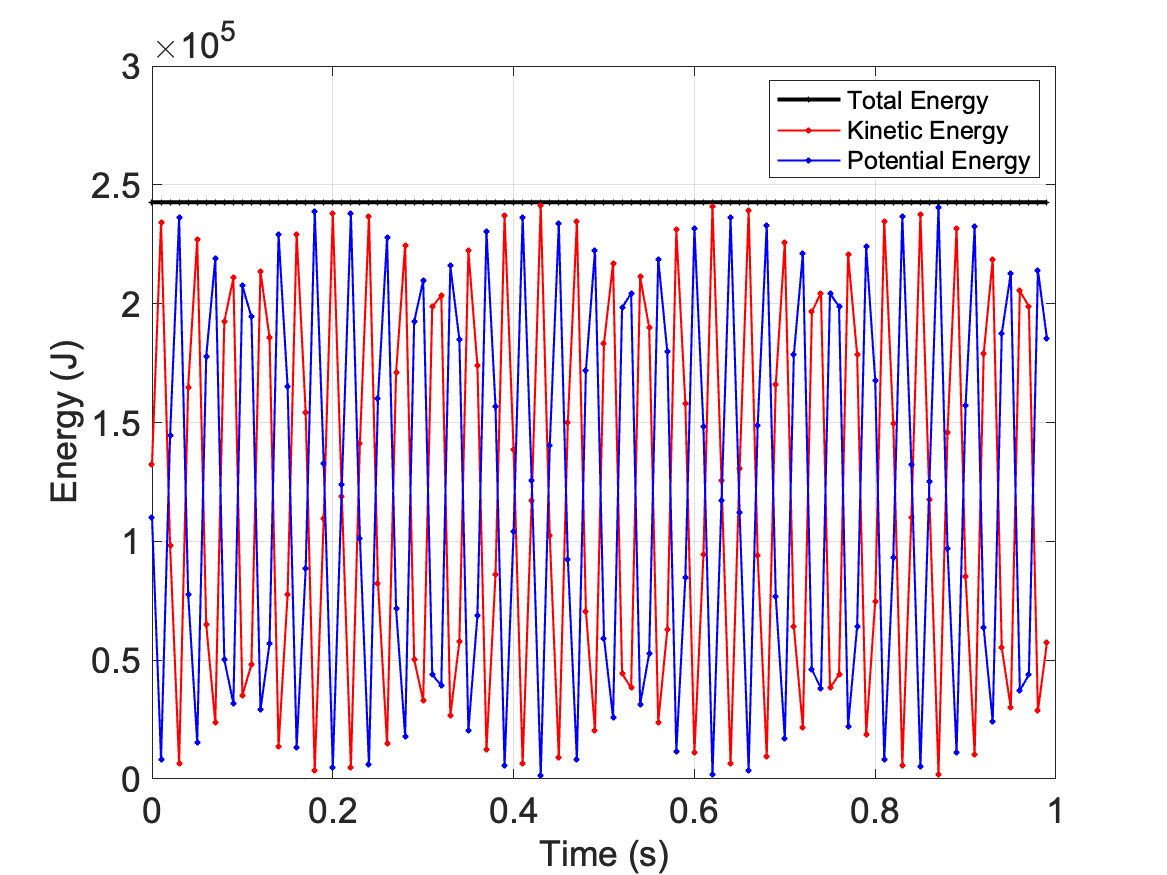}
    \caption{Time histories of the total energy, kinetic energy and potential energy computed using the exponential propagator method.}
    \label{fig:linear energy conservation}
\end{figure}

It is well known that the \Newmark method is energy-conserving for linear elastodynamics~\cite{jog2015conservation}, however, it does not guarantee an accurate solution. This is evident from results in Fig.~\ref{fig:energy conservation comparison}, where we plot the total energy with time at various time-step sizes. At all the time-step sizes, \Newmark method exhibits a constant energy profile, however at larger time-step sizes the predictions are highly inaccurate.  On the contrary, the exponential method also is energy conserving as observed in Fig.~\ref{fig:energy conservation comparison} and Fig.~\ref{fig:linear energy conservation}, but achieves the same accuracy as a refined \Newmark calculation~($\Delta t = 0.0001$) at a time-step size that is 100$\times$ larger~($\Delta t = 0.01$) as observed in  Fig.~\ref{fig:energy conservation comparison}.

%
\begin{figure}[!htbp]
     \centering
     \begin{subfigure}[b]{1.0\textwidth}
         \centering
         \includegraphics[width=0.7\textwidth]{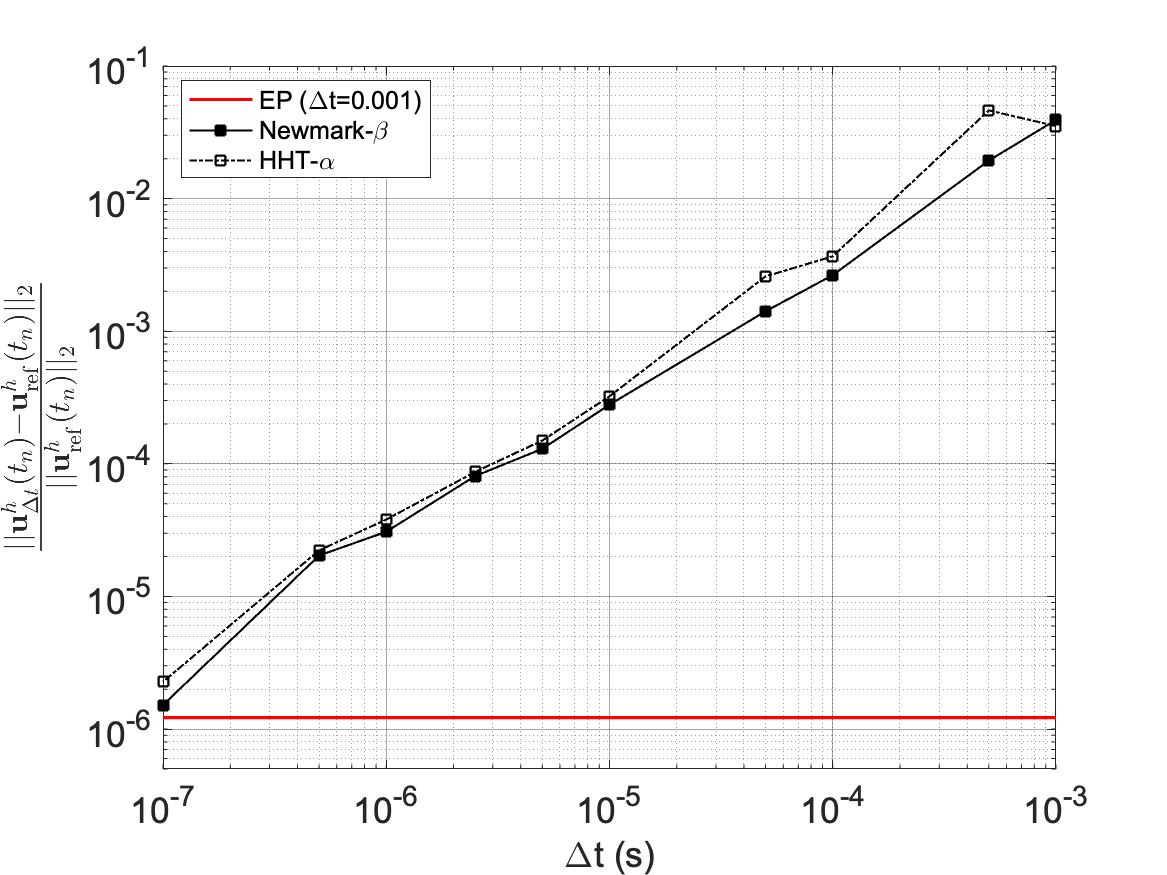}
         \caption{}
         \label{fig:linear numerical convergence - a}
     \end{subfigure}
     \begin{subfigure}[b]{1.0\textwidth}
         \centering
         \includegraphics[width=0.7\textwidth]{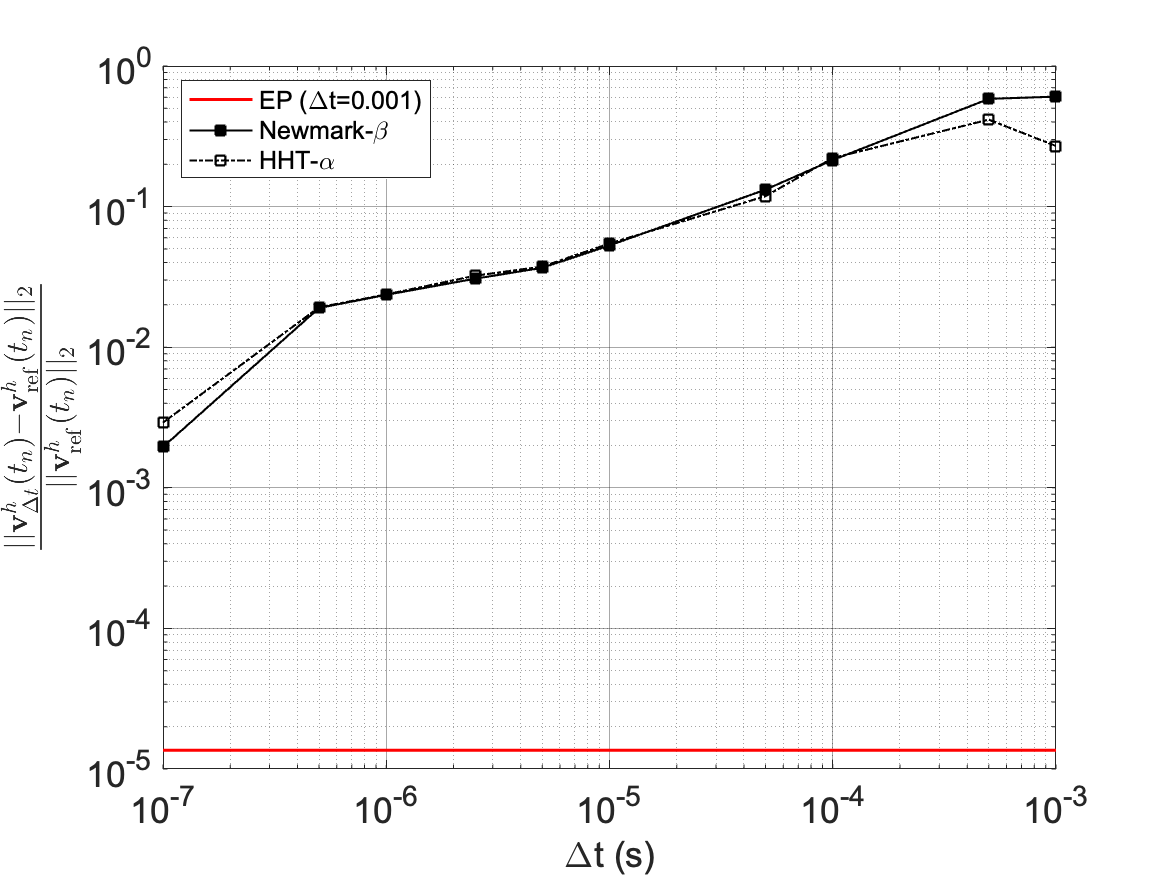}
         \caption{}
         \label{fig:linear numerical convergence - b}
     \end{subfigure} 
        \caption{Relative errors in displacement~(a) and velocity~(b) for the expotential propagator~(EP) method, the \Newmark method and \HHTA method for an elastodynamic simulation of a linear elastic structural steel cantilever beam for various time-step sizes~($\Delta t$), with errors measured at $t_n$~=~1.0~s.}
        \label{fig:linear numerical convergence}
\end{figure}
We note that the exponential propagator is exact in construction for linear elastodynamics problems. Thus, by using the exponential propagator, even numerically, the error incurred from time propagation can be systematically reduced (by using a larger subspace size) for any time-step size. However, as observed in Sec.~\ref{subsec:subspace characteristics}, in terms of computational efficiency, in order to realise maximum computational benefit, it serves best to choose the highest time-step size in the sublinear regime~($\Delta t = 0.001$). In Fig.~\ref{fig:linear numerical convergence}, we show the accuracy of exponential propagator and compare against the numerical convergence of \Newmark method and \HHTA method. It is evident that the exponential propagator is more accurate, in terms of displacement errors, even at a time-step size $10^4\times$ larger than the \Newmark and \HHTA methods. This advantage is even more pronounced in velocity errors. The superior accuracy places the exponential propagator as a more suitable choice for accurate \longtime simulations. The superior accuracy of the exponential propagator also translates to improved computational efficiency, as can be observed in Table \ref{table:linear computational efficiency}, due to the larger time-step sizes that can be utilized with this method.

\begin{table}[!htbp]
\centering
\scalebox{0.85}{

 \begin{tabular}{ | M{0.2\linewidth} | M{0.3\linewidth} | M{0.15\linewidth} | M{0.15\linewidth} | M{0.15\linewidth} | }
  \hline
  Method & $\Delta t$  & relative error displacement & relative error velocity & CPU-hrs  \\ 
  \hline \hline 
  
  Exponential Propagator & 0.001(subspace size = 200) & 1.22e-06 & 1.36e-05 & 7.16 \\

    \hline 
    
  Exponential Propagator & 0.0005(subspace size = 120) & 2.74e-07 & 8.90e-07 & 6.60 \\

  \hline
   \Newmark & 5e-06 & 1.30e-04 & 3.67e-02 & 22.16  \\
   \hline
   \Newmark & 1e-06 & 3.08e-05 & 2.36e-02 & 71.3 \\
   \hline
   \Newmark & 1e-07 & 1.52e-06 & 1.95e-03 & 467.95 \\
   \hline
   \HHTA & 5e-06 & 1.50e-04  & 3.74e-02 & 19.19  \\
   \hline
   \HHTA & 1e-06 & 3.80e-05 & 2.37e-02 & 91.55 \\
   \hline
   \HHTA & 1e-07 & 2.27e-06 & 2.91e-03 & 668.33 \\
  
  \hline
\end{tabular}}
\caption{Comparison of the relative errors and computational time~(CPU-hrs) at various time-step sizes~($\Delta t$) for the exponential propagator, \Newmark and \HHTA methods for a structural steel cantilever beam modelled using linear elastodynamics. The studies are run till a final time of  $t_n$~=~1.0~s, distributed across 16 processors.}
\label{table:linear computational efficiency}  
\end{table}


             
            \subsubsection{Nonlinear cantilever}\label{subsubsec:computational efficiency:nonlinear cantilever}


            For our first nonlinear example, we study the numerical convergence of a structural steel cantilever beam modeled using the St.~Venant-Kirchhoff model. The H-matrix for the St.~Venant-Kirchhoff model is derived in~\ref{subsec:appendix:SVK_derivation}. 
            In the benchmarking of this nonlinear system, in addition to the fully implicit \Newmark and \HHTA methods, we also compare their linearly implicit counterparts. We note that for the time-step sizes considered in the study, the linearly implicit counterparts give the direct integration schemes a computational advantage. For the same target accuracy, in comparison, the fully implicit schemes demonstrate better energy conservation and improved stability, however, they are computationally more expensive. We also include an energy-momentum conserving method~(EMC)~\cite{simo1992exact} in our studies to put into perspective the superior accuracy of the exponential propagator in \longtime simulations as will be discussed in the following sections.

    \begin{figure}[!htbp]
     \centering
     \begin{subfigure}[b]{1.0\textwidth}
         \centering
         \includegraphics[width=0.7\textwidth]{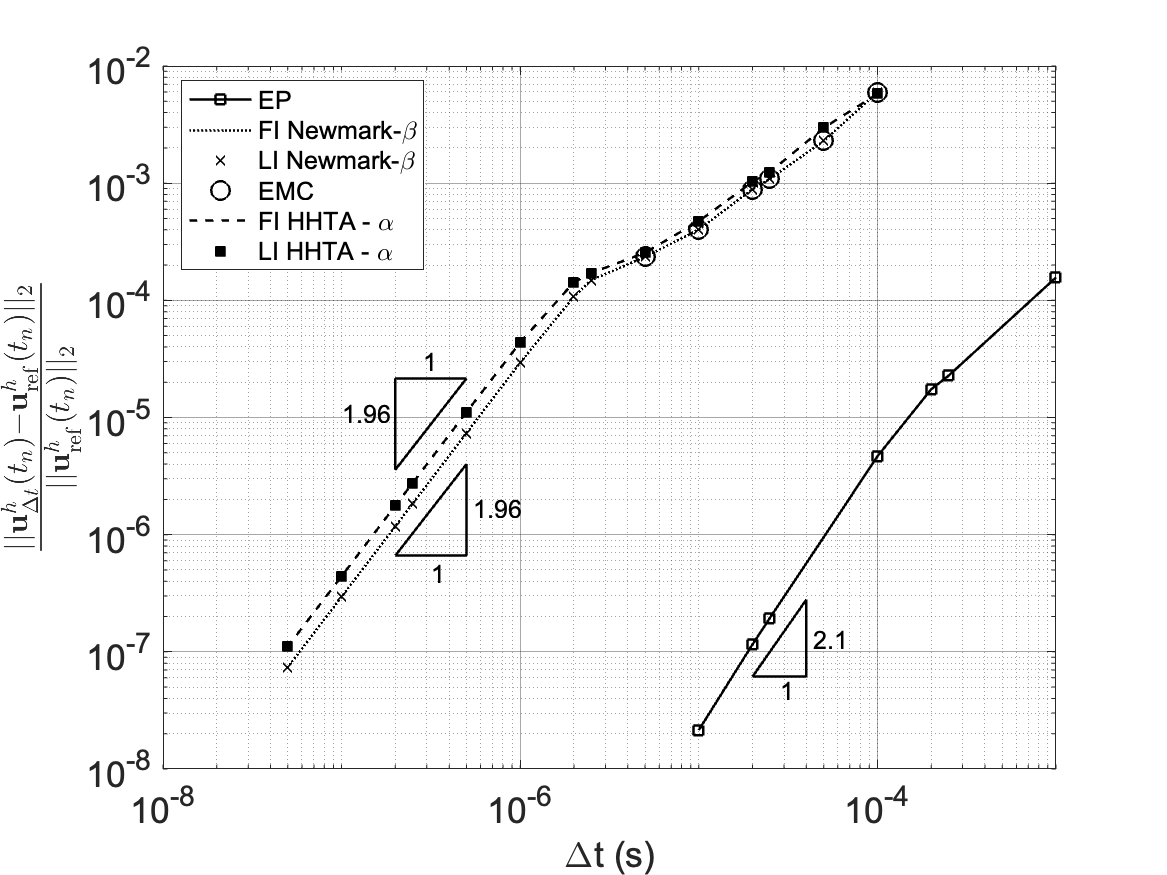}
         \caption{}
         \label{fig:nonlinear SVK numerical convergence - a}
     \end{subfigure}
     \begin{subfigure}[b]{1.0\textwidth}
         \centering
         \includegraphics[width=0.7\textwidth]{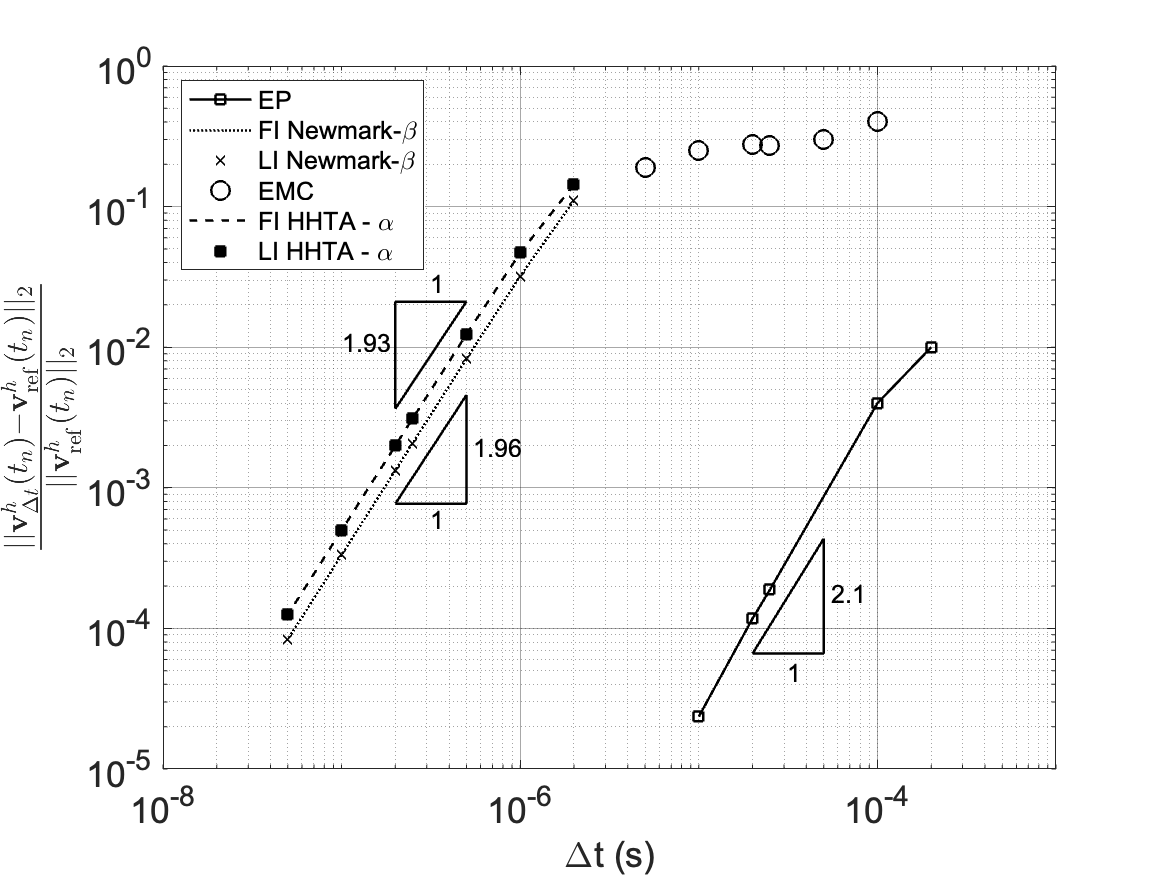}
         \caption{}
         \label{fig:nonlinear SVK numerical convergence - b}
     \end{subfigure} 
        \caption{Rate of convergence in displacement~(a) and velocity~(b) of the exponential propagator~(EP), fully implicit~(FI) \Newmark method, linearly implicit~(LI) \Newmark method, fully implicit~(FI) \HHTA method, linearly implicit~(LI) \HHTA method and EMC method for a structural steel cantilever beam modeled using St.~Venant-Kirchhoff model with errors measured at $t_n$~=~0.02~s.}
        \label{fig:nonlinear SVK numerical convergence}
\end{figure}

    \begin{table}[!htbp]
\centering
\begin{subtable}[t]{\textwidth}
\centering
\scalebox{0.85}{\begin{tabular}{ | M{0.4\linewidth} | M{0.3\linewidth} | M{0.3\linewidth} | }
  \hline
  Method & $\Delta t_{opt}$  & CPU-hrs  \\ 
  \hline \hline 
  
  Exponential  & 2.5e-04 & 0.14  \\
   Propagator & (subspace size 90) & \\
  \hline
  
  Fully implicit \Newmark & 1e-06 & 62.06 \\
  
  \hline
  
  Linearly implicit \Newmark & 1e-06 & 15.77  \\

  \hline

  Fully implicit \HHTA & 1e-06 & 63.95\\

  \hline

  Linearly implicit \HHTA  & 1e-06 & 17.02 \\

\hline
  EMC & 1e-06 & -\\
  
  \hline
\end{tabular}}
\caption{}
\label{table:SVK_T0p02 data - displacement}  
\end{subtable}
\begin{subtable}[t]{\textwidth}
\centering
    \scalebox{0.85}{

 \begin{tabular}{ | M{0.4\linewidth} | M{0.3\linewidth} | M{0.3\linewidth} | }
  \hline
  Method & $\Delta t_{opt}$  & CPU-hrs  \\ 
  \hline \hline 
  
  Exponential  & 1e-05 & 2.06  \\
   Propagator & (subspace size 40) & \\
  
  \hline
  
  Fully implicit \Newmark & 5e-08 & 2233.84  \\
  
  \hline
  
  Linearly implicit \Newmark & 5e-08 & 779.80  \\

  \hline

  Fully implicit \HHTA & 2.5e-08 & 4323.56*\\

  \hline

  Linearly implicit \HHTA & 2.5e-08 & 1337.89* \\

\hline
  EMC & 5e-08 & -\\
  \hline
\end{tabular}}
\caption{}
\label{table:SVK_T0p02 data - velocity}
\end{subtable}
\caption{Comparison of optimal time-step size~($\Delta t_{opt}$) and computational time~(CPU-hrs) required by the exponential propagator, fully implicit \Newmark method, linearly implicit \Newmark method, fully implicit \HHTA method, linearly implicit \HHTA method and EMC method for a structural steel cantilever beam modelled using St.~Venant-Kirchhoff model. The calculations are run till a final time of  $t_n$~=~0.02~s with a target accuracy~(relative error) of $10^{-4}$ in displacement~(a) and velocity~(b), distributed across 16 processors. `\,*\,' indicates the timings are projected estimations as the actual run till $t_n$~=~0.02~s was intractable. `\,-\,' indicates that the target accuracy~(relative error) was not achieved in any tractable computation.}
\label{table:SVK_T0p02 data}  
\end{table}

            Figure~\ref{fig:nonlinear SVK numerical convergence} shows the rates of convergence of the aforementioned methods with errors measured after one oscillation of the cantilever beam which is $\sim$ 0.02s. 
            These numerical results demonstrate the second-order convergence of the exponential propagator. More importantly, at each time-step size the exponential method is many orders of magnitude more accurate than the other methods studied in this benchmarking study. We remark that the exponential propagator is able to achieve a desired accuracy at time-step sizes that are $\sim100\times$ larger than the other methods. This translates to a $\sim100\times$ and $\sim400\times$ computational savings for a target accuracy~(relative error) of $10^{-4}$ in displacement and velocity, respectively, as evident from Table~\ref{table:SVK_T0p02 data - displacement} and Table~\ref{table:SVK_T0p02 data - velocity}.


            \subsubsection{Nonlinear soft tissue plate} \label{subsubsec:computational efficiency:nonlinear soft tissue elastomer} 

\begin{figure}[!htbp]
     \centering
     \begin{subfigure}[b]{1.0\textwidth}
         \centering
         \includegraphics[width=0.7\textwidth]{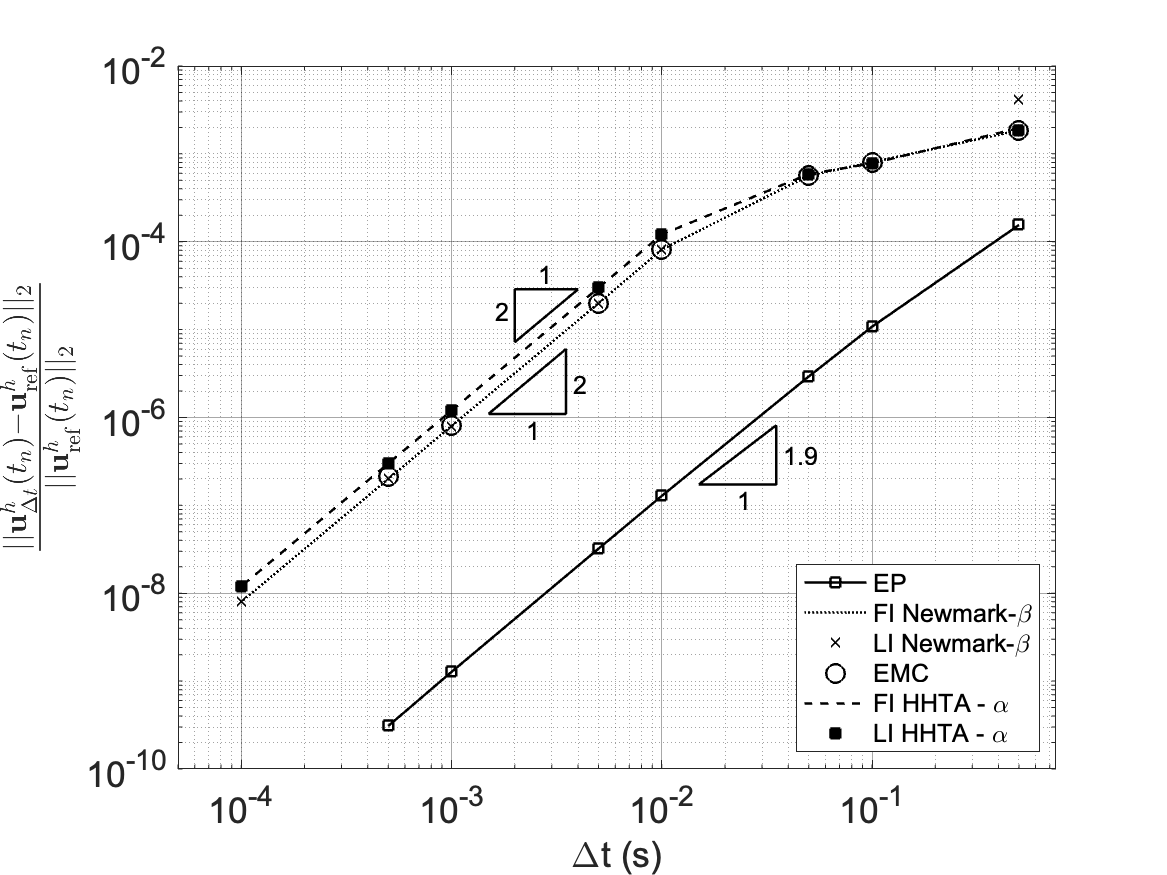}
         \caption{}
         \label{fig:nonlinear STE numerical convergence - a}
     \end{subfigure}
     \begin{subfigure}[b]{1.0\textwidth}
         \centering
         \includegraphics[width=0.7\textwidth]{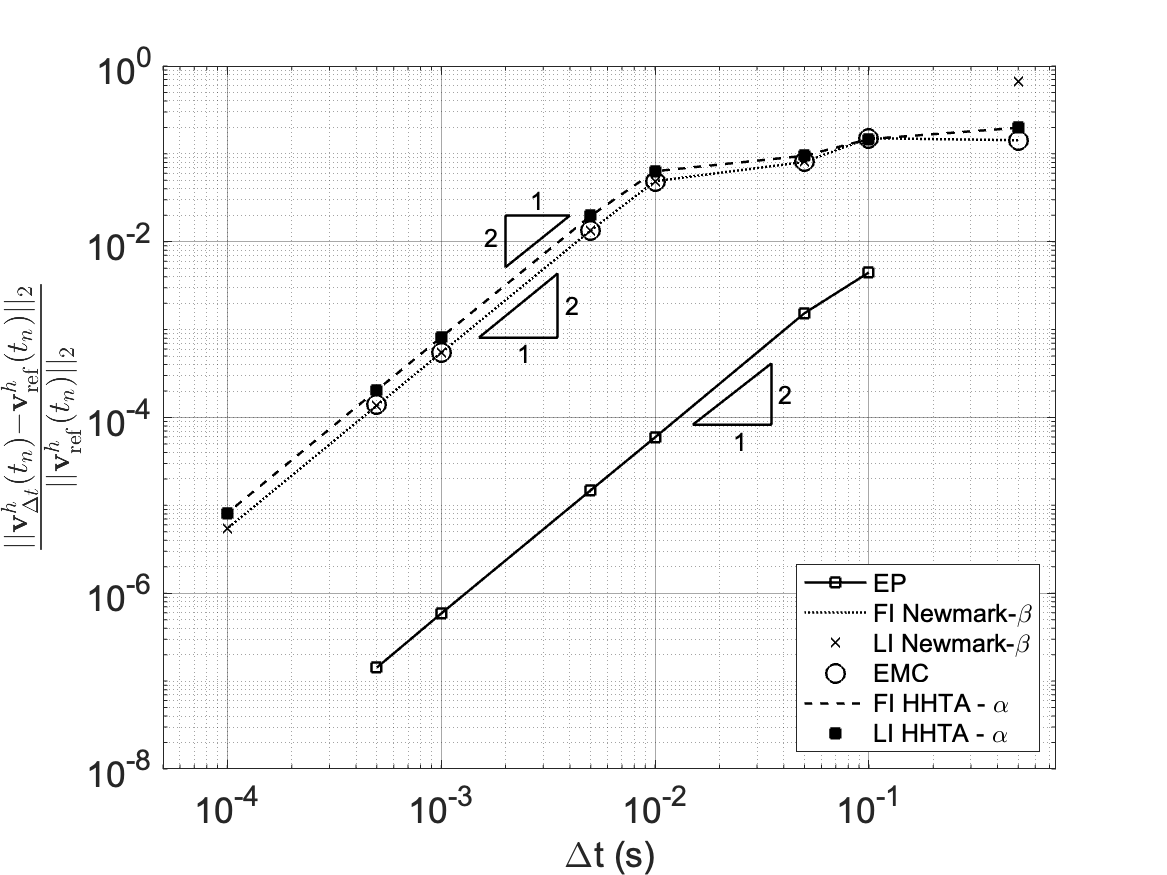}
         \caption{}
         \label{fig:nonlinear STE numerical convergence - b}
     \end{subfigure} 
        \caption{Rate of convergence in displacement~(a) and velocity~(b) of the exponential propagator~(EP), fully implicit~(FI) \Newmark method, linearly implicit~(LI) \Newmark method, fully implicit~(FI) \HHTA method, linearly implicit~(LI) \HHTA method and EMC method for a soft tissue plate modeled using Yeoh model with errors measured at $t_n$=1~s.}
        \label{fig:nonlinear STE numerical convergence}
\end{figure}

          We now present the results of the study of dynamics of a plate of soft tissue elastomer, which is one of the most widely modeled materials using hyperelasticity~\cite{STEParameters}. The parameters for the Yeoh model are obtained from~\cite{STEParameters}. The plate is subjected to an initial load of 10N applied as depicted in the Fig.~\ref{fig:Plate undeformed}, and the evolution in time is studied. We study the numerical convergence by measuring the errors at $t_n$=1~s, which are reported in Fig.~\ref{fig:nonlinear STE numerical convergence}. While the convergence rates are similar for the various methods, the exponential propagator has substantially lower errors. For a desired accuracy, this directly translates to staggering computational savings as demonstrated in Tables~\ref{table:STE data - displacement} and \ref{table:STE data - velocity} for a target accuracy~(relative error) of $10^{-4}$ in displacement and velocity, respectively.

\begin{table}[!htbp]
\centering
\begin{subtable}[t]{\textwidth}
\centering
\scalebox{0.85}{\begin{tabular}{ | M{0.4\linewidth} | M{0.3\linewidth} | M{0.3\linewidth} | }
  \hline
  Method & $\Delta t_{opt}$  & CPU-minutes  \\ 
  \hline \hline 
  
  Exponential  & 2.5e-01 & 2.43  \\
   Propagator & (subspace size 65) & \\
   \hline
  
  Fully implicit \Newmark & 1e-02 & 60.33 \\
  
  \hline
  
  Linearly implicit \Newmark & 1e-02 & 15.31  \\
  
  \hline
  
  Fully implicit \HHTA & 5e-03 & 88.08\\
  
  \hline
  
  Linearly implicit \HHTA & 5e-03 & 31.44 \\
 
 \hline
  EMC & 1e-02 & 440.99\\  
  \hline 
  
  \end{tabular}}
  \caption{}
\label{table:STE data - displacement}  
\end{subtable}
\begin{subtable}[t]{\textwidth}
\centering
\scalebox{0.85}{

 \begin{tabular}{ | M{0.4\linewidth} | M{0.3\linewidth} | M{0.3\linewidth} | }
  \hline
  Method & $\Delta t_{opt}$  & CPU-minutes  \\ 
  \hline \hline 
  
  Exponential  & 1e-02 & 26.41  \\
   Propagator & (subspace size 15) & \\
  
  \hline
  
  Fully implicit \Newmark & 2.5e-04 & 2281.19  \\
  
  \hline
  
  Linearly implicit \Newmark & 2.5e-04 & 627.06  \\

  \hline

  Fully implicit \HHTA & 2.5e-04 & 3672\\

  \hline

  Linearly implicit \HHTA & 2.5e-04 & 706.94 \\

\hline
  EMC & 2.5e-04 & 19469.2\\
  \hline
\end{tabular}}
\caption{}
\label{table:STE data - velocity}
\end{subtable}
\caption{Comparison of optimal time-step size~($\Delta t_{opt}$) and computational time~(CPU-minutes) required by the exponential propagator, fully implicit \Newmark method, linearly implicit \Newmark method, fully implicit \HHTA method, linearly implicit \HHTA method and EMC method for a soft tissue plate modeled using Yeoh model. The studies are run till a final time of  $t_n$~=~1~s with a target accuracy~(relative error) of $10^{-4}$ in displacement~(a) and velocity~(b), distributed across 16 processors.}
\label{table:STE data}  
\end{table}


\subsubsection{Nonlinear Rubber plate} \label{subsubsec:computational efficiency:nonlinear rubber}

\begin{figure}[!htbp]
     \centering
     \begin{subfigure}[b]{1.0\textwidth}
         \centering
         \includegraphics[width=0.7\textwidth]{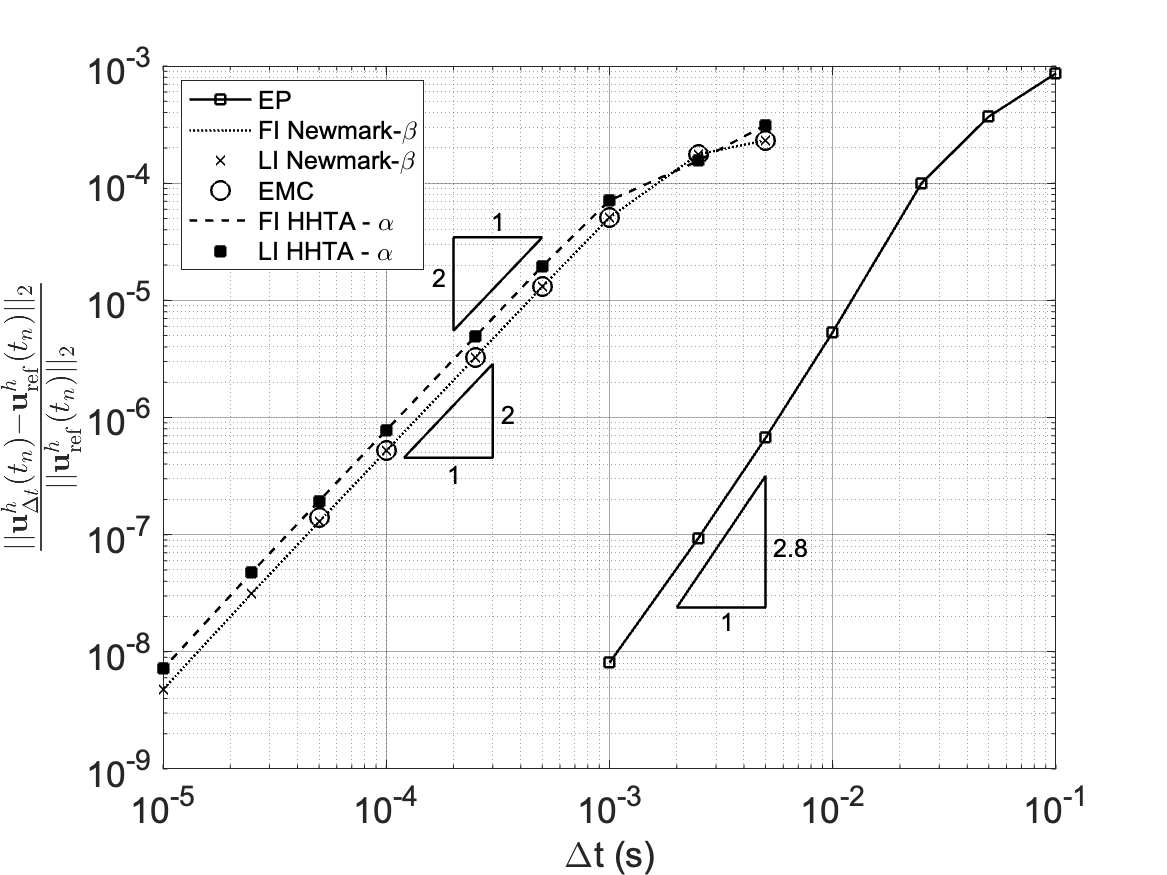}
         \caption{}
         \label{fig:nonlinear TR numerical convergence - a}
     \end{subfigure}
     \begin{subfigure}[b]{1.0\textwidth}
         \centering
         \includegraphics[width=0.7\textwidth]{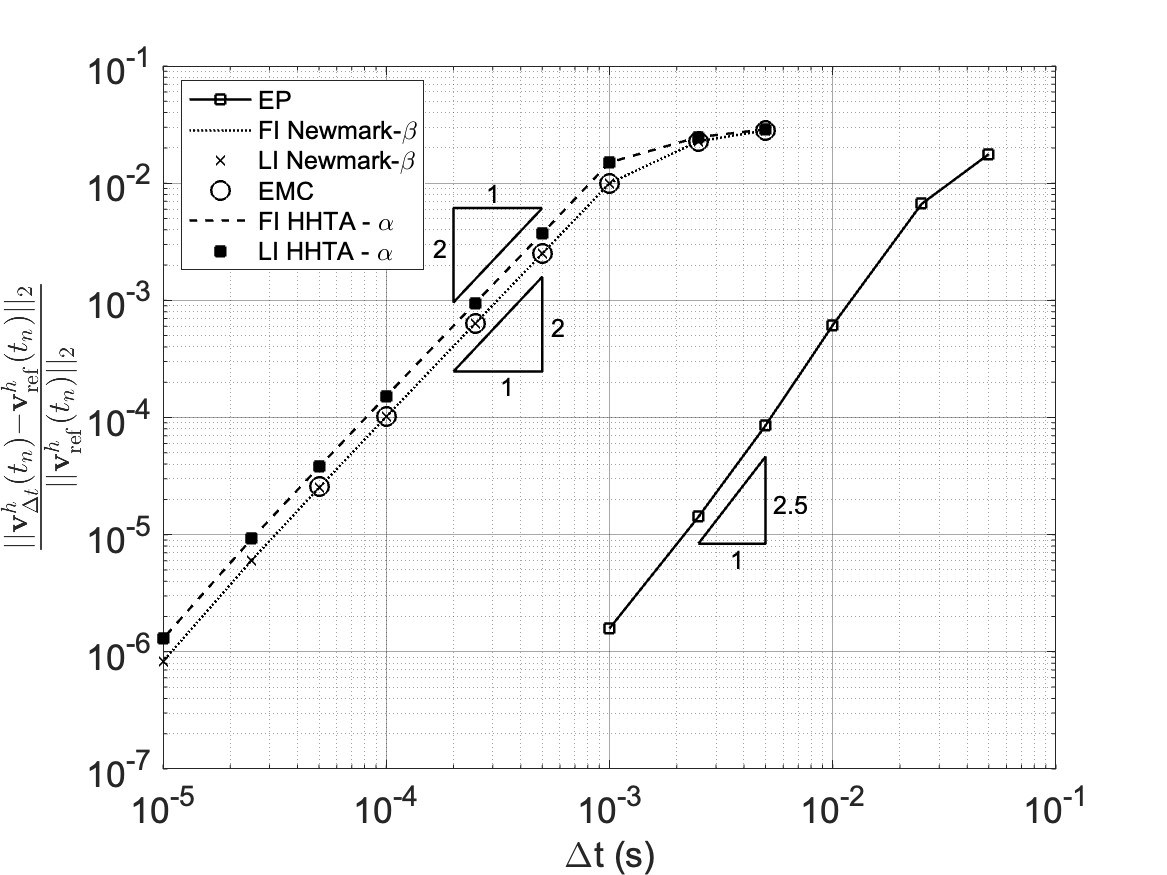}
         \caption{}
         \label{fig:nonlinear TR numerical convergence - b}
     \end{subfigure} 
        \caption{Rate of convergence in displacement~(a) and velocity~(b) of the exponential propagator~(EP), fully implicit~(FI) \Newmark method, linearly implicit~(LI) \Newmark method, fully implicit~(FI) \HHTA method, linearly implicit~(LI) \HHTA method and EMC method for a Treloar rubber plate modeled using Yeoh model with errors measured at $t_n$=1~s.}
        \label{fig:nonlinear TR numerical convergence}
\end{figure}

We also benchmark for another commonly modeled hyperelastic material, a plate of Treloar rubber. The parameters for the Yeoh model are obtained as $C_{10}=100 \text{ Pa}, C_{20}=-1\text{ Pa}, C_{30}=0.01\text{ Pa}, D_1=0.001\text{ (1/Pa)}$ from~\cite{AbaqusManual}. The plate is subjected to an initial load of  1N applied as depicted in the Fig.~\ref{fig:Plate undeformed}, and the evolution in time is studied. We study the numerical convergence by measuring the errors at $t_n$=1~s as depicted in Fig.~\ref{fig:nonlinear TR numerical convergence}. The improved accuracy afforded by the exponential propagator is clearly established, and this directly translates to substantial speed-ups in computational times (cf.~Tables \ref{table:TR data - displacement} and \ref{table:TR data - velocity}).

\begin{table}[!htbp]
\centering
\begin{subtable}[t]{\textwidth}
\centering
\scalebox{0.85}{

 \begin{tabular}{ | M{0.4\linewidth} | M{0.3\linewidth} | M{0.3\linewidth} | }
  \hline
  Method & $\Delta t_{opt}$  & CPU-minutes  \\ 
  \hline \hline 
  
  Exponential  & 0.01 & 29.20  \\
   Propagator & (subspace size 25) & \\
  \hline
  
  Fully implicit \Newmark & 1e-03 & 453.59 \\
  
  \hline
  
  Linearly implicit \Newmark & 1e-03 & 152.54  \\

  \hline

  Fully implicit \HHTA & 1e-03 & 459.89\\

  \hline

  Linearly implicit \HHTA & 1e-03 & 161.33 \\

\hline
  EMC & 1e-03 & 4228.99\\
  
  \hline
\end{tabular}}
\caption{}
\label{table:TR data - displacement}  
\end{subtable}
\begin{subtable}[t]{\textwidth}
\centering
\scalebox{0.85}{

 \begin{tabular}{ | M{0.4\linewidth} | M{0.3\linewidth} | M{0.3\linewidth} | }
  \hline
  Method & $\Delta t_{opt}$  & CPU-minutes  \\ 
  \hline \hline 
  
  Exponential  & 5e-03 & 54.17  \\
   Propagator & (subspace size 20) & \\
  
  \hline
  
  Fully implicit \Newmark & 5e-05 & 24331.12  \\
  
  \hline
  
  Linearly implicit \Newmark & 5e-05 & 5735.95  \\

  \hline

  Fully implicit \HHTA & 5e-05 & 25894.53\\

  \hline

  Linearly implicit \HHTA & 5e-05 & 6553.89 \\

\hline
  EMC & 5e-05 & 155182.13\\
  \hline
\end{tabular}}
\caption{}
\label{table:TR data - velocity}
\end{subtable}
\caption{Comparison of optimal time-step size~($\Delta t_{opt}$) and computational time~(CPU-minutes) required by the exponential propagator, fully implicit \Newmark method, linearly implicit \Newmark method, fully implicit \HHTA method, linearly implicit \HHTA method and EMC method for a Treloar rubber plate modeled using Yeoh model. The studies are run till a final time of  $t_n$~=~1~s with a target accuracy~(relative error) of $10^{-4}$ in displacement~(a) and velocity~(b), distributed across 16 processors.}
\label{table:TR data}  
\end{table}


        \subsection{Remarks on \longtime behaviour}\label{subsec:remarks on \longtime behaviour}

    \begin{figure}[!htbp]
     \begin{subfigure}[b]{1.0\textwidth}
         \centering
         \includegraphics[width=0.7\textwidth]{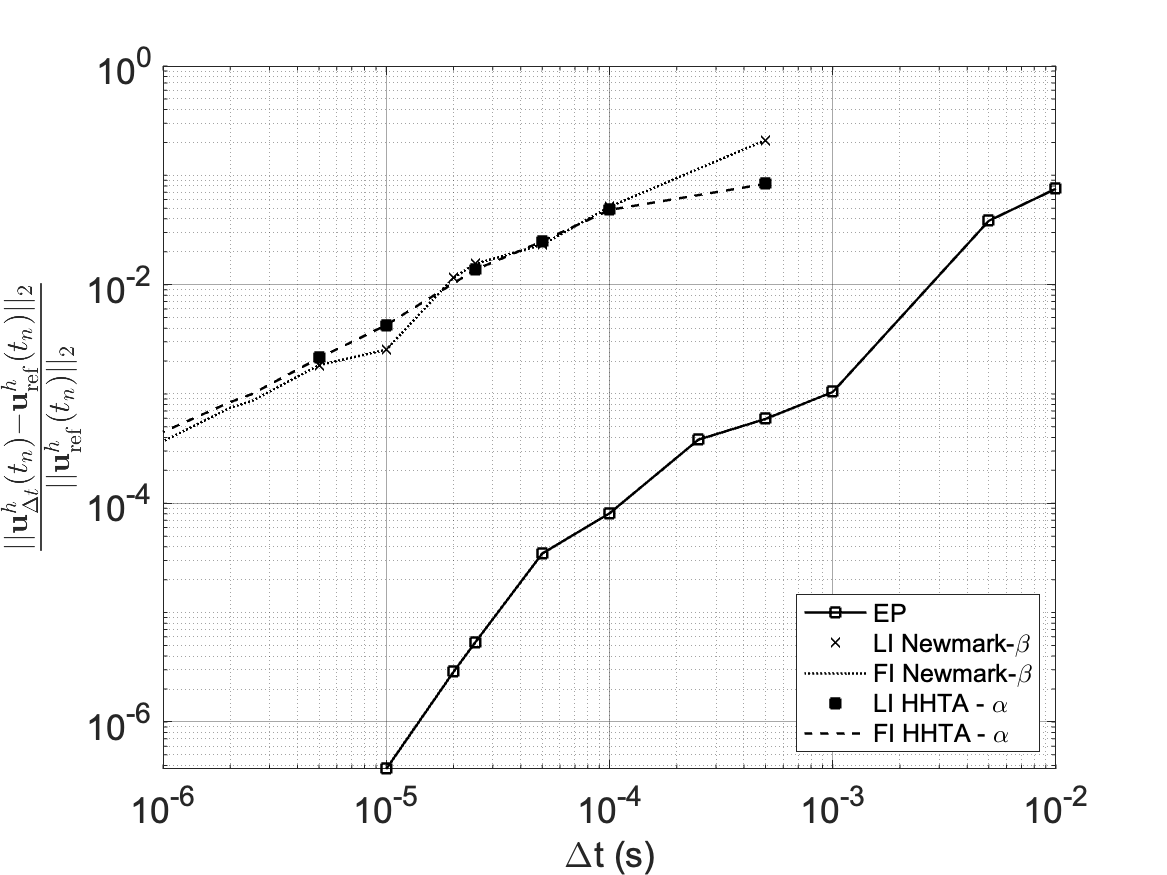}
         \caption{}
         \label{fig:\longtime behaviour numerical convergence - a}
     \end{subfigure}
     \begin{subfigure}[b]{1.0\textwidth}
         \centering
         \includegraphics[width=0.7\textwidth]{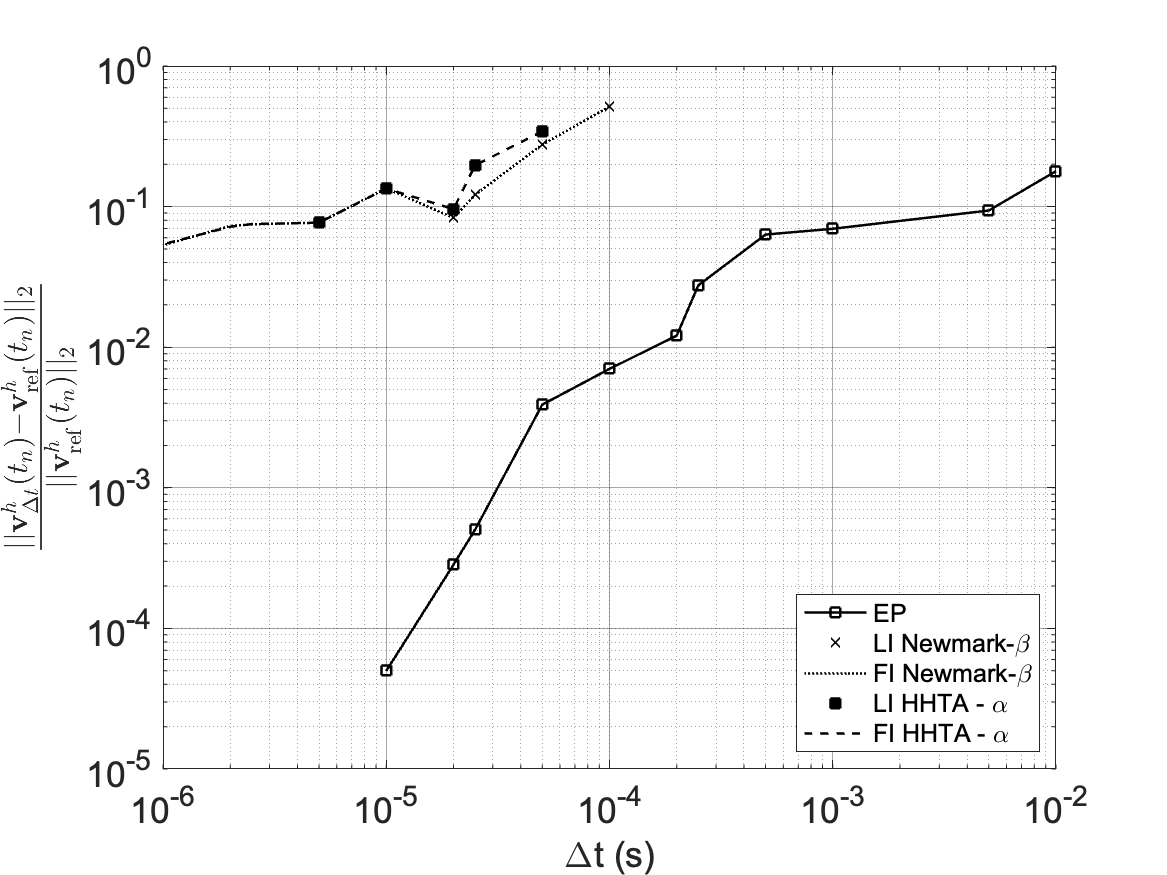}
         \caption{}
         \label{fig:\longtime behaviour numerical convergence - b}
     \end{subfigure} 
        \caption{Rate of convergence in displacement~(a) and velocity~(b) of the exponential propagator~(EP), fully implicit~(FI) \Newmark method, linearly implicit~(LI) \Newmark method, fully implicit~(FI) \HHTA method, linearly implicit~(LI) \HHTA method for a structural steel cantilever beam modeled using St.~Venant-Kirchhoff model with errors measured at $t_n$=1~s.}
        \label{fig:\longtime behaviour numerical convergence}
\end{figure}

%
\begin{figure}[!htbp]
   \begin{subfigure}[b]{1.0\textwidth}
         \centering
         \includegraphics[width=0.7\textwidth]{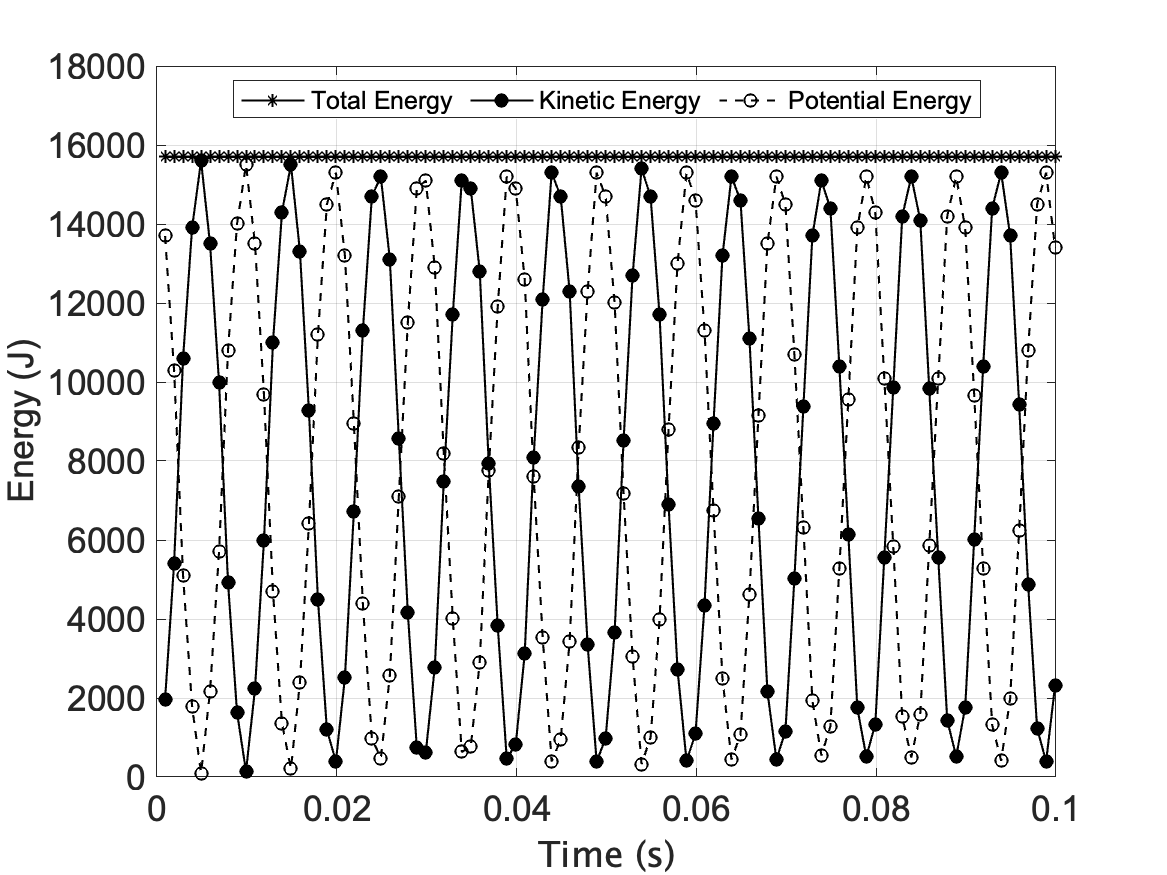}
     \end{subfigure}
     \begin{subfigure}[b]{1.0\textwidth}
         \centering
         \includegraphics[width=0.7\textwidth]{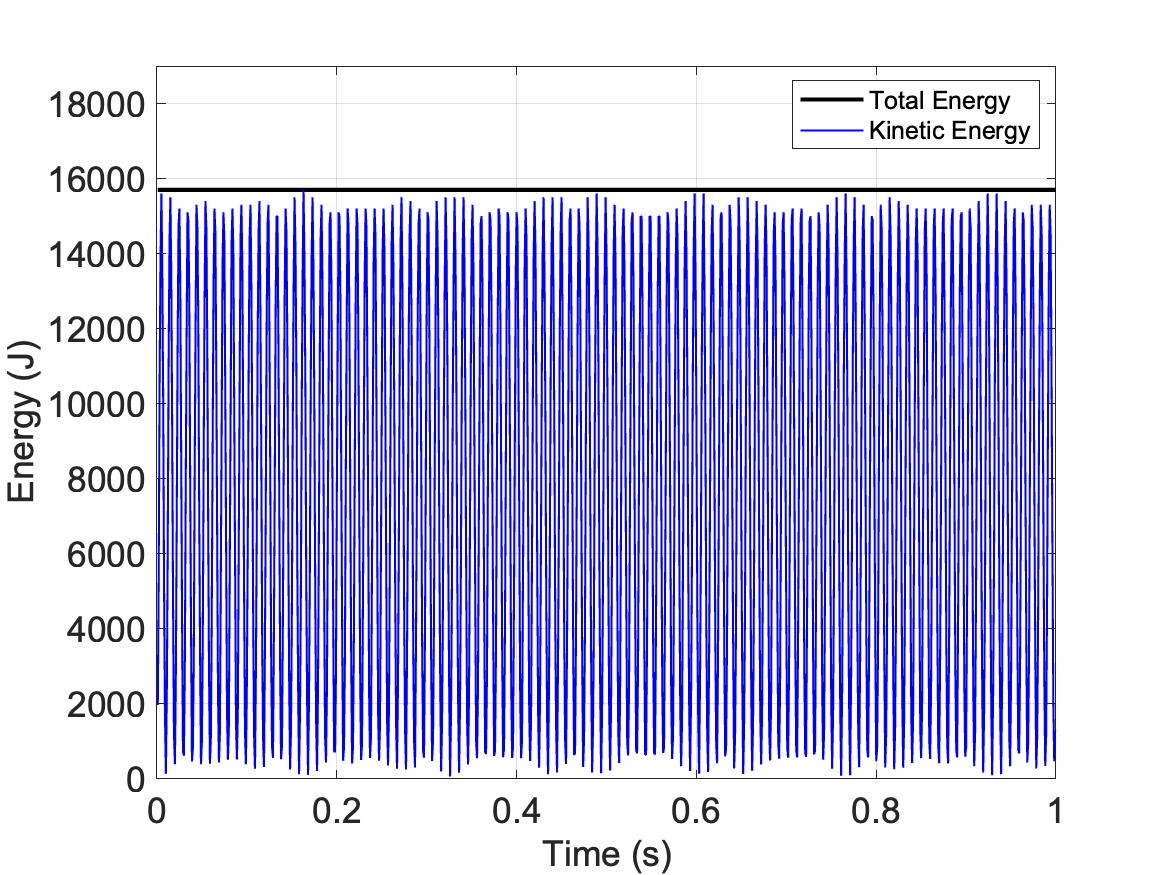}
     \end{subfigure} 
        \caption{Time histories of the total energy, kinetic energy and potential energy computed using the exponential propagator for a structural steel cantilever beam modeled using the St.~Venant-Kirchhoff model.}
        \label{fig:\longtime energy behaviour}
\end{figure}

In all of our benchmark studies, we observe a consistent out-performance by the exponential propagator in terms of accuracy and computational time. An important next question to address is if these results are extendable to \longtime simulations. To this end, we revist our study on the nonlinear cantilever beam in section \ref{subsubsec:computational efficiency:nonlinear cantilever}. We again study the numerical convergence as shown in Fig.~\ref{fig:\longtime behaviour numerical convergence} but till a final time of $t_n$=1~s by which the beam completes $\sim 50$ oscillations. As observed from the Fig. \ref{fig:\longtime behaviour numerical convergence}, the exponential propagator again demonstrates much higher accuracy at any given time-step size. By comparing these results against the numerical convergence in Fig.~\ref{fig:nonlinear SVK numerical convergence} studied after one oscillation, we observe that the advantage realized by the exponential propagator is further increased in a \longtime simulation. Another important observation is that for this benchmark problem, it is practically intractable to achieve a relative accuracy below $10^{-2}$ in velocity with other approaches. On the contrary, the exponential propagator is trivially capable of much higher accuracies in velocity at higher time-step sizes.  We also observe from the total energy plot depicted in Fig.~\ref{fig:\longtime energy behaviour} that the exponential propagator demonstrates a stable energy profile for a \longtime calculation. This stable energy behaviour in combination with the superior accuracy clearly establishes the case for exponential propagators as a more suitable choice for \longtime simulations. The speed-ups in computational time are reported in Tables \ref{table:SVK_T1 data - displacement} and \ref{table:SVK_T1 data - velocity}.

\begin{table}[!htbp]
\centering
\begin{subtable}[t]{\textwidth}
\centering
\scalebox{0.85}{

 \begin{tabular}{ | M{0.4\linewidth} | M{0.3\linewidth} | M{0.3\linewidth} | }
  \hline
  Method & $\Delta t_{opt}$  & CPU-hrs  \\ 
  \hline \hline 
  
  Exponential  & 5e-04 & 7.46  \\
   Propagator & (subspace size 140) & \\
  \hline
  
  Fully implicit \Newmark & 2.5e-06 & 2208.80* \\
  
  \hline
  
  Linearly implicit \Newmark & 2.5e-06 & 655.39  \\

  \hline

  Fully implicit \HHTA & 2.5e-06 & 1961.54*\\

  \hline

  Linearly implicit \HHTA & 2.5e-06 & 683.32 \\

  \hline
\end{tabular}}
\caption{}
\label{table:SVK_T1 data - displacement}  
\end{subtable}

\begin{subtable}[t]{\textwidth}
\centering
\scalebox{0.85}{

 \begin{tabular}{ | M{0.4\linewidth} | M{0.3\linewidth} | M{0.3\linewidth} | }
  \hline
  Method & $\Delta t_{opt}$  & CPU-hrs  \\ 
  \hline \hline 
  
  Exponential  & 1e-04 & 30.07  \\
   Propagator & (subspace size 70) & \\
  
  \hline
  
  Fully implicit \Newmark & - & -  \\
  
  \hline
  
  Linearly implicit \Newmark & - & -  \\

  \hline

  Fully implicit \HHTA & - & -\\

  \hline

  Linearly implicit \HHTA & - & - \\

  \hline
\end{tabular}}
\caption{}
\label{table:SVK_T1 data - velocity}
\end{subtable}
\caption{Comparison of optimal time-step size~($\Delta t_{opt}$) and computational time~(CPU-hrs) required by the exponential propagator, fully implicit \Newmark method, linearly implicit \Newmark method, fully implicit \HHTA method and linearly implicit \HHTA method for a structural steel cantilever beam modeled using St.~Venant-Kirchhoff model. The studies are run till a final time of  $t_n$~=~1~s with a target accuracy~(relative error) of $10^{-3}$ in displacement~(a) and a target accuracy~(relative error) of $10^{-2}$ in velocity~(b), distributed across 16 processors. `\,*\,' indicates the timings are projected estimations as the actual run till $t_n$~=~1~s was intractable. `\,-\,' indicates that the target accuracy~(relative error) was not achieved in any tractable time-step size.}
\label{table:SVK_T1 data}  
\end{table}

        \subsection{Parallel scalability}\label{subsec:parallel scalability}

\begin{figure}[h]
    \centering
    \includegraphics[width=0.7\textwidth]{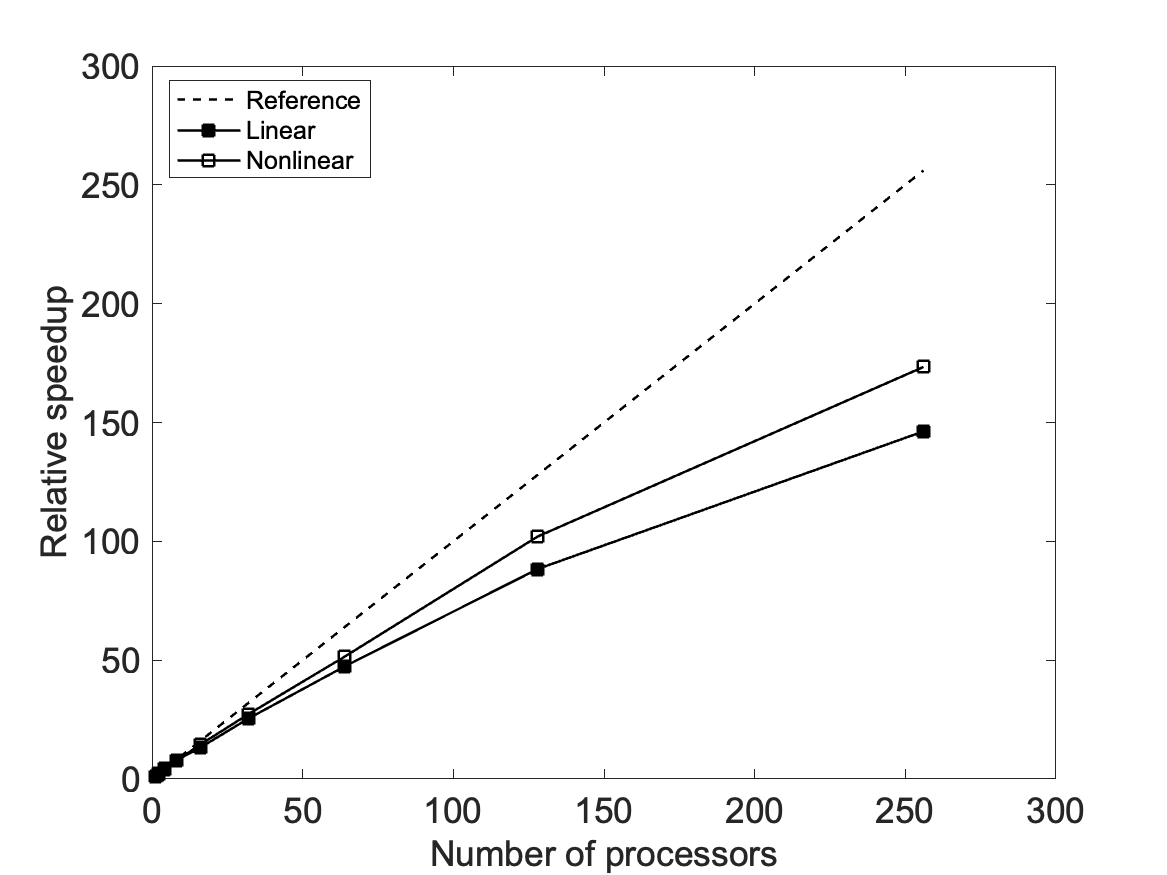}
    \caption{Parallel scalability of the exponential time propagator.}
    \label{fig:Parallel scalability}
\end{figure}

We now demonstrate the parallel scalability of the exponential time propagator in Fig.~\ref{fig:Parallel scalability}. In our study, we choose a cantilever beam of dimensions 4m$\times$1m$\times$1m with a mesh-size of 0.0625m containing $\sim4 \times 10^5$ degrees of freedom. 
We study the relative speedup of the propagator in the linear and nonlinear problems, with the nonlinear problem modeled using the St.~Venant-Kirchhoff model. In both the cases, we measure the wall time spent by the exponential propagator to run 100 time-steps at a subspace size of 50. The relative speedups are measured against the wall time required by a single processor calculation.

From the Fig.~\ref{fig:Parallel scalability}, we observe that the linear problem demonstrates good scalability till 128 processors with an efficiency of 70\%. At 128 processors, the numer of degrees of freedom per processor is $\sim$3,200, which is on the lower end to provide good parallel scalability with further parallelization. The nonlinear problem also exhibits good parallel scalability till 128 processors, with a better parallel efficiency of 80\%. The difference in the parallel efficiencies can be attributed to the fact that the linear problem is dominated by matrix-vector multiplications, which also involves communication. On the other hand, for the nonlinear problem, matrix constructions, which are embarrasingly parallel, also constitute a significant part of the calculation, thus providing an overall better parallel efficiency.


\section{Conclusion and future work} \label{sec:conclusion and future work}


In summary, we have developed an accurate, computationally efficient, and scalable exponential propagator for elastodynamics. In formulating the exponential propagator, we recast the internal force vector, discretized on a finite element basis, in terms of a H-matrix acting on the discrete displacement vector. This, in turn, allowed for a reformulation of the second-order dynamical system into an equivalent first-order system. We expressed the solution to the first-order initial value problem in terms of a Magnus expansion. Specifically, we truncated the Magnus series after the first term and evaluated it using the mid-point rule. A predictor-corrector technique was employed to evaluate the H-matrix at the midpoint of the time-step. Through error analysis, we established second-order convergence for the proposed propagator. We also proved unconditional stability and energy conservation for linear elastodynamics and demonstrated symplecity for both linear and nonlinear elastodynamics. To accurately and efficiently evaluate the action of the exponential propagator on a vector, we employed an adaptive Krylov subspace-based projection, where the optimal subspace size is determined via \textit{apriori} error estimate. 

We demonstrated the accuracy and performance of the proposed exponential propagator using both linear and nonlinear benchmark systems. Numerically, we obtained close to second-order convergence, corroborating our error estimates. In the case of linear and nonlinear elastodynamics, we numerically demonstrated good energy conservation. We established the computational advantage in using the Krylov subspace projection and also provided guidance on how to choose an optimal combination of the subspace size and time-step size to realize maximum computational advantage. We observed stable dynamics in nonlinear problems even in the long time-scale regime. In all the benchmark systems, the exponential propagator attains the required accuracy at orders of magnitude larger time-steps as compared to the widely used Newmark-$\beta$ and HHT-$\alpha$ methods, amounting to a staggering $10-100\times$ speedup in computational time. The benefits were more pronounced for \longtime scale simulations where the exponential propagator gave orders of magnitude more accurate solutions at much lower computational costs compared to Newmark-$\beta$ and HHT-$\alpha$. In particular, in the \longtime scale simulations, the accuracies achieved through the exponential propagator were intractable for the conventional methods even with more refined temporal discretizations.

Among the various approximations involved in the exponential propagator, we observe that the dominant error arises from the truncation of the Magnus series. In our energy studies, we noticed that the minor deviations from energy conservation were directly attributable to the error accrued by Magnus truncation. However, one can construct higher-order propagators by including more terms in the Magnus expansion. These high-order propagators can simultaneously boost the accuracy and efficiency over the proposed second-order propagator, and will be a subject of our future work. 
This can improve the qualitative characteristics of the solution while maintaining the same superior accuracy without any significant addition to the computational cost. The exponential propagator is a good substitute for conventional methods, in particular, to perform \longtime stable simulations in problems arising in elastodynamics. 
The methods developed here are not however restricted to problems only in elastodynamics and can be extended to any dynamical system, such as imaginary-time propagation in quantum mechanics~\cite{Davies1980, Lehtovaara2007, Bader2013}, molecular dynamics \cite{nettesheim1996explicit, hochbruck1999bunch, hochbruck1999exponential}, Maxwell equations \cite{botchev2006gautschi,botchev2009application,verwer2009unconditionally}, advection-diffusion equations \cite{friesner1989method, bergamaschi2004relpm, caliari2004interpolating}, Navier Stokes equations\cite{schulze2009exponential, saad1991application, edwards1994krylov, newman2004exponential}, to list a few.

\section{Acknowledgements}
We thank Michigan Institute for Computational Discovery in Engineering at the University of Michigan for supporting this work. V.G. also gratefully acknowledges the support of the Army Research Office through the DURIP grant W911NF1810242, which also provided computational resources for this work.

\newpage
\appendix

\newcommand{\xt}{\left(\mathbf{X},t\right)}
\newcommand{\xo}{\left(\mathbf{X},0\right)}
\newcommand{\x}{\mathbf{X}}
\newcommand{\gradX}{\grad_\mathbf{X}}
\newcommand{\domain}{\Omega_0}
\newcommand{\elemDomain}{\Omega_e}
\newcommand{\totalB}{\partial \Omega_0}
\newcommand{\DirichletB}{\partial \Omega_0^D}
\newcommand{\NeumannB}{\partial \Omega_0^N}
\renewcommand{\u}{\bm{u}}
\renewcommand{\v}{\bm{v}}
\renewcommand{\w}{\bm{w}}
\newcommand{\mbu}{\mathbf{u}}
\newcommand{\mbv}{\mathbf{v}}
\newcommand{\mbw}{\mathbf{w}}
\newcommand{\numElem}{n_e}
\newcommand{\elemUnion}{\bigcup_{e=1}^{\numElem}}
\newcommand{\shapeFuncSumI}{\sum\limits_{I=1}^{n}}
\newcommand{\shapeFuncSumJ}{\sum\limits_{J=1}^{n}}
\newcommand{\E}{\mathbf{E}}
\newcommand{\F}{\mathbf{F}}
\newcommand{\I}{\mathbf{I}}
\renewcommand{\S}{\mathbf{S}}
\renewcommand{\tr}{\text{tr}}
\newcommand{\wc}{\text{$\bm{w}$}}
\newcommand{\wAe}{\mbw_{A,e}}
\newcommand{\wBe}{\mbw_{B,e}}
\newcommand{\uAe}{\mbu_{A,e}}
\newcommand{\uBe}{\mbu_{B,e}}
\renewcommand{\D}{\mathbf{D}}
\renewcommand{\P}{\mathbf{P}}
\newcommand{\C}{\mathbf{C}}
\newcommand{\detF}{\det(\mathbf{F}_e)}
\newcommand{\detFtt}{(\det(\mathbf{F}_e))^{-\frac{2}{3}}}
\newcommand{\trFTF}{\text{tr}(\F_e^T\F_e)}
\newcommand{\Cinv}{\F_e^{-1}\F_e^{-T}}
\newcommand{\Finve}{\Finv_e}
\newcommand{\FinvTe}{\FinvT_e}

\section{\texorpdfstring{$\mathbf{H}(u)$}{H(u)} construction}\label{sec:appendix:H formulation}
 In this work, we study the dynamics of benchmark systems in linear and nonlinear hyperelasticity. In this section, we provide the formulation of H-matrix~($\mathbf{H(u)}$) for the St.~Venant-Kirchhoff and Yeoh models. These formulations can be trivially extended to the Neo-Hookean and Mooney-Rivlin models. The notation in this section is consistent with the notation used in section~\ref{sec:formulation}.
 %
\subsection{St.~Venant-Kirchhoff model}\label{subsec:appendix:SVK_derivation}
St.~Venant-Kirchhoff model is a nonlinear system exhibiting purely geometric nonlinearity. The strain energy function for the St.~Venant-Kirchhoff model is given as follows
\begin{equation}
    \Psi(\mathbf{E}) = \frac{\lambda}{2} [\text{tr}(\mathbf{E})]^2 + \mu\text{tr}(\mathbf{E}^2)
\end{equation}
where $\lambda$,$\mu$ are the Lam\'e constants, tr(.) denotes the trace of the matrix. Here $\mathbf{E}$ is the Green-Lagrange strain and is defined as
\begin{equation}\label{eq: Green-Lagrange strain and deformation gradient}
     \E = \frac{\F^T\F - \I}{2}\text{; } \F = \I + \mathbf{\grad_\mathbf{X}\text{$\bm{u}$}}
\end{equation}
where $\F$ is the deformation gradient, $\bm{u}$ is the displacement vector, $\I$ is the identity tensor, $\mathbf{X}$ is any arbitrary point in the reference configuration, and $\mathbf{\grad_\mathbf{X}}$ is the gradient operator in the reference configuration.
The corresponding second Piola-Kirchoff stress tensor expressed in the Lagrangian framework is given as 
\begin{equation}
    \S = \frac{\partial{\Psi(\E)}}{\partial{\E}} = \lambda(\tr(\E))\I + 2\mu \E \, .
\end{equation}
The internal virtual work of the system is given as
\begin{equation}\label{eq:internal virtual work - SVK}
    \int_{\Omega_0} \delta\E:\S \, dV = \int_{\Omega_0}\left[\frac{1}{2}\brac{\F^T(\mathbf{\grad_X\wc}) + (\mathbf{\grad_X\wc})^T\F}\right]:\S \,dV
\end{equation}
where $\Omega_0$ represents the spatial domain in the reference configuration, $dV$ is the differential volume element, $\bm{w}(\mathbf{X})$ is the test function, and $\delta \E$ represents the directional derivative tensor of the Green-Lagrange strain~($\E$) along the direction of $\bm{w}$. In the finite element framework, the spatial domain~($\Omega_0$) is discretized into $n_e$ non-overlapping  subdomains($\Omega_e$), and equation \eqref{eq:internal virtual work - SVK} is correspondingly written as
\begin{equation}\label{eq:internal virtual work discretized - SVK}
\begin{split}
    \int_{\Omega_0} \delta\E \colon \S dV & =  \bigcup\limits_{e=1}^{n_e} \int_{\Omega_e} \frac{1}{2} \brac{\F_e^T\brac{\mathbf{\grad_X}\wc_e} + \brac{\mathbf{\grad_X}\wc_e}^T\F_e} \colon \brac{\frac{\lambda}{2}\brac{\text{tr}\brac{\F_e^T\F_e - \I}}\I + \mu\brac{\F_e^T\F_e - \I}} \, dV_e \\
    \end{split}
\end{equation}
where $dV_e$, $\F_e$, $\wc_e$ are the differential volume, deformation gradient, test function, respectively, in the finite element~($\Omega_e$). The corresponding interpolation within each finite element~($\Omega_e$) is given by
\begin{equation}\label{eq: u,w,X discretization}
    \begin{split}
        \u_e(\xi) & \approx \u^h_e(\xi) = \sum_{A=1}^{n}N_A(\xi) \mbu_{A,e}\,, \\
        \w_e(\xi) & \approx \w^h_e(\xi) = \sum_{A=1}^{n}N_A(\xi) \mbw_{A,e}\,, \\
        \x_e(\xi) & \approx \x^h_e(\xi) = \sum_{A=1}^{n}N_A(\xi) \x_{A,e}\,, \\
    \end{split}
\end{equation}
where $\xi$ is an arbitrary point within the element $\elemDomain$; $N_A(\xi)$ represents the value of the $A^{th}$ shape function at $\xi$; $\u_e^h, \w_e^h, \x_e^h$ are the elemental-level approximations of the fields $\u, \w, \x$ in the finite element space; $\mathbf{u}_{A,e}, \mathbf{v}_{A,e}, \mathbf{x}_{A,e}$ represent the $n$ nodal  values of $\u, \v, \x$, respectively. Correspondingly, $\F_e(\xi)$ and $\gradX\wc_e(\xi)$ are expressed as
\begin{align}
        \F_e(\xi) & \approx \I + \sum_{A=1}^n \mbu_{A,e}\otimes \gradX N_A(\xi)\,, \label{eq:Fe in FE space} \\
        \gradX \wc_e(\xi) & \approx \sum_{A=1}^n \mbw_{A,e}\otimes \gradX N_A(\xi)\,. \label{eq:gradX in FE space}
\end{align}
Using equations~\eqref{eq:Fe in FE space}, \eqref{eq:gradX in FE space}, the internal virtual work in equation~\eqref{eq:internal virtual work discretized - SVK} can be approximated in the finite element space as
\begin{align}\label{eq:internal virtual work in finite-element space - SVK}
          \int_{\Omega_0} \delta\E \colon \S dV & =  \bigcup\limits_{e=1}^{n_e} \int_{\Omega_e} \frac{1}{2} \brac{\F_e^T\brac{\mathbf{\grad_X}\wc_e} + \brac{\mathbf{\grad_X}\wc_e}^T\F_e} \colon \brac{\frac{\lambda}{2}\brac{\text{tr}\brac{\F_e^T\F_e - \I}}\I + \mu\brac{\F_e^T\F_e - \I}} \, dV_e \nonumber \\
          & := \Bigg[ \frac{1}{2}\sum_A \Bigg( \brac{\wAe \otimes \gNa} + \brac{\gNa \otimes \wAe} + \sum_B \Big( \brac{\brac{\gNb \otimes \uBe}\brac{\wAe \otimes \gNa}} \nonumber \\
          & \qquad + \brac{\brac{\gNa \otimes \wAe}\brac{\uBe \otimes \gNb}} \Big) \Bigg) \Bigg] \colon \Bigg[ \sum_A\Bigg( \lambda \brac{\gNa \cdot \uAe} + \mu \big(\uAe \otimes \gNa \nonumber \\
          & \qquad + \gNa \otimes \uAe\big) + \sum_B \Big(\frac{\lambda}{2}\brac{\gNa \cdot \gNb}\brac{\uBe \cdot \uAe} \nonumber \\
          & \qquad + \mu \brac{\gNa \otimes \gNb} \brac{\uAe \cdot \uBe} \Big) \Bigg) \Bigg]  \,dV_e \,.
\end{align}
Let us define the following for ease of representation:
\begin{align}
 \D_e & = \sum_A \uAe \otimes \gNa \,, \label{eq: D_e}\\
\S_{\text{SVK},e} & = \sum_A\Bigg( \lambda \brac{\gNa \cdot \uAe} + \mu \brac{\uAe \otimes \gNa + \gNa \otimes \uAe} \nonumber \\
          & \qquad +  \sum_B \Big(\frac{\lambda}{2}\brac{\gNa \cdot \gNb}\brac{\uBe \cdot \uAe} + \mu \brac{\gNa \otimes \gNb} \brac{\uAe \cdot \uBe} \Big) \Bigg)\,. \label{eq:SVKSExpression}
\end{align}
Applying the contraction on equation~\eqref{eq:internal virtual work in finite-element space - SVK} and simplying the subsequent expression using equations~\eqref{eq: D_e} and \eqref{eq:SVKSExpression}, we get the following
\begin{align}
        \int_{\Omega_0} \delta\E \colon \S dV & := \bigcup\limits_{e=1}^{n_e} \int_{\Omega_e} \frac{1}{2}\Bigg[\sum_A \sum_B \wAe \Bigg( \lambda \brac{\gNa \otimes \gNb} \brac{2\I + \D_e^T} + 2\mu \brac{\gNb \cdot \gNa} \nonumber \\ 
        & \quad + \mu \brac{\gNb \otimes \gNa} \brac{2\I + \D_e^T} + \mu \D_e^T \brac{\gNa \cdot \gNb} + 2 \brac{\gNb \cdot \S_{\text{SVK},e} \cdot \gNa}\Bigg) \uBe \Bigg] \,dV_e \nonumber\\
        & = \mbw^T \Bigg[\bigcup\limits_{e=1}^{n_e} \int_{\Omega_e} \frac{1}{2}\sum_A \sum_B \Bigg( \lambda \brac{\gNa \otimes \gNb}\brac{2\I + \D_e^T} + 2\mu \brac{\gNb \cdot \gNa} \nonumber \\ 
        & \quad + \mu \brac{\gNb \otimes \gNa}\brac{2\I + \D_e^T} + \mu \D_e^T \brac{\gNa \cdot \gNb}  + 2 \brac{\gNb \cdot \S_{\text{SVK},e} \cdot \gNa}\Bigg) \,dV_e \Bigg] \mbu \nonumber\\
        & = \mbw^T\,\Hu\,\mbu \label{eq:H(u)SVK derivation} 
\end{align}
where $\mbu$, $\mbw$ represent the finite element nodal values of the fields $\u$, $\w$ assembled across all elements.
Thus, from equation~\eqref{eq:H(u)SVK derivation}, we obtain the H-matrix for the St.~Venant-Kirchhoff model as
\begin{align}
    \mathbf{H(u)} & = \bigcup\limits_{e=1}^{n_e} \int_{\Omega_e} \frac{1}{2}\sum_A \sum_B \Bigg( \lambda \brac{\gNa \otimes \gNb} \brac{2\I + \D_e^T} + 2\mu \brac{\gNb \cdot \gNa} \nonumber\\ 
        & \qquad + \mu \brac{\gNb \otimes \gNa} \brac{2\I + \D_e^T} + \mu \D_e^T \brac{\gNa \cdot \gNb}  + 2 \brac{\gNb \cdot \S_{\text{SVK},e} \cdot \gNa}\Bigg) \,dV_e \,. \nonumber
\end{align}
\subsection{Yeoh model}\label{subsec:appendix:Yeoh_derivation}
In this section, we outline the development of the H-matrix specific to the Yeoh model, employing an approach analogous to the one utilized in the construction of the H-matrix for the St.~Venant-Kirchhoff model. We build on the framework developed for the St.~Venant-Kirchhoff model and to this end, we continue to employ equations~\eqref{eq: Green-Lagrange strain and deformation gradient}, \eqref{eq:internal virtual work - SVK}, \eqref{eq: u,w,X discretization}, \eqref{eq:Fe in FE space}, \eqref{eq:gradX in FE space}, \eqref{eq: D_e} and also persist using the notation from the previous exercise~(\ref{subsec:appendix:SVK_derivation}). \\
The strain energy function for the Yeoh model is given as follows
\begin{equation}
    \Psi(\E) = C_{10}\brac{\Ibar - 3} + C_{20} \brac{\Ibar - 3}^2 + C_{30} \brac{\Ibar - 3}^3 + \frac{1}{D_1}\brac{J-1}^2
\end{equation}
where $J = \det(\F)$, $I_1 =\text{tr}(2\E + \I)$, $\Ibar = I_1 \Jtt$; $C_{10}$, $C_{20}$, $C_{30}$, ${D_1}$ represent the material constants of the model. 
The corresponding second Piola-Kirchoff stress tensor is thus given as
\begin{align}\label{eq: second Piola-Kirchhoff stress tensor - Yeoh}
    \S & = \frac{\partial{\Psi(\E)}}{\partial{\E}} \nonumber\\
    & = 2\,\Jtt\Bigg[C_{10} + 2\,C_{20}\,\brac{\trC\,\Jtt - 3} + 3 \, C_{30}\,\brac{\trC\,\Jtt - 3}^2\Bigg]\brac{\I-\frac{\C^{-1}\trC}{3}} \nonumber\\
    & \qquad + 2\,\brac{\frac{1}{D_1}}\,J\,(J-1)\,\C^{-1}
\end{align}
where $\C = 2\E + \I$.
Using equations~\eqref{eq: Green-Lagrange strain and deformation gradient}, \eqref{eq:internal virtual work - SVK}, \eqref{eq: second Piola-Kirchhoff stress tensor - Yeoh}, the internal virtual work for the Yeoh model is then given as
%
\begin{equation}\label{eq: internal virtual work discretized - Yeoh}
\begin{split}
    \int_{\Omega_0} \delta\E \colon \S dV  & =  \bigcup\limits_{e=1}^{n_e} \int_{\Omega_e} \frac{1}{2} \brac{\F_e^T\brac{\mathbf{\grad_X}\wc_e} + \brac{\mathbf{\grad_X}\wc_e}^T\F_e} \colon \Bigg[ 2\,\detFtt\Bigg(C_{10} + 2\,C_{20}\,\brac{\trFTF\,\detFtt - 3} \\ 
    & \qquad + 3 \, C_{30}\,\brac{\trFTF\,\detFtt - 3}^2\Bigg)\brac{\I-\frac{\Cinv\trFTF}{3}} \\
    & \qquad + 2\,\brac{\frac{1}{D_1}}\,\detF\,(\detF-1)\,\Cinv \Bigg] \, dV_e\,. \\
    \end{split}
\end{equation}
Now, let us define the following for ease of representation:
\begin{align}
    \alpha_e & = 2\,\detFtt\,\Bigg(C_{10} + 2\,C_{20}\,\brac{\trFTF\,\detFtt - 3} \nonumber \\
    & \qquad + 3 \, C_{30}\,\brac{\trFTF\,\detFtt - 3}^2\Bigg) \, \label{eq: alpha_e}\\
      \beta_e & = 2\,\brac{\frac{1}{D_1}}\,\detF\,. \label{eq: beta_e}
\end{align}
Thus, equation~\eqref{eq: internal virtual work discretized - Yeoh} can be rewritten using equations~\eqref{eq: alpha_e} and \eqref{eq: beta_e} as follows
\begin{equation}\label{eq: internal virtual work discretized - Yeoh - alpha_e and beta_e}
\begin{split}
    \int_{\Omega_0} \delta\E \colon \S dV  & =  \bigcup\limits_{e=1}^{n_e} \int_{\Omega_e} \frac{1}{2} \brac{\F_e^T\brac{\mathbf{\grad_X}\wc_e} + \brac{\mathbf{\grad_X}\wc_e}^T\F_e} \colon \Bigg[ \alpha_e\,\brac{\I-\frac{\Cinv\trFTF}{3}} \\
    & + \beta_e\,(\detF-1)\,\Cinv \Bigg] \, dV_e\, \\
    \end{split}
\end{equation}
Next, consider a matrix $\P_e$ defined as follows:
\begin{equation}\label{eq: P matrix definition}
    \P_e = \left[\begin{array}{c|c|c} \begin{matrix}1 + \dudxp{2}{2}{3}{3} \\ + \dudxm{2}{2}{3}{3} - \dudxm{2}{3}{3}{2}\end{matrix} & \begin{matrix} -\frac{1}{2}\dudx{2}{1} \\ +\dudxm{2}{3}{3}{1} - \dudxm{2}{1}{3}{3} \end{matrix} & \begin{matrix} -\frac{1}{2}\dudx{3}{1} \\ +\dudxm{2}{1}{3}{2} - \dudxm{3}{1}{2}{2} \end{matrix} \\\hline \begin{matrix} -\frac{1}{2}\dudx{1}{2} \\ +\dudxm{1}{3}{3}{2} - \dudxm{1}{2}{3}{3} \end{matrix} & \begin{matrix}1 + \dudxp{1}{1}{3}{3} \\ + \dudxm{1}{1}{3}{3} - \dudxm{1}{3}{3}{1}\end{matrix} & \begin{matrix} -\frac{1}{2}\dudx{3}{2} \\ +\dudxm{1}{2}{3}{1} - \dudxm{1}{1}{3}{2} \end{matrix} \\\hline \begin{matrix} -\frac{1}{2}\dudx{1}{3} \\ +\dudxm{1}{2}{2}{3} - \dudxm{1}{3}{2}{2} \end{matrix} & \begin{matrix} -\frac{1}{2}\dudx{2}{3} \\ +\dudxm{1}{3}{2}{1} - \dudxm{1}{1}{2}{3} \end{matrix} & \begin{matrix}1 + \dudxp{1}{1}{2}{2} \\ + \dudxm{1}{1}{2}{2} - \dudxm{1}{2}{2}{1}\end{matrix} \end{array}\right]\,.
\end{equation}
We then have the following result:
\begin{align}\label{eq: detF - 1}
    \detF - 1 = \text{tr}\brac{\begin{bmatrix} \dudx{1}{1} & \dudx{1}{2} & \dudx{1}{3}\\ \dudx{2}{1} & \dudx{2}{2} & \dudx{2}{3} \\ \dudx{3}{1} & \dudx{3}{2} & \dudx{3}{3}\end{bmatrix}^T \P_e} = \sum_A \gNa \cdot \P_e^T \cdot \uAe
\end{align}
Substituting equations~\eqref{eq:Fe in FE space}, \eqref{eq:gradX in FE space}, \eqref{eq: detF - 1} in equation~\eqref{eq: internal virtual work discretized - Yeoh - alpha_e and beta_e}, we get
\begin{equation}\label{eq:H(u)Yeoh derivation}
\begin{split}
    \int_{\Omega_0} \delta\E \colon \S dV  & :=  \bigcup\limits_{e=1}^{n_e} \int_{\Omega_e}\Bigg[\sum_A \sum_B \wAe\Bigg( 2\alpha_e\brac{\gNb \cdot \Finve \cdot \gNa}\I + 2\alpha_e\brac{\FinvTe\cdot\brac{\gNb \otimes \gNa}} \\
    & \qquad - \alpha_e\, \brac{\gNb \cdot \brac{\Finve \FinvTe \mathbf{D}_e^{T}}\cdot\gNa}\I - \alpha_e\,\mathbf{D}_e \,\Finve \FinvTe \brac{\gNb \otimes \gNa} \\
    & \qquad - \frac{2}{3} \alpha_e\,\Finve \FinvTe \brac{\gNa \otimes \gNb}\brac{2\mathbf{I} + \mathbf{D}_e^T} \\
    & \qquad + 2\alpha_e\,\brac{\gNb \cdot \left[ \I - \frac{\Finve \FinvTe \text{tr}(\mathbf{F}^T_e\mathbf{F}_e)}{3} \right] \cdot \gNa}\I + \brac{2 \beta_e} \Finve \FinvTe \brac{\gNa \otimes \gNb}\mathbf{P}_e^T \\ & \qquad + \brac{2\beta_e}\brac{\detF-1} \brac{\gNa \cdot \brac{\Finve \FinvTe} \cdot \gNb}\mathbf{I}\Bigg)\uBe \Bigg] \,dV_e \\
    & =  \mbw^T \Bigg[\bigcup\limits_{e=1}^{n_e} \int_{\Omega_e}\sum_A \sum_B \Bigg( 2\alpha_e\brac{\gNb \cdot \Finve \cdot \gNa}\I + 2\alpha_e\brac{\FinvTe\cdot\brac{\gNb \otimes \gNa}} \\
    & \qquad - \alpha_e\, \brac{\gNb \cdot \brac{\Finve \FinvTe \mathbf{D}_e^{T}}\cdot\gNa}\I - \alpha_e\,\mathbf{D}_e \,\Finve \FinvTe \brac{\gNb \otimes \gNa} \\
    & \qquad - \frac{2}{3} \alpha_e\,\Finve \FinvTe \brac{\gNa \otimes \gNb}\brac{2\mathbf{I} + \mathbf{D}_e^T} \\
    & \qquad + 2\alpha_e\,\brac{\gNb \cdot \left[ \I - \frac{\Finve \FinvTe \text{tr}(\mathbf{F}^T_e\mathbf{F}_e)}{3} \right]\cdot \gNa}\I + \brac{2 \beta_e} \Finve \FinvTe \brac{\gNa \otimes \gNb}\mathbf{P}_e^T \\ & \qquad + \brac{2\beta_e}\brac{\detF-1} \brac{\gNa \cdot \brac{\Finve \FinvTe} \cdot \gNb}\mathbf{I}\Bigg) \,dV_e \Bigg]\mbu \\
    & = \mbw^T \, \Hu \,\mbu \,.
    \end{split}
\end{equation}
where $\mbu$,$\mbw$ represent the finite element nodal values of the fields $\u$, $\w$ assembled across all elements.
Thus, from equation~\eqref{eq:H(u)Yeoh derivation}, we obtain the H-matrix for the Yeoh model as
\begin{align}
    \mathbf{H(u)} & = \bigcup\limits_{e=1}^{n_e} \int_{\Omega_e}\sum_A \sum_B \Bigg( 2\alpha_e\brac{\gNb \cdot \Finve \cdot \gNa}\I + 2\alpha_e\brac{\FinvTe\cdot\brac{\gNb \otimes \gNa}} \nonumber\\
    & \qquad - \alpha_e\, \brac{\gNb \cdot \brac{\Finve \FinvTe \mathbf{D}_e^{T}}\cdot\gNa}\I - \alpha_e\,\mathbf{D}_e \,\Finve \FinvTe \brac{\gNb \otimes \gNa} \nonumber\\
    & \qquad - \frac{2}{3} \alpha_e\,\Finve \FinvTe \brac{\gNa \otimes \gNb}\brac{2\mathbf{I} + \mathbf{D}_e^T} \nonumber\\
    & \qquad + 2\alpha_e\,\brac{\gNb \cdot \left[ \I - \frac{\Finve \FinvTe \text{tr}(\mathbf{F}^T_e\mathbf{F}_e)}{3} \right]\cdot \gNa}\I + \brac{2 \beta_e} \Finve \FinvTe \brac{\gNa \otimes \gNb}\mathbf{P}_e^T \nonumber \\ & \qquad + \brac{2\beta_e}\brac{\detF-1} \brac{\gNa \cdot \brac{\Finve \FinvTe} \cdot \gNb}\mathbf{I}\Bigg) \,dV_e \,. \nonumber
\end{align}

{
 }

 \bibliographystyle{elsarticle-num} 


\end{document}